\documentclass[12pt]{article} 

\usepackage[applemac]{inputenc}
\usepackage{amssymb}
\usepackage{amstext} 
\usepackage{amsmath} 
\usepackage{amscd}
\usepackage{hyperref}
\usepackage[all]{xypic}
\usepackage{rotating}
\usepackage{color}
\usepackage{url} 
\usepackage[english,francais]{babel}
\usepackage{theorem}
\usepackage{mathrsfs}

\newcommand{\modif}[1]{{#1}}
\newcommand{\modiff}[1]{#1} 
\newlength{\secskip}
\newlength{\thsecskip}
\newlength{\ssecskip}
\newlength{\thssecskip}
\newlength{\sssecskip}
\newlength{\thsssecskip}
\makeatletter
\renewcommand\section{\@startsection{section}{1}{\z@}%
                                   {\secskip}
                                   {2.3ex \@plus.2ex}
                                   {\normalfont\Large\bfseries}}
\renewcommand\subsection{
			  \@startsection{subsection}{2}{\z@}%
                                     {\ssecskip}
                                     {1.5ex \@plus .2ex}
									 {\normalfont\normalsize\bfseries}}
\renewcommand\subsubsection{
              \@startsection{subsubsection}{3}{\z@}%
                                     {\sssecskip}
                                     {\z@}
									 {\normalfont\normalsize\bfseries}}
\makeatother
\newcommand\Subsection[1]{\subsection{{\hskip-0.5em}#1}}
\makeatletter
\renewcommand\l@subsection{\@dottedtocline{2}{1.5em}{2.7em}}
\makeatother
\newcommand\Subsubsection[1]{\subsubsection{{\hskip-0.5em}#1}}
\makeatletter
\renewcommand\l@subsubsection{\@dottedtocline{3}{3.8em}{3.5em}}
\makeatother
\newcommand\rref[1]{{\rm\ref{#1}}}
\theoremstyle{change}
\theorembodyfont{\slshape}
\newtheorem{fteostar}{Th\'eor\`eme. --}
\newtheorem{fpropstar}[fteostar]{Proposition. --}
\newtheorem{flemstar}[fteostar]{Lemme. --}

\newtheorem{fteo}[subsection]{Th\'eor\`eme. --}
\newtheorem{fsteo}[subsubsection]{Th\'eor\`eme. --}

\newtheorem{fsprop}[subsubsection]{Proposition. --}
\newtheorem{fssprop}[paragraph]{Proposition. --}
\newtheorem{flem}[subsection]{Lemme. --}
\newtheorem{fslem}[subsubsection]{Lemme. --}
\newtheorem{fsslem}[paragraph]{Lemme. --}

\newtheorem{fscor}[subsubsection]{Corollaire. --}

\newtheorem{fdefi}[subsection]{D\'efinition. --}
\newtheorem{fsdefi}[subsubsection]{D\'efinition. --}

\theorembodyfont{\rmfamily}

\newtheorem{fsrem}[subsubsection]{Remarque. --}
\newtheorem{fssrem}[paragraph]{Remarque. --}

\newtheorem{fsrems}[subsubsection]{Remarques. --}
\newtheorem{fssrems}[paragraph]{Remarques. --}
\newtheorem{fcond}[subsection]{Condition}
\newtheorem{fscond}[subsubsection]{Condition}
\newtheorem{fsscond}[paragraph]{Condition}

\newcommand{\Aut}{{\rm Aut}}

\newcommand{\Spec}{\mathop{\mathrm{Spec}}}
\newcommand{\Card}{\mathrm{Card}}

\newcommand{\dem}{{\sl \noindent D\'emonstration\hskip0.1em: }} 
\newcommand{\rH}{\mathrm{H}}

\newcommand{\GL}{{\rm GL}}
\newcommand{\SL}{{\rm SL}}

\newcommand{\id}{{\mathrm{Id}}}

\newcommand{\Hilb}{{\rm Hilb}}
\newcommand{\Gal}{{\rm Gal}}

\newcommand{\W}{\mathop{\textstyle\prod}}
\newcommand{\PC}[1]{\mathop{{\scriptstyle\wedge}}\limits^{#1}}

\newcommand{\im}{\mathrm{Im}}

\newcommand{\Top}{\mathrm{top}}

\newcommand{\FF}{{\mathbb F}}
\newcommand{\CC}{{\mathbb C}}
\newcommand{\RR}{{\mathbb R}}

\newcommand{\QQ}{{\mathbb Q}}
\newcommand{\NN}{{\mathbb N}}

\newcommand{\PP}{{\mathbb P}}

\newcommand{\GG}{{\mathbb G}}

\newcommand{\AAA}{{\mathbb A}}

\newcommand{\VAR}{\mathrm{VAR}}

\newcommand{\TOP}{\mathrm{TOP}}
\newcommand{\ENS}{\mathrm{ENS}}

\newcommand{\ul}[1]{\underline{#1}}
\newcommand{\Aa}{\mathbb{A}}

\renewcommand{\setminus}{\smallsetminus}

\def\FG{\mathfrak{G}}

\def\FI{\mathfrak{I}}

\def\Fm{\mathfrak{m}}

\newcommand{\cC}{\mathscr{C}} 

\newcommand{\cE}{\mathscr{E}}

\newcommand{\cX}{\mathscr{X}}
\newcommand{\cY}{\mathscr{Y}}
\newcommand{\cZ}{\mathscr{Z}}

\newcommand{\calL}{\mathscr{L}}
\newcommand{\cA}{\mathscr{A}}

\newcommand{\cO}{\mathscr{O}}

\newcommand{\ol}{\overline}
\newcommand{\til}{\widetilde}

\newcommand{\OOO}{\mathscr{O}}

\newcommand{\bec}{{\natural}}

\newcommand{\fl}{\rightarrow} 
\newcommand{\ffl}{\longrightarrow} 
\newcommand{\inj}{\hookrightarrow}

\newcommand{\sfl}[1]{\mathop{\fl}\limits^{#1}}
\newcommand{\sffl}[1]{\mathop{\ffl}\limits^{#1}}
\newcommand{\flis}{\sfl{\sim}} 
\newcommand{\fflis}{\sffl{\sim}}

\newcommand{\surj}{\rlap{$\rightarrow$}\kern-2pt\rightarrow}
\newcommand{\ssurj}{\rlap{$\longrightarrow$}\kern-2pt\longrightarrow}
\newcommand{\surjgauche}{\rlap{$\leftarrow$}\kern-2pt\leftarrow}
\newcommand{\ssurjgauche}%
{\rlap{$\longleftarrow$}\kern-2pt\longleftarrow}

\newcommand{\mapdown}[1]
{\big\downarrow\rlap{$\vcenter{\hbox{$\scriptstyle#1$}}$}}

\newcommand{\mapup}[1]
{\big\uparrow\rlap{$\vcenter{\hbox{$\scriptstyle#1$}}$}}

\newcommand{\Mapdown}[1]
{\Big\downarrow\rlap{$\vcenter{\hbox{$\scriptstyle#1$}}$}}

\newcommand{\varfl}[1]{%
\setbox0=\hbox{$\;\;{\scriptstyle#1}\;\;\;$}%
\setbox1=\hbox to\wd0{$\;$\rightarrowfill$\;$}%
{\mathop{\box1}\limits^{\box0}}%
}
\newcommand{\varflspecial}[2]{%
\setbox0=\hbox{$\;\;{\scriptstyle#1}\;\;\;$}%
\setbox1=\hbox to\wd0{$\;$\rightarrowfill$\;$}
{\mathop{\box1}\limits^{#2}}%
}
\newcommand{\varflman}[2]{
\setbox0= \hbox to #1{$\;$\rightarrowfill$\;$}
{\mathop{\box0}\limits^{#2}}
}
\newcommand{\varmapstoman}[2]{%
\setbox0= \hbox to #1{$\;\mapstochar$\rightarrowfill$\;$}
{\mathop{\box0}\limits^{#2}}
}

\newcommand{\pr}{\mathrm{pr}}

\newcounter{nc}
\renewcommand{\thenc}{{\rm(\roman{nc})}}
\newenvironment{romlist}%
{\begin{list}{\thenc}{
\usecounter{nc} 
\parsep=0pt
\setlength  \labelwidth{\leftmargin}
\addtolength\labelwidth{-\labelsep}
}
}{\end{list}}
\newcounter{nnc}
\renewcommand{\thennc}{{\rm(\alph{nnc})}}
\newenvironment{subromlist}%
{\begin{list}{\thennc}{
\usecounter{nnc}
\parsep=0pt
\setlength  \labelwidth{\leftmargin}
\addtolength\labelwidth{-\labelsep}
}
}{\end{list}} 
\newcommand{\pauseromlist}%
{\global\edef\savecount{\arabic{nc}}\end{romlist}}
\newcommand{\finpauseromlist}%
{\begin{romlist}\setcounter{nc}{\savecount}}

\newcounter{ctnum}
\renewcommand{\thectnum}{\textup{(\arabic{ctnum})}}
\newenvironment{numlist}%
{\begin{list}{\thectnum}{
\usecounter{ctnum} 
\parsep=0pt
\leftmargin=0pt%
\setlength{\itemindent}{\labelwidth}%
\addtolength{\itemindent}{\labelsep}%
}
}{\end{list}}
\newcommand{\cqfd}{ \phantom{$\square$}\hfill\llap{$\square$}}

\newcommand{\red}{{\mathrm{red}}}
\newcommand{\wh}[1]{\widehat{#1}}

\def\lgr{\longrightarrow}
\def\la{\longleftarrow}
\def\simlgr{\buildrel\sim\over\lgr}
\def\simla{\buildrel\sim\over\la}

\title{ Fibr\'es principaux sur les corps valu\'es hens\'eliens}
\author{Ofer Gabber\thanks{C.N.R.S. et I.H.\'E.S., Le Bois-Marie, 35 route de Chartres,
F-91440 Bures sur Yvette. \href{mailto:gabber@ihes.fr}{gabber@ihes.fr}}\\
Philippe  Gille\thanks{UMR 5208 du CNRS -
Institut Camille Jordan - Universit\'e Claude Bernard Lyon 1,
43 boulevard du 11 novembre 1918,
F-69622 Villeurbanne cedex. \href{mailto:gille@math.univ-lyon1.fr}{gille@math.univ-lyon1.fr}}\, \footnote{L'auteur a b\'en\'efici\'e du soutien du projet ANR Gatho,  ANR-12-BS01-0005.} 
\\
Laurent  Moret-Bailly\thanks{IRMAR, Universit\'{e} de Rennes 1, Campus de Beaulieu, F-35042 Rennes Cedex.  \href{mailto:laurent.moret-bailly@univ-rennes1.fr}{laurent.moret-bailly@univ-rennes1.fr}}\ %
}

\begin{document}
\date{}
\maketitle
\begin{center}\slshape
Accept\'e pour publication dans \textup{Algebraic Geometry} le 4 juin 2014.

\end{center}
\noindent{\bf R\'esum\'e.} Soit $K$ le corps des fractions d'un anneau de valuation hens\'elien $A$. 
On suppose que le compl\'et\'e $\wh{K}$ est une extension s\'eparable de $K$. Soient $Y$ une $K$-vari\'et\'e, $G$ un $K$-groupe
alg\'ebrique et $f:X\to Y$ un $G$-torseur au-dessus de $Y$. 
On consid\`ere l'application induite $X(K)\to Y(K)$, continue pour les topologies d\'eduites de la valuation.
Si $I$ d\'esigne son image, nous montrons que $I$ est localement ferm\'ee dans $Y(K)$; de plus la surjection induite
$X(K)\to I$ est un fibr\'e principal sous le groupe topologique $G(K)$.
\medskip\selectlanguage{english}

\noindent{\bf Abstract:} Let $K$  be the fraction field of a henselian valuation ring $A$. Assume that 
the completion $\wh{K}$ is a separable extension of $K$. Let $Y$ be a $K$-variety, $G$ an algebraic group over $K$, and 
$f:X \to Y$ a $G$-torsor over $Y$. We consider the induced map $X(K) \to Y(K)$, 
which is  continuous for the topologies deduced from the valuation. If $I$ 
denotes the image of this map, we prove that $I$ is locally closed in $Y(K)$; moreover, the induced surjection $X(K)\to I$ 
is a principal bundle with group $G(K)$ (also topologized by the valuation).\medskip


\medskip


\setcounter{tocdepth}{2}
\selectlanguage{francais}
\tableofcontents


\bigskip

\addtocounter{section}{-1}

\section{Introduction}
\Subsection{Notations}\label{NotIntro}
Soit  $K$ le corps des fractions d'un anneau de valuation $A$; on notera $v$ la
valuation associ\'ee et $\Gamma$ son groupe. La donn\'ee de $v$ d\'etermine une structure 
de corps topologique s\'epar\'e sur $K$; nous supposerons toujours qu'il n'est pas discret, 
c'est-\`a-dire que $\Gamma\neq0$. On notera  $\wh{K}$ le compl\'et\'e de $K$.

Dans cette introduction, nous supposerons $(K,v)$ \emph{admissible}, au sens suivant:
\begin{fsdefi}\label{def:admissible} Avec les notations ci-dessus, on dit que $(K,v)$ (ou $A$) est 
\emph{admissible} si $A$ est \emph{hens\'elien} et si l'extension $\wh{K}/K$ est \emph{s\'eparable}.
\end{fsdefi}

Pour toute  $K$-\emph{vari\'et\'e} (c'est-\`a-dire tout $K$-sch\'ema de type fini) $X$, on note $X_\Top$ 
l'ensemble $X(K)$ muni de la topologie d\'eduite de  la topologie de $K$. Tout $K$-morphisme $f: X\to Y$ de 
$K$-vari\'et\'es induit une application continue $f_{\Top}:X_{\Top}\to Y_{\Top}$. 

Nous nous int\'eressons dans cet article au cas o\`u un $K$-groupe alg\'ebrique $G$ agit \`a droite  sur $X$ et 
o\`u $f:X\to Y$ est un $G$-torseur pour cette action. Noter qu'alors le groupe topologique $G_{\Top}$ agit librement 
et contin\^ument sur $X_{\Top}$ et que l'on a une \emph{bijection continue} $\ol{f_{\Top}}:X_{\Top}/G(K)\ffl I:=\im(f_{\Top})$. 
L'objet de ce travail est l'\'etude topologique des applications  $f_{\Top}$ et $\ol{f_{\Top}}$.\smallskip

Il est en fait naturel (et \`a certains \'egards plus simple) de consid\'erer la situation plus g\'en\'erale o\`u $X$ et $Y$ sont 
des \emph{espaces alg\'ebriques de type fini} (et quasi-s\'epar\'es) sur $K$, que nous appellerons dor\'enavant 
\og $K$-espaces\fg. Tout d'abord, un $G$-torseur sur une $K$-vari\'et\'e $Y$ --- d\'efini, comme il se doit, 
comme faisceau fppf sur $Y$ --- n'est pas n\'ecessairement un sch\'ema (en-dehors du cas important o\`u $G$ est affine),
alors que c'est automatiquement un $K$-espace, selon un th\'eor\`eme d'Artin; \`a partir de l\`a, il devient \'egalement 
judicieux d'envisager des $G$-torseurs $X\to Y$ o\`u $X$ \emph{et $Y$} sont des $K$-espaces, afin de travailler dans une 
cat\'egorie stable par les op\'erations usuelles, produits fibr\'es notamment.

D'autre part, si un $K$-groupe alg\'ebrique $G$ op\`ere  librement \`a droite sur une $K$-vari\'et\'e $X$, le faisceau 
fppf quotient $X/G$ est toujours un $K$-espace (voir \ref{QuotLibre}) mais, de nouveau, n'est pas n\'ecessairement une $K$-vari\'et\'e. 

Ces consid\'erations nous ont conduits \`a formuler syst\'ematiquement nos r\'esultats dans le cadre des $K$-espaces; si
le plan g\'en\'eral des d\'emonstrations n'en est pas affect\'e, il nous a fallu revenir sur des r\'esultats bien connus 
pour les sch\'emas mais dont l'extension aux espaces alg\'ebriques n'est pas suffisamment document\'ee: voir par exemple \ref{sec:espalg}
pour la d\'efinition et les propri\'et\'es de $X_{\Top}$ lorsque $X$ est un $K$-espace. (En ce qui concerne $G$, rappelons qu'un 
espace alg\'ebrique en groupes \modif{quasi-s\'epar\'e} de type fini sur un corps est toujours un sch\'ema \cite[Lemma 4.2]{A1}).

\smallskip

Notre r\'esultat principal est le suivant:

\begin{fteo}\label{ThPpal} Soient $(K,v)$ un corps valu\'e admissible, $G$ un $K$-groupe 
alg\'ebrique \emph{(c'est-\`a-dire un $K$-sch\'ema en groupes de type fini)}, $Y$ un $K$-espace, 
$f:X\to Y$ un $G$-torseur au-dessus de $Y$. D\'efinissons $f_{\Top}:X_{\Top}\to Y_{\Top}$ et
$\ol{f_{\Top}}:X_{\Top}/G(K)\ffl I:=\im(f_{\Top})$ 
comme ci-dessus. Alors:
\begin{numlist} 

\item\label{ThPpal1} en tant que sous-espace de $Y_{\Top}$,  $I$ est:
	\begin{subromlist}
	\item\label{ThPpal1a} localement ferm\'e (dans tous les cas);
	\item\label{ThPpal1b} ouvert et ferm\'e si $G$ est lisse ou si $K$ est parfait;
        \item \label{ThPpal1c} ferm\'e si $G$ satisfait la condition $(*)$ \textup{(voir  \ref{condition:star})}, et en
        particulier si  $G^\circ_{\red}$ est lisse, ou  si $G^\circ$ est commutatif, ou si $G$ est de rang r\'eductif nul 
        (i.e. si $G_{\ol{K}}$ n'a pas de sous-tore non trivial);
\end{subromlist}
\item\label{ThPpal2} L'application ${f_{\Top}}$ est ouverte sur son image $I$; en particulier, la bijection $\ol{f_{\Top}}$ est 
un hom\'eomorphisme.
\item\label{ThPpal3} Si $Y$ est \emph{localement s\'epar\'e (par exemple une vari\'et\'e, cf. \ref{conventions})}, 
alors ${f_{\Top}}$ fait de $X_{\Top}$ un $G_{\Top}$-fibr\'e principal au-dessus de $I$.
\end{numlist}
\end{fteo} 

\Subsection{Plan de la d\'emonstration. }\label{sous-sect-plan} 
Pour ne pas alourdir l'introduction, nous supposons ici que $Y$ est localement s\'epar\'e, par exemple une $K$-vari\'et\'e.

Le cas, sans doute bien connu, o\`u le groupe $G$ est lisse est trait\'e au \S\,\ref{sec:CorpsTop}. Dans le cas g\'en\'eral, on
 note $G^\bec$ le plus grand $K$-sous-groupe lisse   de 
$G$ (\S\,\ref{ssec:pgsgl}). On d\'ecompose alors le $G$-torseur $f: X \to Y$ en 
\[
\xymatrix@M=1.5mm@C=16mm{
X \ar[dr]^{f} \ar[d]_{\pi} & 
\\
Z:= X/G^\bec \ar@{^{}->}[r]^-{h} & Y.
}
\] 
Le morphisme $\pi$ est un $G^\bec$-torseur et $h:Z \to Y$ est une \og fibration en $G/G^\bec$\fg. Les propri\'et\'es de $G^\bec$ 
impliquent que l'application $h_K: Z(K) \to Y(K)$  est \emph{injective}. D'un point de vue topologique, on a donc le diagramme 
\[
\xymatrix@M=1.5mm@C=16mm{
X_\Top \ar[dr]^{f_\Top } \ar[d]_{\pi_\Top } & 
\\
Z_\Top  \ar@{^{(}->}[r]^-{h_\Top } & Y_\Top  .
}
\] 
D'apr\`es le cas lisse (\S\,\ref{tors-lisse}), $\pi_\Top$ fait de $X_{\Top}$ un $G^\bec_\Top$-fibr\'e principal sur l'image de
$\pi_\Top$, qui est ouverte et ferm\'ee dans $Z_{\Top}$. 

L'application $h_\Top$ 
est nettement plus d\'elicate \`a analyser. Les deux outils cl\'es pour cette \'etude sont:
\begin{itemize}
\item le th\'eor\`eme d'approximation de Greenberg (g\'en\'eralis\'e aux corps valu\'es admissibles dans \cite{MB2});
\item la \og bonne compactification\fg\ (\S\,\ref{section-compactification}) de $G/G^\bec$, et la compactification
relative de $h$ qui s'en d\'eduit.
\end{itemize}
On montre ainsi que  $h_\Top$ est un hom\'eomorphisme sur son image, et que celle-ci 
s'\'ecrit $F_1 \setminus F_2$ pour des ferm\'es remarquables $F_1$, $F_2$ de  $Y_\Top$ (lemme \ref{clef}).

\Subsection{Application aux espaces homog\`enes. } Un cas particulier important est celui o\`u $X=H$ est un $K$-groupe 
al\-g\'e\-brique contenant $G$ comme sous-groupe, et o\`u $f:H\to Y:=H/G$ est le morphisme de passage au quotient. L'image 
$I$ est alors l'orbite sous $H(K)$ de la classe neutre $y_{0}\in Y(K)$ dans l'ensemble $Y(K)=(H/G)(K)$: on voit donc par \ref{ThPpal} 
qu'elle est localement ferm\'ee dans $Y_{\Top}$ et qu'elle s'identifie (avec sa topologie) \`a l'espace quotient $H_{\Top}/G(K)$ (avec, 
le cas \'ech\'eant, les compl\'ements \ref{ThPpal1}\,\ref{ThPpal1b}, 
\ref{ThPpal1}\,\ref{ThPpal1c} de l'\'enonc\'e) et que  $H_{\Top}$ est m\^eme un $G_{\Top}$-fibr\'e principal au-dessus de $I$.

Si de plus $G$ est \emph{distingu\'e} dans $H$, il en r\'esulte que $I$ est toujours ferm\'ee puisque c'est un sous-groupe localement 
ferm\'e du groupe topologique $Y_{\Top}$.

\Subsection{Application aux orbites. }Consid\'erons un $K$-groupe alg\'ebrique $H$ op\'erant \`a gauche sur une 
$K$-vari\'et\'e $S$ (ou plus g\'en\'eralement sur un $K$-espace); soit $s_{0}$ un point de $S(K)$, de stabilisateur 
$G\subset H$. On sait alors que le morphisme d'orbite $\omega_{s_{0}}:h\mapsto h.s_{0}$ se factorise en 
$$ H\buildrel {f}\over \longrightarrow H/G \fflis Y \inj S$$
o\`u $f$ est le morphisme canonique, o\`u la deuxi\`eme fl\`eche est un isomorphisme et la troisi\`eme une immersion
(voir \cite[III.3.5.2]{DG} lorsque $S$ est une vari\'et\'e, et \ref{ThOrbite} plus bas pour les $k$-espaces). D'autre part,
l'image $I$ de $H(K)$ par $\omega_{s_{0}}$ est \'evidemment l'orbite de $s_{0}$ sous l'action de $H(K)$, dans l'espace $S(K)$, et 
le stabilisateur de $s_{0}$ pour  cette action est $G(K)$. 

Le th\'eor\`eme \ref{ThPpal} nous dit donc que cette orbite est localement ferm\'ee dans $Y_{\Top}$ \emph{et donc dans
$S_{\Top}$}, et que la surjection continue canonique $H_{\Top}\surj I$ est une fibration principale de groupe $G_{\Top}$; en particulier 
elle induit un hom\'eomorphisme de $H_{\Top}/G(K)$ sur $I$.

Lorsque $K$ est un corps \emph{local}, on trouve en partie ce r\'esultat dans l'article \cite{BZ}; plus pr\'ecis\'ement, les auteurs
montrent que $I$ est localement ferm\'e dans $S_{\Top}$ et s'identifie topologiquement \`a $H_{\Top}/G(K)$  (dans \cite{BZ}, ces 
assertions sont dispers\'ees entre 1.5, 1.6 et 6.8).

\Subsection{Plan de l'article. }
Le \S\,\ref{sec:actions} est consacr\'e \`a des rappels sur les groupes alg\'ebriques: produit contract\'e d'un torseur et d'un
sch\'ema \`a op\'erateurs, propri\'et\'es du plus grand sous-groupe lisse d'un groupe alg\'ebrique $G$ et la technique de 
compactification partielle d'espaces homog\`enes via les sch\'emas de Hilbert ponctuels. 

Au \S\,\ref{sec:CorpsTop}, apr\`es des rappels sur les vari\'et\'es sur un corps topologique, on introduit la notion de
corps  \emph{topologiquement hens\'elien}, caract\'eris\'ee par la validit\'e du th\'eor\`eme des fonctions implicites; 
le cas des torseurs sous un groupe lisse (cf. \ref{ThPpal}\,\ref{ThPpal1}\,\ref{ThPpal1b}) est valable dans ce cadre, et
est \'etabli au \S\,\ref{tors-lisse}.

C'est aussi sur un corps topologiquement hens\'elien que le foncteur $X\mapsto X_{\Top}$ (des vari\'et\'es vers les espaces topologiques) 
s'\'etend aux espaces al\-g\'e\-briques avec de bonnes propri\'et\'es. 

Au \S\,\ref{sec:admissible}, nous donnons quelques propri\'et\'es des corps valu\'es admissibles, notamment 
le th\'eor\`eme d'approximation fort (g\'en\'eralisant celui de Greenberg). Nous en d\'eduisons un r\'esultat
topologique (th\'eor\`eme \ref{th:propre}) concernant les morphismes propres, qui sera essentiel dans la preuve de 
\ref{ThPpal}. Un cas particulier de \ref{th:propre} est le fait (assez facile \`a d\'emontrer directement) que si 
$K$ est admissible et si $f$ est un morphisme \emph{fini} de $K$-vari\'et\'es, l'application induite $f_{\Top}$ est ferm\'ee (\ref{fini}). 

Le \S\,\ref{section-compactification}, ind\'ependant des pr\'ec\'edents, est consacr\'e \`a la d\'emonstration d'un 
th\'eor\`eme de compactification d\^u \`a Gabber (\ref{ofer}, cas particulier d'un th\'eor\`eme annonc\'e dans \cite{Ga}):
soient $G$  un groupe alg\'ebrique sur un corps $k$ quelconque, et $G^\bec$ son plus grand sous-groupe lisse. Alors le 
$G$-espace homog\`ene $G/G^\bec$ admet une compactification \'equivariante dont le seul point \`a corps r\'esiduel s\'eparable sur $k$ est l'origine. 

Le \S\,\ref{caslocal} contient la d\'emonstration,  esquiss\'ee plus haut, des parties restantes \ref{ThPpal1}\,\ref{ThPpal1a} et \ref{ThPpal2} 
du th\'eor\`eme \ref{ThPpal}. 
Enfin le \S\,\ref{sec:exemples} contient divers (contre-)exemples et compl\'ements.

\begin{fsrem}
Des questions voisines (dans le cas d'une valuation compl\`ete de rang $1$) sont abord\'ees dans \cite{BTh}; il est \`a noter cependant
que l'assertion (1)\,(a) du th\'eor\`eme 1 de cet article est contredite par l'exemple donn\'e au \S\,\ref{ex:OrbiteNonFermee} du pr\'esent travail.
\end{fsrem}

\Subsection{Conventions}\label{conventions}
Si $k$ est un corps, une \emph{$k$-vari\'et\'e} est un $k$-sch\'ema de type fini; un \emph{$k$-groupe alg\'ebrique} est un $k$-sch\'ema 
en groupes de type fini; un \emph{$k$-espace} est un $k$-espace alg\'ebrique de type fini et quasi-s\'epar\'e. On renvoie \`a \cite{Kn}
 pour les propri\'et\'es des espaces alg\'ebriques. Rappelons simplement ici qu'un $k$-espace $X$ est \emph{localement s\'epar\'e} si le
 monomorphisme diagonal $X\to X\times_{k}X$ est une immersion; cette condition est toujours v\'erifi\'ee si $X$ est une vari\'et\'e.

Si $k'$ est une extension finie de $k$, et $V$ un $k'$-espace, nous noterons $\W_{k'/k}(V)$ sa restriction de Weil \`a $k$. On rappelle
 qu'en tant que foncteur sur la cat\'egorie des $k$-espaces $T$, elle est d\'efinie par
$$\bigl(\W_{k'/k}(V)\bigr)(T):= V(T\times_{k}k')$$
et qu'elle est repr\'esentable par un $k$-espace \modiff{(\cite[Lemma 5.10]{ConAdelic} ou \cite[Theorem 1.5]{Ol})}, et par une $k$-vari\'et\'e lorsque $V$  est une 
$k$-vari\'et\'e quasi-pro\-jec\-tive \modiff{(\cite[Theorem 7.6/4 et Proposition 7.6/5]{BLR}, ou \cite[A.5.8]{CGP})}.

Sauf mention contraire, la cohomologie des faisceaux utilis\'ee est la cohomologie fppf. En particulier, si $G$ est un $k$-groupe 
alg\'ebrique et $Y$ un $k$-espace, $\rH^1(Y,G)$ d\'esigne l'ensemble des classes d'isomorphie de $G_{Y}$-torseurs pour la topologie fppf.

Une application continue $f:X\to Y$ entre espaces topologiques est dite \emph{stricte} si la topologie induite sur $\im(f)\subset Y$ co\"{\i}ncide avec
 la topologie quotient d\'eduite de la surjection canonique $X\to \im(f)$ (en d'autres termes, si la bijection canonique $\mathrm{Coim}(f)\to\im(f)$ est un hom\'eomorphisme).

Une application continue $f$ entre espaces topologiques est dite \emph{propre} si elle est universellement ferm\'ee 
(conform\'ement \`a \cite[I, \S \,10, \no 1, d\'efinition 1]{BTG}, nous ne faisons pas d'hypoth\`ese de s\'eparation). 
Il revient au m\^eme de dire que $f$ est ferm\'ee \`a fibres quasi-compactes \cite[I, \S \,10, \no 2, th\'eor\`eme 1]{BTG}.
\medskip

\noindent{\bf Remerciements.}  Nous tenons \`a remercier Brian Conrad et Bertrand Lemaire pour leurs suggestions bienvenues, et le rapporteur pour sa lecture attentive du manuscrit et ses nombreuses remarques.
\section{Actions de groupes alg\'ebriques: rappels et compl\'ements}\label{sec:actions}
Dans cette section, $k$ d\'esigne un corps; on fixe une cl\^oture alg\'ebrique $\ol{k}$ de $k$, et l'on note $k_{s}\subset\ol{k}$ 
la cl\^oture s\'eparable correspondante.

\Subsection{Orbites dans un espace alg\'ebrique \`a groupe d'op\'erateurs}\label{ThOrbite} 
On fixe un $k$-groupe alg\'ebrique $G$. Soit $Y$  un $k$-faisceau fppf  muni d'une action \`a gauche de $G$. Nous dirons que
 l'action est transitive (ou que $Y$ est \emph{transitif} sous $G$) si $Y\to\Spec(k)$ est un \'epimorphisme  (de fa\c{c}on \'equivalente, 
il existe une extension finie $k'$ de $k$ telle que $Y(k')\neq\emptyset$) et si le morphisme 
$$ \begin{array}{rcl}
G\times Y& \ffl &Y\times Y\\
(g,y)&\longmapsto&(g.y,y)
\end{array}
$$
est un  \'epimorphisme.
\begin{fsteo}\label{th:orbite}
Soit $Y$ un $k$-faisceau fppf muni d'une action \emph{transitive} de $G$. On suppose qu'il existe une extension finie $k'$ de $k$ et
 un \'el\'ement $y$ de $Y(k')$ tel que le sous-faisceau de $G_{k'}$ stabilisateur de  $y$ soit un \emph{sous-sch\'ema en groupes ferm\'e} de $G_{k'}$. Alors:
\begin{numlist}
\item\label{th:orbite1} $Y$ est un $k$-sch\'ema quasi-projectif.
\item\label{th:orbite2} Soit $X$ un $k$-espace muni d'une action \`a gauche de $G$, et soit $f:Y\to X$ un morphisme $G$-\'equivariant. 
Alors $f$ se factorise de mani\`ere essentiellement unique en $Y\xrightarrow{\pi} Z\xrightarrow{j} X$, o\`u:
\begin{romlist}
\item $Z$ est un $k$-sch\'ema quasi-projectif;
\item $\pi$ est fid\`element plat;
\item $j$ est une immersion.
\end{romlist}
\end{numlist}
\end{fsteo}
\dem \ref{th:orbite1} Le $k'$-faisceau $Y\times_{k}k'$ obtenu par changement de base est clairement isomorphe \`a $G/H$, 
o\`u $H$ est le stabilisateur de $y'$. C'est donc un $k'$-sch\'ema \modiff{quasi-projectif: si $G$ est lisse et connexe ce r\'esultat est d\^u \`a Chow \cite{Ch}, et le cas g\'en\'eral s'en d\'eduit, cf. \cite[VI.2.6]{R}}. L'assertion  \ref{th:orbite1} en r\'esulte par  \cite[VIII, 7.6]{SGA1} 
(descente de sch\'ema quasi-projectifs par un morphisme fini localement libre).

Montrons \ref{th:orbite2}. L'unicit\'e est claire, puisque $Z$ est n\'ecessairement le sous-faisceau de $X$ image de $f$. Montrons d'abord que $Z$, 
ainsi d\'efini, est un $k$-sch\'e\-ma quasi-projectif: choisissons une extension finie $k'$ de $k$ et un point $y\in Y(k')$. Alors le stabilisateur
 de $z:=\pi(y)$ est aussi celui de $j(z)$ ($j$ est un monomorphisme) et est donc un sous-sch\'ema en groupes ferm\'e de $G_{k'}$: on conclut par \ref{th:orbite1}. 

Pour voir que $j$ est une immersion et que $\pi$ est fid\`element plat, on peut donc, par descente fpqc, supposer $k$ alg\'ebriquement clos, ce que
 nous ferons d\'esormais. Fixons encore un point $y\in Y(k)$, de stabilisateur $S\subset G$. Alors $Y$ s'identifie \`a $G/S$ et $Z$ \`a 
$G/H$ o\`u $H$ est le stabilisateur de $z=f(y)$; $\pi$ s'identifie \`a la projection 
$G/S \to G/H$, qui est fid\`element plate. Il reste \`a voir que $j:Z\inj X$ est une immersion, $Z$ \'etant l'orbite de $z$ dans $X$. 

Si $X$ est une vari\'et\'e, cela r\'esulte de \cite[III.3.5.2]{DG}; nous allons nous ramener \`a ce cas.

Supposons d'abord $G$ \emph{r\'eduit}, de sorte que $G/H$ l'est \'egalement. On sait \cite[II, Proposition 6.7]{Kn} 
que $X$ admet un plus grand sous-espace ouvert $U$ qui est un sch\'ema, et que $U$ est dense dans $X$. Il est clair 
que $U$ est stable par l'action du groupe $G(k)$, et donc par l'action de $G$ puisque $k$ est al\-g\'e\-bri\-que\-ment clos. 
Donc, si $z\in U(k)$,  $j$ se factorise par $U$ et l'assertion r\'esulte du cas des vari\'et\'es. Sinon, comme $G$ est r\'eduit, 
$j$ se factorise par le sous-espace r\'eduit $Y$ compl\'ementaire de $U$, qui est \'egalement stable par $G$. La den\-si\-t\'e 
implique que $\dim(Y)<\dim(X)$ et l'on conclut  par r\'ecurrence sur la dimension.

Dans le cas g\'en\'eral, posons $Z=G/H$. Comme $Z_{\red}=G_{\red}/(G_{\red}\cap H)$, le cas pr\'ec\'edent montre que 
la restriction de $j$ \`a $Z_{\red}$ est une immersion. On conclut donc par le lemme \ref{CritImmersion} qui suit.\cqfd

\begin{fsslem}\label{CritImmersion} Soit $j:Z\to X$ un monomorphisme de type fini d'espaces alg\'ebriques. On suppose que la
restriction $j_{0}$ de $j$ \`a $Z_{\red}$ est une immersion (resp. une immersion ferm\'ee). Alors $j$ est une immersion (resp. une immersion ferm\'ee). 
\end{fsslem}
\dem $j_{0}$ se d\'ecompose en une immersion ferm\'ee $Z_{\red}\to U$  suivie d'une immersion ouverte $U\to X$; comme $j$ se
factorise automatiquement par $U$, on peut remplacer $X$ par $U$, et il suffit de traiter le cas d'une immersion ferm\'ee. 
Alors $j$ est propre puisque $j_{0}$ l'est; 
\modif{d'autre part $j$ est s\'epar\'e et quasi-fini donc est quasi-affine d'apr\`es \cite[A.2]{L-MB} (ou \cite[II, 6.15]{Kn} lorsque $j$ est de pr\'e\-sen\-ta\-tion finie). 
Un morphisme propre et quasi-affine est fini, et un monomorphisme fini est une immersion ferm\'ee.} \cqfd

\begin{fssrem} Comme nous l'a signal\'e un rapporteur, on peut aussi utiliser le crit\`ere valuatif 
\cite[Corollary 2.13 p. 102]{Mo} pour d\'e\-mon\-trer \ref{CritImmersion}.
\end{fssrem}

\begin{fscor}\label{cor:orbites} Soit $X$ un $k$-espace muni d'une action \`a gauche de $G$, et soit $Y$ un sous-faisceau 
de $X$, stable par $G$ et transitif. Alors $Y$ est un sous-espace localement ferm\'e de $X$, et un $k$-sch\'ema quasi-projectif.
\end{fscor}
\dem il existe une extension finie $k'$ de $k$ et un point $y\in Y(k')$. Le stabilisateur de $y$ est le m\^eme dans $Y$et dans $X$; c'est donc un sous-sch\'ema
 en groupes ferm\'e de $G$, 
et la conclusion r\'esulte facilement de \ref{th:orbite}.\cqfd

\begin{fsdefi}\label{def:orbite}
Soit $X$ un $k$-espace muni d'une action \`a gauche de $G$. Une \emph{$k$-orbite} de $X$ sous $G$ est un sous-espace de $X$, stable par $G$ et transitif.
\end{fsdefi}

\begin{fsrems} \label{rem:orbites}  Soit $X$ un $k$-espace  muni d'une action \`a gauche de $G$. 
\begin{romlist}
\item  \label{rem:orbites2} (orbite d'un $k$-point) Fixons un point  $x \in X(k)$. Le morphisme $i:\Spec(k)\to X$ correspondant est alors une immersion ferm\'ee (v\'e\-ri\-fi\-ca\-tion laiss\'ee au lecteur) de sorte que le stabilisateur  $G_x$ est un $k$-sous-groupe ferm\'e de $G$. 
On consid\`ere le morphisme d'orbite $f_x: G \to X$, $g \mapsto g.x$.
Alors le $k$-sch\'ema quotient $G/G_x$ repr\'esente  le faisceau fppf image $T_x$ de $f_x$. Le th\'eor\`eme \ref{th:orbite} montre que $T_x \to X$
 est une immersion, ainsi $T_x$ est une $k$-orbite de $X$ sous $G$.

\item   \label{rem:orbites3} Si $X$ est une $k$-vari\'et\'e, la notion de $k$-orbite sous $G$ co\"{\i}ncide avec celle de \cite[\S 10.2, d\'efinition 10.4]{BLR}: ceci r\'esulte du corollaire \ref{cor:orbites}. 
\end{romlist}

\end{fsrems}

\Subsection{Quotient par une action libre}
Soit $G$ un $k$-groupe alg\'ebrique op\'erant \`a droite  sur un $k$-espace $X$. 
On note $\alpha:X\times_{k}G\to X$ l'action en question. On suppose  que l'action $\alpha$ est \emph{libre}, c'est-a-dire que le morphisme 
$$\begin{array}{rcl}
\rho=(\pr_{1},\alpha):\;X\times_{k}G&\ffl&X\times_{k}X\\
(x,g)&\longmapsto& (x,x.g)
\end{array}
$$
est un monomorphisme. 
\label{QuotLibre}\begin{fteostar}
Sous les hypoth\`eses ci-dessus, le faisceau quotient $X/G$ est un $k$-espace, et la projection naturelle $\pi:X\to X/G$ fait de $X$ un $G$-torseur au-dessus de $X/G$.

De plus, $X/G$ est localement s\'epar\'e (resp. s\'epar\'e) si et seulement si $\rho:\;X\times_{k}G\to X\times_{k}X$ est une immersion (resp. une immersion ferm\'ee).
\end{fteostar}
\dem par d\'efinition, $X/G$ est le quotient (au sens des faisceaux fppf) de $X$ par la relation d'\'equivalence $R\subset X\times_{k}X$ image de $\rho$. 
Comme $\rho$ est un monomorphisme, $R$ s'identifie via $\rho$ \`a $X\times_{k}G$ et l'on voit donc que $R$ est une relation d'\'equivalence \emph{plate}. Le
 fait que $X/G$ soit un $k$-espace r\'esulte donc de \cite[Corollary 6.3]{A2}. D'autre part nos hypoth\`eses impliquent que les diagrammes 
$$\xymatrix{X\times_{k}G\ar[r]^-{\alpha}\ar[d]_{\pr_{1}}& X\ar[d]^{{\pi}}\\
X\ar[r]^-{\pi}&X/G
}
\quad
\xymatrix{\quad\ar@{}[d]|(.6){\text{\normalsize et }}\\ \quad}
\quad
\xymatrix{X\times_{k}G\ar[r]^-{\rho}\ar[d]_{\pi\,\circ\,\pr_{1}}& X\times_{k}X\ar[d]^{{\pi\,\times\,\pi}}\\
X/G\ar[r]^-{\Delta}&(X/G)\times_{k}(X/G)
}$$
sont \emph{cart\'esiens}, $\Delta$ d\'esignant le morphisme diagonal. On en d\'eduit, gr\^ace au premier diagramme, que $X$ est un $G$-torseur sur $X/G$;
 le second donne les assertions de s\'eparation (remarquer que les fl\`eches verticales sont fid\`element plates).
 \cqfd
\Subsection{Produit contract\'e d'un torseur et d'un espace \`a op\'erateurs}\label{subsec-prodcontract}
\Subsubsection{Produits contract\'es. }Soit $G$ un faisceau en groupes sur $\Spec(k)$, op\'e\-rant \`a droite sur un $k$-faisceau $X$ et \`a gauche sur 
un $k$-faisceau $U$. Rappelons que le \emph{produit contract\'e} $X\PC{G} U$  est le faisceau quotient de $X\times_{k}U$ par l'action
 \`a droite de $G$  donn\'ee par $((x,u),g)\mapsto (xg,g^{-1}u)$. 

Nous supposerons dans la suite que $G$ est un $k$-groupe alg\'ebrique, que $X$ et $U$ sont des $k$-espaces, et de plus que $G$ op\`ere \emph{librement} sur $X$, de sorte que  $X$ est un $G$-torseur au-dessus du $k$-espace $Y:=X/G$ (\ref{QuotLibre}). 
Dans ce cas, $G$ op\`ere aussi librement sur  $X\times_{k}U$ via l'action ci-dessus, de sorte que $X\PC{G} U$ est un $k$-espace. En outre, 
on a un morphisme canonique $X\PC{G} U\to Y$, localement isomorphe
 (pour la topologie fppf) \`a la premi\`ere projection $Y\times U\to Y$;  cette construction commute \`a tout 
changement de $k$-espace de base $Y'\to Y$.

\Subsubsection{Torsion. }\label{torsion} 
\modif{Notons  $G'$ le $Y$-faisceau $\ul{\Aut}_{G}(X/Y)$ des automorphismes du $G_{Y}$-torseur $X$: c'est une forme tordue de $G_{Y}$, qui op\`ere naturellement \`a gauche sur $X$ par $Y$-morphismes, et aussi sur $X\PC{G}Z$ pour tout $Y$-faisceau $Z$ muni d'une action \`a gauche de $G_{Y}$. On obtient de cette fa\c{c}on un foncteur (\og torsion par $X$\fg) de la cat\'egorie des $Y$-faisceaux avec action \`a gauche de $G_{Y}$ vers celle des $Y$-faisceaux avec action \`a gauche de $G'$. Ce foncteur est une \emph{\'equivalence de cat\'egories} \cite[III, remarque 1.6.7]{Gi}. 

Lorsque $Y=\Spec(k)$, il s'ensuit formellement, en particulier, que l'on a une bijection naturelle entre les $k$-\emph{orbites} de $U$ sous $G$ et les $k$-orbites de $X\PC{G}U$ sous $G'$.
}

\Subsubsection{Faisceaux inversibles. } Si 
$L$ est un fibr\'e en droites $G$-lin\'earis\'e sur $U$, alors $M:=X\PC{G} L$ est de fa\c{c}on naturelle
 un fibr\'e en droites $G'$-lin\'earis\'e sur $X\PC{G} U$. 
 Supposons en outre que $Y$ et $U$ soient des vari\'et\'es, et  que $L$ soit ample sur $U$; alors  $M$ est ample sur $X\PC{G} U$ relativement \`a $Y$, et en particulier $X\PC{G} U$ 
est un $k$-sch\'ema quasi-projectif sur $Y$ (et m\^eme projectif si  $U$ est projectif sur $k$);  
voir {\cite[\S\,10.2, Lemma 6]{BLR}}.

\Subsection{Le plus grand sous-groupe lisse d'un groupe alg\'ebrique}\label{ssec:pgsgl}
Rappelons \cite[Lemma C.4.1]{CGP} qu'une $k$-vari\'et\'e $X$ admet un sous-sch\'ema ferm\'e canonique, que nous noterons
$X^\bec$,  
caract\'eris\'e comme le plus petit des sous-sch\'emas ferm\'es $X'$ tels que $X'(L)=X(L)$
 pour toute extension \emph{s\'eparable} $L$ de $k$; on peut  construire $X^\bec$ comme l'adh\'erence sch\'ematique de
 l'ensemble des
 points de $X$ \`a corps r\'esiduel s\'eparable (resp. fini s\'eparable) sur $k$. Ce sous-sch\'ema 
est fonctoriel en $X$ (pour les morphismes de $k$-vari\'et\'es), et sa formation commute aux extensions  s\'eparables des 
scalaires et 
aux produits de $k$-vari\'et\'es; en particulier, si $X$ est un $k$-groupe alg\'ebrique, alors $X^\bec$ est un
sous-groupe lisse de $X$ 
(et donc son plus grand sous-groupe lisse); on prendra garde qu'il n'est pas n\'ecessairement distingu\'e. 
\modif{Ceci \'etant, si on consid\`ere un produit semi-direct $N \rtimes G$ de $k$-groupes alg\'ebriques, l'isomorphisme
de $k$-sch\'emas $N^\bec \times_k G^\bec \simlgr (N \times G)^\bec$ indique que l'action de 
$G^\bec$ sur $N$ stabilise $N^\bec$ et que l'on a un isomorphisme $N^\bec \rtimes_k G^\bec \simlgr (N \rtimes_k G)^\bec$.
}

\begin{fsprop}\label{sous-groupe-lisse} Soient $G$ un $k$-groupe alg\'ebrique et $T$ un $G$-torseur \`a droite sur $\Spec{(k)}$.
Notons 
$\varepsilon\in (G/G^\bec)(k)$ la classe neutre, $Q$ le quotient $T/G^\bec$ 
et $L$ une extension \emph{s\'eparable} de $k$. Alors:\smallskip

\noindent\textup{(1)} On a $(G/G^\bec)(L)=(G/G^\bec)(k)=\{\varepsilon\}$.\smallskip

\noindent\textup{(2)} On a les \'equivalences:
$$
Q(k_{s})\neq\emptyset \;\Longleftrightarrow \; T(k_{s})\neq\emptyset\; \Longleftrightarrow \; \Card\,(Q(k_{s}))=1 \; \Longleftrightarrow \;\Card\,(Q(k))=1.
$$
\end{fsprop}
\dem (1) Soit $x\in (G/G^\bec)(L)$: comme le morphisme canonique $\pi:G\to G/G^\bec$ est lisse et surjectif, la fibre $G_{x}$ de $\pi$ en $x$ 
admet des points \`a corps r\'esiduel s\'eparable sur $L$, et donc sur $k$; ces points appartiennent donc \`a $G^\bec_{L}$, donc $x=\varepsilon$.\smallskip

\noindent\textup{(2)} Supposons $Q(k_{s})$ non vide; consid\'erant le morphisme lisse surjectif $T\to Q$, on voit comme en (1) que $T(k_{s})\neq\emptyset$. 
Ensuite, si $T(k_{s})\neq\emptyset$, alors  $T_{k_{s}}\cong G_{k_{s}}$ et donc $Q_{k_{s}}\cong (G/G^\bec)_{k_{s}}$, d'o\`u $\Card\,(Q(k_{s}))=1$ d'apr\`es (1).

Si $Q(k_{s})$ est r\'eduit \`a un point, alors ce point est $k$-rationnel par descente, donc $\Card\,(Q(k))=1$. Enfin, l'implication 
$\Card\,(Q(k))=1\Rightarrow Q(k_{s})\neq\emptyset$ est triviale.\cqfd

\begin{fscor}\label{quotient-torseur} Soient $G$ un $k$-groupe alg\'ebrique,  $Y$ une $k$-vari\'et\'e et $f:X\to Y$ un $G$-torseur sur $Y$. 
Posons $Z=X/G^\bec$ et soit $h:Z\to Y$ le morphisme canonique. Alors, 
pour toute extension s\'eparable $k'$ de $k$, l'application $Z(k')\to Y(k')$ induite par $f$ est injective.
\end{fscor}
\dem cela r\'esulte de \ref{sous-groupe-lisse} (2) appliqu\'e aux fibres de $h$.\cqfd

\medskip

La condition suivante intervient dans les compactifications de groupes \cite{Ga}:

\begin{fsdefi} \label{condition:star}
On dit qu'un $k$-groupe alg\'ebrique $G$ satisfait la condition $(*)$ 
si tous les $\overline k$-tores de $ G_{\overline k}$ sont 
des $\overline k$-tores de ${(G^\bec)}_{\overline k}$.  
\end{fsdefi}
Cette propri\'et\'e est \og insensible\fg\ aux extensions s\'eparables de 
$k$; plus pr\'e\-ci\-s\'e\-ment:
\begin{fslem}\label{fsExtSep} Soit $G$ un $k$-groupe alg\'ebrique, et soit $L$ une extension de $k$. Si $G$ v\'erifie  $(*)$, le $L$-groupe $G_{L}$ v\'erifie  $(*)$. 
La r\'eciproque est vraie si l'extension $L/k$ est s\'eparable.
\end{fslem}
\dem on sait que $(G^\bec)_{L}\subset(G_{L})^\bec$, avec \'egalit\'e lorsque $L/k$ est s\'eparable. Les assertions r\'esultent donc du lemme \ref{souslemmeExtSep} 
qui suit, appliqu\'e avec $H=G^\bec$.\cqfd
\begin{fsslem}\label{souslemmeExtSep} Pour un $k$-groupe alg\'ebrique $G$ et un sous-groupe lisse $H$ de $G$, d\'esignons par $*(k,H,G)$ la propri\'et\'e \og tout  
sous-$\overline k$-tore de $ G_{\overline k}$ est contenu dans $H_{\overline k}$\fg.

Alors, si $L$ est une extension de $k$, les propri\'et\'es $*(k,H,G)$ et $*(L,H_{L},G_{L})$ sont \'equivalentes. 
\end{fsslem}
\dem il est trivial que $*(L,H_{L},G_{L})$ entra\^{\i}ne $*(k,H,G)$. 

Pour la r\'eciproque on peut supposer $k$ alg\'ebriquement clos. Supposons que $*(L,H_{L},G_{L})$ ne soit pas v\'erifi\'ee:   soient donc ${\ol{L}}$  
une cl\^oture alg\'ebrique de $L$ et $T$ un sous-tore de $G_{\ol{L}}$ non contenu dans $H_{\ol{L}}$. Alors il existe un $k$-sch\'ema $Y$ affine, 
int\`egre et de type fini, dont le corps des fonctions est une sous-extension (de type fini) de $\ol{L}$, et un sous-$Y$-tore $\mathscr{T}$ de $G_{Y}$ 
qui n'est pas contenu dans $H_{Y}$. Comme $k$ est alg\'ebriquement clos, il existe un point $y\in Y(k)$ tel que la fibre $\mathscr{T}_{y}$ de $\mathscr{T}$ en $y$ 
(vue comme sous-tore de $G$) ne soit pas contenue dans $H$. Ainsi,  $*(k,H,G)$ n'est pas satisfaite.\cqfd
\begin{fslem} \label{lem:star} 
Soit $G$ un $k$-groupe alg\'ebrique.
\begin{numlist}
\item\label{lem:star1} Si $G$ est lisse, unipotent ou commutatif, il v\'erifie $(*)$.
\item\label{lem:star2} Soit $G'$ un $k$-sous-groupe de $G$ contenant la composante neutre $G^{\bec\circ}$ de $G^\bec$.
Si $G$ v\'erifie $(*)$, alors $G'$ v\'erifie $(*)$.
\item\label{lem:star3} Soit $G'$ un $k$-sous-groupe de $G$ tel que $G/G'$ soit fini. Pour que  $G$ v\'erifie $(*)$, il faut et il suffit que  $G'$
v\'erifie $(*)$.
\item\label{lem:star4} On suppose que $G$ est un $k$-sous-groupe distingu\'e 
d'un $k$-groupe lisse  $H$. Alors $G$ v\'erifie $(*)$.
\end{numlist}
\end{fslem}

\dem pour \ref{lem:star1}, le cas lisse et le cas unipotent sont \'evidents; pour le cas commutatif, remarquer que l'unique tore maximal de $G_{\ol{k}}$ est d\'efini sur $k$ \cite[C.4.4]{CGP} et donc contenu dans $G^\bec$. L'assertion \ref{lem:star2} est imm\'ediate, ainsi que \ref{lem:star3} car $G_{\ol{k}}$ et $G'_{\ol{k}}$ ont les m\^emes sous-tores.

Montrons \ref{lem:star4}. Le fait que $G$ soit distingu\'e dans $H$, et $H$ lisse,  entra\^ {\i}ne
que  $G^\bec$ est un sous-groupe \emph{distingu\'e} de $H$: en effet $H(k_{s})$ normalise $G^\bec_{k_{s}}$ et est dense dans $H$ puisque $H$ est lisse, donc le normalisateur de $G^\bec$ dans $H$ est \'egal \`a $H$.

Sans perte de g\'en\'eralit\'e, nous pouvons supposer $k$ s\'eparablement clos. En particulier le groupe  lisse $H$ admet un $k$-tore maximal 
$E_{0}$ et $T_0=(E_{0} \cap G)_{\red}^0$ est un $k$-tore maximal de $G$ et  a fortiori de  $G^\bec$.  
Soit $T$ un tore maximal de $G_{\overline k}$. Alors $T= (E \cap G_{\ol{k}})^0_{\red}$  pour un tore maximal $E$ de $H_{\overline k}$.
On sait qu'il existe $h \in H(\overline k)$ satisfaisant $E = h \,  E_{0 , \overline k} \,  h^{-1}$ \cite[C.4.5(1)]{CGP}.
Ainsi $T = h \, T_{0 , \overline k} \, h^{-1} \subset h\,  (G^\bec)_{\overline k}\, \,  h^{-1} = (G^\bec)_{\overline k}$.
\cqfd
\begin{fsrem}
Il existe des $k$-groupes r\'esolubles de dimension $1$  qui ne satisfont pas la condition
$(*)$, voir \S\,\ref{ex:OrbiteNonFermee}. 
\end{fsrem}
Dans ce travail, la condition $(*)$ interviendra via le lemme suivant, variante du lemme de Rosenlicht:
\begin{fslem}\label{Rosenlicht} Soit $\ol{k}$ une cl\^oture alg\'ebrique de $k$. Soit $\Gamma$ un $k$-groupe alg\'ebrique affine
 op\'erant sur une $k$-vari\'et\'e \emph{quasi-affine} $Z$. Soit $s\in Z(k)$, et soit $S\subset \Gamma$ le stabilisateur de $s$.
On suppose que tous les $\ol{k}$-tores de $\Gamma$ sont contenus dans $S$. Alors l'orbite $\Gamma.s$ de $s$ est ferm\'ee dans $Z$.
\end{fslem}
\dem on peut supposer $k$ alg\'ebriquement clos et $\Gamma$ lisse (remplacer $\Gamma$ par $\Gamma_{\mathrm{red}}$ ne change pas 
l'orbite ensemblistement). Alors $\Gamma.s$ est r\'eunion finie de transform\'es de $\Gamma^\circ
.s$ par des \'el\'ements de $\Gamma(k)$: il suffit donc de voir que l'orbite $\Gamma^\circ
.s$ sous $\Gamma^\circ$ est ferm\'ee. On supposera donc aussi $\Gamma$ connexe. Rempla\c{c}ant  $Z$ par l'adh\'erence de $\Gamma.s$, 
on supposera  $\Gamma.s$ dense dans $Z$. Dans ces conditions, le sous-groupe  $\Gamma_{t}$ de $\Gamma$ engendr\'e par ses sous-tores 
(qui est distingu\'e) op\`ere trivialement sur $Z$, de sorte que l'action se factorise par $U:=\Gamma/\Gamma_{t}$\,,
 groupe affine lisse connexe qui n'admet pas de sous-tore non trivial \cite[A.2.8]{CGP} et est donc \emph{unipotent}. 
Vu l'hypoth\`ese quasi-affine, le lemme de Rosenlicht \cite[XVII.5.7.3]{SGA3} montre alors que les orbites de $U$ dans $Z$ sont ferm\'ees.\cqfd

\Subsection{Compactifications partielles dans les sch\'emas de Hilbert ponctuels}
\label{orbite:bord}

\Subsubsection{Notations. }\label{notationsHilb}On fixe une $k$-alg\`ebre finie $k'$, de dimension  $d\geq1$, et  un $k'$-espace alg\'ebrique  $Q'$, s\'epar\'e et de type fini. On note $\sigma:Q'\to\Spec(k')$ le morphisme structural et:
\begin{itemize}
\item $V$ le $k$-espace induit $Q'\to\Spec(k')\to\Spec(k)$ (de sorte que $\sigma$ peut \^etre vu 
comme un morphisme \emph{de $k$-espaces} $V\xrightarrow{\sigma}\Spec(k')$);
\item $W= \W\limits_{k'/k} Q'$ le foncteur de restriction de Weil de $Q'$ relativement \`a $k'/k$.
\end{itemize}
Pour tout $k$-sch\'ema $S$, on a donc 
$$W(S)=Q'(k'\times_{k}S).$$
On consid\`ere d'autre part le foncteur de Hilbert $V^{[d]}=\Hilb^d_{V/k}$ des 
sous-$k$-espaces (ferm\'es) finis  de longueur $d$ de $V$. Explicitement, pour tout $k$-sch\'ema $S$, on a 
 $$\begin{array}{rcl}
V^{[d]}(S) &=&\Bigl\{  \text{sous-espaces }Z \subset V \times_k S,\\
&& \qquad\text{ finis localement libres
de rang $d$ sur $S$} \Bigr\}.
\end{array}
$$  
Un tel $Z$ est automatiquement un sch\'ema (affine sur $S$); on lui associe le  $S$-morphisme compos\'e
$$\varphi_{Z}: Z\inj V\times_{k}S \xrightarrow{\;\sigma\times_{k}\id_{S}\;} k'\times_{k}S.$$
Noter que $\varphi_{Z}$ est un morphisme de $S$-sch\'emas finis localement libres de m\^eme rang $d$.
\begin{fssrem}
Nous consid\'erons ici $W$ et $V^{[d]}$ comme des faisceaux sur le site fppf de $\Spec(k)$; nous n'aurons \`a utiliser aucun des r\'esultats profonds de repr\'esentabilit\'e connus, comme le fait que $V^{[d]}$ est un $k$-espace \cite[\S\,6]{A1} et est m\^eme  quasi-projectif si $V$ l'est (\cite[\S\,4]{TDTE4}, \cite[\S\,2]{B}, ou \cite[\S\,5.5]{Ni}). Le seul cas qui  nous servira dans l'article est celui o\`u $Q'$ est fini sur  $k'$, pour lequel la repr\'esentabilit\'e (et m\^eme la projectivit\'e de $V^{[d]}$) est tr\`es facile: voir le lemme \ref{lem:linearisation} plus bas.
\end{fssrem}
\Subsubsection{}\label{ImmOuverte} Rappelons que l'on a un morphisme de foncteurs
$$
u=u_{Q'/k'/k}: W  \, \to \, V^{[d]}
$$
d\'efini ainsi: si $S$ est un $k$-sch\'ema, on associe  \`a $w \in W(S)$
le $k'$-morphisme $w': S \times_k k' \to Q'$. Son graphe
$$
\Gamma_{w'} \subset Q' \times_{k'} (S \times_k k') = V \times_k S
$$
est un sous $S$-sch\'ema ferm\'e  de $V \times_k S$, 
isomorphe \`a $S \times_{k} k'$ et donc fini et libre de rang  $d$ sur $S$: c'est le point $u_{Q'/k'/k}(S)(w)\in V^{[d]}(S)$ voulu. 
\label{LemImmOuverte}\begin{flemstar}
\begin{numlist}
\item\label{LemImmOuverte1} Pour tout $k$-sch\'ema $S$, l'application $u_{Q'/k'/k}(S)$ ci-dessus induit une 
bijection de $W(S)$ sur l'ensemble des $Z\subset Q'  \times_k S$ tels que le $S$-morphisme $\varphi_{Z}$
 d\'efini en \rref{notationsHilb} soit un isomorphisme.

En particulier, le morphisme de foncteurs $u: W  \, \to \, V^{[d]}$ est re\-pr\'e\-sen\-table par une immersion ouverte.
\item\label{LemImmOuverte2} Si $k'$ est un corps, $u$ induit une bijection de $W(k)$ sur $V^{[d]}(k)$.
\end{numlist}
\end{flemstar}
\dem Un  sous-sch\'ema $Z\subset Q'  \times_k S$ tel que $\varphi_{Z}$ est un isomorphisme
d\'efinit une section de $Q' \times_k S= V \times_k S \to k' \times_k S$, c'est-\`a-dire  
un point de $W(S)=\Bigl(\prod_{k'/k} Q'\Bigr)(S)$. La partie \ref{LemImmOuverte1} est alors imm\'ediate et laiss\'ee au lecteur. 
L'assertion \ref{LemImmOuverte2} en r\'esulte: si $k'$ est un corps et si $Z\subset V$ est 
fini de rang $d$ sur $k$, alors $\varphi_{Z}:Z\to\Spec(k')$ est n\'ecessairement un isomorphisme.\cqfd

\begin{fssrems}De fa\c{c}on imag\'ee, on peut formuler \ref{LemImmOuverte2} en disant que $V^{[d]}$, vu comme 
compactification partielle de $W$, n'a pas de point $k$-rationnel \`a l'infini. Le m\^eme argument montre d'ailleurs 
qu'il n'a pas de point $K$-rationnel \`a l'infini, d\`es que $K$ est une extension de $k$ telle que $K\otimes_{k}k'$ soit un corps: 
le cas utile pour nous sera celui o\`u $k'$ (resp. $K$) est une extension radicielle (resp. s\'eparable) de $k$. \smallskip

Dans la suite, nous allons g\'en\'eraliser cette propri\'et\'e en rempla\c{c}ant les points par les orbites sous l'action 
d'un groupe alg\'ebrique.
\end{fssrems}
\Subsubsection{Version \'equivariante: orbites \`a l'infini. }\label{notations: pas d'orbite} On garde les notations
 de \ref{notationsHilb} et \ref{ImmOuverte}, et l'on se donne de plus un $k$-groupe alg\'ebrique $G$ op\'erant \`a gauche 
sur $V$ \emph{par $k'$-automorphismes}, c'est-\`a-dire que $\sigma:V\to\Spec(k')$ est $G$-invariant). Il reviendrait au m\^eme de
 se donner une action du $k'$-groupe $G_{k'}$ sur $Q'$; l'action de $G$ sera plus commode \`a utiliser). 

On en d\'eduit formellement des actions de $G$ (vu comme foncteur en groupes) sur les foncteurs $W$ et $V^{[d]}$, et 
 l'immersion ouverte $u$ est \'equivariante pour ces actions.

\begin{fsteo}\label{th:pas-d'orbite} Avec les hypoth\`eses et notations de \rref{notations: pas d'orbite}, soient $K$ 
une extension de $k$ et $J \subset V^{[d]}_{K}$ une $K$-orbite  pour l'action de $G_{K}$ \emph{(cf. \ref{def:orbite})}. 

Alors, si $K\otimes_{k}k'$ est un  corps, on a $J\subset W_{K}$.
\end{fsteo}
\dem  en rempla\c{c}ant $k'$ par  $K\otimes_{k}k'$, $Q'$ par $Q'_{K}$, etc., on peut sup\-po\-ser que $K=k$
 (et que $k'$ est un corps). Alors $y:J\to V^{[d]}$ d\'efinit  un sous-espace $Z$ de $Q'\times_{k}J$, fini 
localement libre de rang $d$ sur $J$, et  stable sous $G$ (op\'erant sur $J\times_{k}Q'$ par l'action produit). 
Par \ref{LemImmOuverte}, il s'agit de montrer que le $J$-morphisme
$$\varphi_{Z}: Z \ffl k'\times_{k}J$$
d\'efini en \ref{notationsHilb} est un isomorphisme. Il est d\'ecrit par un morphisme $G$-\'equivariant
$$\psi_{Z}:k'\otimes_{k}\cO_{J}\ffl \cA:=\pr_{2*}(\cO_{Z})$$
de $\cO_{J}$-Alg\`ebres finies  localement libres de rang $d$, lin\'earis\'ees pour l'action de $G$ sur $J$. En particulier 
le conoyau $\cC$ de $\psi_{Z}$ est un $\cO_{J}$-module coh\'erent $G$-lin\'earis\'e. Consid\'erons, pour chaque 
$s\in\NN$, la strate de Fitting (de rang $s$)  $F_{s}(\cC)\subset J$: 
c'est un sous-sch\'ema localement ferm\'e de $J$, tel que la restriction de $\cC$ \`a $F_{s}(\cC)$ soit localement libre de rang $s$, 
et universel pour cette propri\'et\'e (cf. \cite[\S 11.8]{GW})  Puisque $\cC$ est $G$-lin\'earis\'e, chaque $F_{s}(\cC)$ est stable par $G$; 
comme $G$ op\`ere transitivement sur $J$, on a donc $F_{s}(\cC)=\emptyset$ sauf pour
 une valeur $r$ de $s$, pour laquelle $F_{r}(\cC)=J$. En d'autres termes, $\cC$ est \emph{localement libre de rang constant $r$}. 

Il reste \`a voir que $r=0$. Puisque $\cC$ est localement libre de rang $r$, l'image de $\psi_{Z}$ est une alg\`ebre quotient 
de $k'\otimes_{k}\cO_{J}$, localement libre de 
rang $d-r$, et d\'efinit donc un sous-sch\'ema $T\subset k'\times_{k}J$, fini localement libre sur $J$ de rang $d-r$ et $G$-invariant. 
Comme $G$ op\`ere transitivement sur $J$ et \emph{tri\-via\-le\-ment} sur $\Spec(k')$, $T$ provient par descente d'un sous-$k$-sch\'ema 
$T_{0}$ de rang $d-r$ de $\Spec(k')$. Puisque $k'$ est un corps, on a soit $T_{0}=\emptyset$ et $r=d$, soit $T_{0}=\Spec(k')$ et $r=0$; 
mais le premier cas est exclu car $\psi_{Z}$ est un morphisme d'alg\`ebres non nulles, donc n'est pas nul.\cqfd

\begin{fslem}\label{lem:linearisation} Soient $Y$ un $k$-sch\'ema \emph{fini}, et $d\in\NN$. Alors le foncteur $Y^{[d]}:=\Hilb^d_{Y/k}$
des sous-sch\'emas finis de longueur $d$ de $Y$ est repr\'esentable par un $k$-sch\'ema projectif, muni d'un faisceau ample $\ul{\Aut}_{k}(Y)$-lin\'earis\'e.
\end{fslem}

\dem  On pose $Y=\Spec{(A)}$ o\`u $A$ est une $k$-alg\`ebre finie; on note $A_{\mathrm{lin}}$ le $k$-espace vectoriel sous-jacent \`a $A$, 
et $G$ le $k$-groupe $\ul{\Aut}_{k}(Y)$; c'est \'evidemment un sous-sch\'ema en groupes ferm\'e de $\GL(A_{\mathrm{lin}})$, et en 
particulier un $k$-groupe alg\'ebrique affine. 

La projectivit\'e de $Y^{[d]}$ est un fait g\'en\'eral \cite[prop 2.13]{B}, mais se voit facilement ici: en effet, 
pour tout $k$-sch\'ema $S$, un point de $Y^{[d]}(S)$ est la m\^eme chose qu'une $\OOO_{S}$-alg\`ebre 
quotient de $\OOO_{S}\otimes_{k}A$, localement libre de rang $d$ comme $\OOO_{S}$-module. Ainsi, $Y^{[d]}$ est de 
fa\c{c}on naturelle un sous-sch\'ema ferm\'e de la grassmannienne $\mathrm{Gr}_{d}(A_{\mathrm{lin}})$ des quotients
 de rang $d$ de $A_{\mathrm{lin}}$. De plus, l'action de $G$ sur $Y^{[d]}$ est induite par son action naturelle sur 
 $\mathrm{Gr}_{d}(A_{\mathrm{lin}})$, laquelle se 
factorise par celle de $\mathrm{GL}_{k}(A_{\mathrm{lin}})$, de sorte que le faisceau ample
 canonique sur $\mathrm{Gr}_{d}(A_{\mathrm{lin}})$ (qui est $\mathrm{GL}_{k}(A_{\mathrm{lin}})$-lin\'earis\'e) 
induit un faisceau ample $G$-lin\'earis\'e sur $Y^{[d]}$.\cqfd

\section{Corps topologiquement hens\'eliens; le cas des torseurs sous un groupe lisse}\label{sec:CorpsTop}

On d\'esigne par $F$ un corps topologique s\'epar\'e; 
on notera $\wh{F}$ son compl\'et\'e (qui est une $F$-alg\`ebre topologique, et est 
un corps si la topologie de $F$ est d\'efinie par une valuation). 

\Subsection{Vari\'et\'es sur un corps topologique; corps topologiquement hen\-s\'e\-liens}\label{TopVar} 

On note $\VAR_F$ la cat\'egorie des $F$-vari\'et\'es, $\TOP$ celle des espaces 
topologiques, $\ENS$ celle des ensembles.

On peut (voir \cite[app. III]{W}, \cite[p. 256]{KS}) 
d'une mani\`ere et d'une seule, associer \`a toute $F$-vari\'et\'e $X$  (ou plus g\'en\'eralement \`a tout $F$-sch\'ema localement de type fini)
une topologie sur l'ensemble $X(F)$, de telle sorte que les conditions suivantes 
soient v\'erifi\'ees (dans lesquelles $X$ et $Y$ d\'esignent des $F$-vari\'et\'es quelconques):

\begin{romlist}
\item\label{TopVar1} si $X=\Aa^1_F$, la bijection naturelle de $X(F)$ sur $F$ est un hom\'eo\-mor\-phisme;
\item\label{TopVar2} si $f: X\to Y$ est un $F$-morphisme, l'application induite $f(F): X(F)\to Y(F)$ est continue;
 si de plus $f$ est une immersion ouverte (resp. ferm\'ee), alors $f(F)$ est un plongement 
topologique ouvert (resp. ferm\'e);
\item\label{TopVar3} la bijection naturelle de $(X\times_F Y)(F)$ sur $X(F)\times Y(F)$ est un 
hom\'eo\-mor\-phisme;
\item\label{TopVar4} si $X$ est un sch\'ema s\'epar\'e, alors $X(F)$ est un espace s\'epar\'e.
\end{romlist}
(Ces conditions ne sont pas ind\'ependantes: par exemple, \ref{TopVar4} est con\-s\'e\-quence de  \ref{TopVar3} et de  \ref{TopVar2} appliqu\'ees au morphisme diagonal $X\to X\times_{F}X$.)

On obtient ainsi un foncteur, not\'e $X\mapsto X(F)_\Top$ 
(ou encore $X\mapsto X_\Top$),
 de $\VAR_F$ dans $\TOP$.  Il r\'esulte
 facilement des propri\'et\'es ci-dessus que ce foncteur transforme les immersions en plongements 
topologiques et commute aux produits fibr\'es.

\begin{fsdefi}\label{DefHens} Un corps topologique s\'epar\'e $F$ est
 \emph{topologiquement hens\'elien} si  pour tout morphisme \emph{\'etale} 
$f:X\to Y$ de $F$-vari\'et\'es, l'application induite $f_{\Top}$ est un \emph{hom\'eomorphisme local}.
\end{fsdefi}
\begin{fssrems}
Si l'on se limite aux topologies d\'eduites de valuations ou de valeurs absolues (\og V-topologies\fg), 
les corps topologiquement hens\'eliens sont appel\'es \emph{t-hens\'eliens} dans l'article \cite{PZ}, qui
en donne diverses caract\'erisations. 

Outre les corps valu\'es hens\'eliens (proposition \ref{Hens=hens} plus bas), les corps $\RR$ et $\CC$ sont 
topologiquement hens\'eliens (pour leur topologie classique), de m\^eme que leurs sous-corps $\ol{\QQ}$ et $\ol{\QQ}\cap\RR$ et
que tout corps r\'eel clos \modif{muni de la topologie de l'ordre}.
\end{fssrems}
\begin{fslem}\label{lisse-ouv} Soit $F$ un corps  topologiquement hens\'elien, et soit 
$f:Y\to X$ un morphisme \emph{lisse} de $F$-vari\'et\'es. 

Alors, pour tout  $y \in Y(F)$, il existe un voisinage ouvert $\Omega$ de
$f(y)$  dans $X_\Top$ tel que l'application induite $f_\Top:  Y_\Top \to X_\Top$ 
admette une section continue $\Omega \to Y_\Top$. 

En particulier, $f_{\Top}$ est \emph{ouverte}.
\end{fslem}
\dem remarquer que si $y\in Y(F)$, il existe un morphisme  $h:Z\to Y$ et un point $z\in Z(F)$ tels que $h(z)=y$ et que le compos\'e 
$f\circ h:Z\to X$ soit \'etale \cite[\S\,2.2, prop. 14]{BLR}. \cqfd

\begin{fsrem}
R\'eciproquement, soit $F$ un corps topologique s\'epar\'e tel que $f_{\Top}$ soit ouverte pour tout morphisme lisse $f:X\to Y$ de $F$-vari\'et\'es.
Alors $F$ est 
topologiquement hens\'elien: le lecteur le v\'erifiera en utilisant le fait  que si $f$ est \'etale, le morphisme diagonal 
$X\to X\times_{Y} X$ est une immersion ouverte. (Cette remarque n'est pas utilis\'ee dans la suite).
\end{fsrem}

Nous allons justifier cette terminologie  par l'\'enonc\'e suivant (sans doute bien connu, mais difficile \`a trouver sous cette forme dans la litt\'erature).

\begin{fsprop} \label{Hens=hens} Tout corps valu\'e hens\'elien est topologiquement hen\-s\'e\-lien. 
\end{fsprop}

\noindent{\it D\'emonstration.} Soient $F$ un corps valu\'e hens\'elien et $f: X \to Y$ un morphisme \'etale
de $F$-vari\'et\'es. Soit $x \in X(F)$  et soit $y$ son image dans $Y(F)$: montrons que $f_{\Top}$ est un 
hom\'eomorphisme local au point $x$. Bien entendu, il nous suffira pour cela de supposer que $f$ est \'etale \emph{au point $x$}, 
ce que nous ferons syst\'ematiquement dans les r\'eductions qui vont suivre.%

On peut supposer que $Y= \Spec{(A)}$ est affine et  que $X$ est un ouvert
de $\Spec{\left(A[T]/P(T)\right)}$ o\`u $P\in A[T]$ est un polyn\^ome (unitaire, si l'on veut) tel que
$P_y(T)\in F[T]$ soit s\'eparable \cite[\S\,2.3, prop. 3]{BLR}. Sans perte de g\'en\'eralit\'e, on peut supposer
que $X=\Spec{\left(A[T]/P(T) \right)}$.
 On plonge $Y$ dans un espace affine 
$\widetilde Y= \Aa^n_F= \Spec{(\widetilde A)}$ et on rel\`eve $P(T)$ en un polyn\^ome unitaire 
$\widetilde P(T) \in \widetilde A[T]$. Rempla\c{c}ant $Y$ par $\til{Y}$
 et $P$ par $\til{P}$, nous sommes ramen\'es \`a la situation o\`u 
$Y=\Aa^n_{F}=\Spec{F[\ul{Z}]}=\Spec{F[Z_{1},\dots,Z_{n}]}$, $X=\Spec{\left(F[\ul{Z},T]/(P)\right)}$ o\`u 
$P$ est unitaire en $T$, et o\`u $y$ (resp. $x$) est l'origine de $\Aa^n_{F}$ (resp. de $\Aa^{n+1}_{F}$). 
En outre, la projection $X\to Y$ \'etant \'etale en $x$, nous pouvons choisir les coordonn\'ees de mani\`ere
que l'hyperplan tangent en $x$ \`a $X$ ait pour \'equation $T=0$, de sorte que $P$ est 
(\`a un scalaire inversible pr\`es) de la forme  $$P(\ul{Z},T)=T+\sum\limits_{\vert I\vert+j\geq2}a_{I,j}\,\ul{Z}^I T^j \qquad(a_{I,j}\in F).$$
Notons $R$ l'anneau de la valuation $v$, et $\Fm$ son id\'eal maximal. Pour $\alpha\in F^\times$, nous pouvons 
remplacer $P$ par  $P_{\alpha}(\ul{Z},T):=\frac{1}{\alpha}P(\alpha\ul{Z},\alpha T)$. 
Le coefficient de $\ul{Z}^I T^j$ dans $P_{\alpha}$ est $\alpha^{\vert I\vert+j-1}a_{I,j}$: dans cette formule 
l'exposant de $\alpha$ est $>0$ donc,  prenant  $\alpha$ assez proche de $0$, on peut supposer que les 
coefficients $ a_{I,j}$ sont dans $\Fm$. La propri\'et\'e de Hensel montre alors que pour chaque
 $\ul{z}\in R^n$, le polyn\^ome $P(\ul{z},T)\in R[T]$ admet une unique racine $t(\ul{z})$ 
dans $\Fm$. En d'autres termes, $f_{\Top}$ 
induit une bijection entre $f_{\Top}^{-1}(R^n)\cap (R^n\times \Fm)$ (voisinage de $x$ dans $X_{\Top}$) et 
$R^n$ (voisinage de $y$ dans $F^n$). 

Il reste \`a voir que l'application $\ul{z}\mapsto t(\ul{z})$ est continue \`a l'origine, et pour cela il suffit 
d'\'etablir que $v(t(\ul{z}))\modiff{\geq}\min\limits_{i=1,\dots,n}v(z_{i})$. En d'autres termes, nous devons montrer que 
si $\ul{z}=(z_{1},\dots,z_{n})\in R^n$ et $t\in\Fm$ v\'erifient $v(t)\modiff{<}\min\limits_{i=1,\dots,n}v(z_{i})$ \modiff{(ce qui implique que $t\neq0$)}, alors on a $P(\ul{z},t)\neq0$. 
Il suffit pour cela de voir que chaque terme $\ul{z}^I\,t^j$, pour $\vert I\vert+j\geq2$, est de valuation $\modiff{>} v(t)$ (rappelons que  $v(a_{I,j})>0$).
Or c'est vrai si $\vert I\vert\geq1$ vu l'hypoth\`ese sur les $v(z_{i})$, et c'est vrai si  $I=0$ car alors $j\geq2$ et $t\in\Fm\modiff{\setminus\{0\}}$.
\cqfd

\Subsection{Extension aux espaces alg\'ebriques}\label{sec:espalg}

\`A partir de \ref{prop:caractTopEspAlg}, le corps topologique $F$ sera toujours suppos\'e \emph{topologiquement hens\'elien}.

On identifie la cat\'egorie $\VAR_{F}$ \`a une sous-cat\'egorie pleine de la cat\'egorie  $\mathrm{ALG}_{F}$ 
des $F$-espaces (c'est-\`a-dire, rappelons-le, des $F$-espaces alg\'ebriques quasi-s\'epar\'es  de type fini). 
On se propose de d\'efinir, pour tout $F$-espace $X$, une topologie sur $X(F)$, qui co\"{\i}ncide avec celle d\'ej\`a
d\'efinie  lorsque $X$ est un sch\'ema. 
Le cas d'un corps valu\'e complet de rang $1$ est trait\'e dans \cite[\S 5]{ConAdelic}. 
Pour ne pas allonger d\'emesur\'ement ce travail, nous omettons ici la plupart des d\'emonstrations.

\begin{fsdefi}\label{def:TopEspAlg} Soit $X$ un $F$-espace. On munit $X(F)$ de la topologie suivante: une partie $\Omega\subset X(F)$ est ouverte si et seulement si, pour toute $F$-vari\'et\'e $Z$ et tout $F$-morphisme $\varphi:Z\to X$, l'ensemble $\varphi^{-1}(\Omega)\subset Z(F)$ est un ouvert de $Z_{\Top}$.

L'espace topologique ainsi obtenu sera not\'e $X_{\Top}$ ou $X(F)_{\Top}$.
\end{fsdefi}

\begin{fsrem}\label{remDefTopEspAlg}
En d'autres termes, on munit $X(F)$ de la topologie la plus fine rendant continues toutes les applications $\varphi(F)$ lorsque $\varphi$ 
parcourt tous les $F$-morphismes d'une $F$-vari\'et\'e vers $X$. Noter qu'il est clair que l'on peut se limiter aux vari\'et\'es $Z$ qui sont \emph{affines}.
\end{fsrem} 

\begin{fsprop}\label{prop:soritesTopEspAlg1} 
\ 

\begin{numlist} 
\item\label{prop:soritesTopEspAlg11}  Si $X$ est une $F$-vari\'et\'e, la topologie d\'efinie en \rref{def:TopEspAlg} (o\`u $X$ est vue comme $F$-espace) co\"{\i}ncide 
avec celle d\'ej\`a d\'efinie en \rref{TopVar}.
\item\label{prop:soritesTopEspAlg12} Soit $f:X\to Y$ un morphisme de $F$-espaces. Alors:
\begin{romlist}
\item\label{prop:soritesTopEspAlg121} l'application $f_{\Top}:X_{\Top}\to Y_{\Top}$ induite par $f$ sur les points $F$-rationnels est continue;
\item\label{prop:soritesTopEspAlg122} si $f$ est une immersion  (resp. une immersion ouverte, resp. une immersion ferm\'ee), alors $f_{\Top}$ est un hom\'eomorphisme 
sur son image, qui est localement ferm\'ee (resp. ouverte, resp. ferm\'ee) dans $Y_{\Top}$.
\end{romlist}
\end{numlist}  \end{fsprop}
\dem les assertions \ref{prop:soritesTopEspAlg11} et  \ref{prop:soritesTopEspAlg12}\,\ref{prop:soritesTopEspAlg121} sont imm\'ediates \`a partir de 
la d\'efinition \ref{def:TopEspAlg}. 
Montrons  \ref{prop:soritesTopEspAlg12}\,\ref{prop:soritesTopEspAlg122}: si $f$ est une immersion alors $f_{\Top}$ est injective, donc il suffit de 
voir que si $f$ est une immersion ouverte (resp. ferm\'ee) alors $f_{\Top}$ est une application ouverte (resp. ferm\'ee), le cas d'une
 immersion arbitraire en r\'esultant par composition. 

Traitons le cas d'une immersion ouverte (celui d'une immersion ferm\'ee est enti\`erement analogue). Soit  $\Omega$ un ouvert de $X_{\Top}$, et
 montrons que $f_{\Top}(\Omega)$ est ouvert dans $Y_{\Top}$: soit donc $h:Y'\to Y$ un $F$-morphisme, o\`u $Y'$ est une vari\'et\'e. Il s'agit
 de voir que $h_{\Top}^{-1}(f_{\Top}(\Omega))$ est un ouvert de $Y'_{\Top}$. Posons $X'=X\times_{Y} Y'$ et soient $f':X'\to Y'$ et $h':X'\to X$ les
 morphismes \'evidents. Alors $f'$ est une immersion ouverte, et en particulier $X'$ est un sch\'ema. De plus 
$h_{\Top}^{-1}(f_{\Top}(\Omega))=f'_{\Top}({h^{\prime-1}_{\Top}}(\Omega))$ qui est bien ouvert dans $Y'_{\Top}$ puisque
 ${h^{\prime-1}_{\Top}}(\Omega)$ est ouvert dans $X'_{\Top}$ et que $f'$ est une immersion ouverte entre vari\'et\'es.\cqfd

\medskip

Lorsque $F$ est topologiquement hens\'elien, on dispose d'une autre caract\'erisation de la topologie de \ref{def:TopEspAlg}, en termes de 
recouvrements \'etales (cf. \ref{prop:caractTopEspAlg} plus bas). Pour y parvenir, nous aurons besoin  de l'\'enonc\'e suivant, variante d'un
 r\'esultat de Gruson et Raynaud \cite[I, proposition 5.7.6]{GR}, et de son corollaire:

\begin{fsprop}\textup{\cite[Theorem 3.1.1]{CLO}\label{prop:FiltrationCLO}} Soit $X$ un espace alg\'ebrique quasi-compact et quasi-s\'epar\'e. 
Il existe une suite croissante finie
$$\emptyset=U_{0}\subset U_{1}\subset\dots\subset U_{r}=X$$
d'ouverts quasi-compacts de $X$, avec la propri\'et\'e suivante: pour chaque $i>0$, si l'on note $Z_{i}$ le sous-espace ferm\'e r\'eduit de $U_{i}$ 
compl\'ementaire de $U_{i-1}$, il existe un morphisme \'etale surjectif $\pi_{i}:Y_{i}\to U_{i}$ qui est un isomorphisme au-dessus de $Z_{i}$, et 
o\`u $Y_{i}$ est un sch\'ema quasi-compact et s\'epar\'e.\cqfd
\end{fsprop}

\begin{fscor}\label{cor:FiltrationCLO} Soit $X$ un espace  alg\'ebrique quasi-compact et quasi-s\'epar\'e. Il existe un sch\'ema affine $Y$ et un 
morphisme \'etale surjectif $\pi:Y\to X$ tel que pour tout anneau artinien $A$ l'application $Y(A)\to X(A)$ induite par $\pi$ soit surjective.\cqfd
\end{fscor}

\begin{fsprop}\label{prop:caractTopEspAlg} On suppose $F$ topologiquement hens\'elien.

Soit $Y$ un $F$-espace. Fixons  un $F$-morphisme \'etale $\pi:X\to Y$, o\`u $X$ est une  $F$-vari\'et\'e s\'epar\'ee. 

Alors $\pi_{\Top}:X_{\Top}\to Y_{\Top}$ est \emph{ouverte}.

En particulier, si de plus $\pi$ induit une \emph{surjection} de $X(F)$ sur $Y(F)$ \emph{(\'etant donn\'e $Y$, 
un tel morphisme $\pi$ existe, d'apr\`es \ref{cor:FiltrationCLO})}, alors 
$Y_{\Top}$ s'identifie au quotient de $X_{\Top}$ par la relation d'\'equivalence d\'efinie par $\pi_{\Top}$.\cqfd
\end{fsprop}

\subsubsection{Produits fibr\'es. }\label{prop:ProdFibEspAlg}Consid\'erons un diagramme \emph{cart\'esien}
$$\xymatrix{X\times_{S}Y\ar[r]^(.6){g'}\ar[d]_{f'}&X\ar[d]^{f}\\ 
Y\ar[r]^{g}&S}$$
de $F$-espaces. 
Notons $j:X\times_{S}Y\to X\times_{F}Y$ le monomorphisme naturel: il induit une injection continue
$$j_{\Top}:(X\times_{S}Y)_{\Top}\inj (X\times_{F}Y)_{\Top}.$$
On a d'autre part une application continue naturelle
$$\alpha:(X\times_{S}Y)_{\Top}\to X_{\Top}\times_{S_{\Top}}Y_{\Top}$$
d\'eduite des propri\'et\'es universelles, et \'evidemment \emph{bijective}.

\begin{fpropstar} Avec les hypoth\`eses et notations ci-dessus, on suppose $F$ topologiquement hens\'elien. 
\begin{numlist}
\item\label{prop:ProdFibEspAlg1} Les conditions suivantes sont \'equivalentes:
\begin{romlist}
\item\label{prop:ProdFibEspAlg11} $\alpha$ est un hom\'eomorphisme;
\item\label{prop:ProdFibEspAlg12} $j_{\Top}$ est un hom\'eomorphisme  sur son image.
\pauseromlist
\item\label{prop:ProdFibEspAlg2} Les conditions de \rref{prop:ProdFibEspAlg1} sont satisfaites dans les deux cas suivants:
\finpauseromlist
\item\label{prop:ProdFibEspAlg21} l'un des morphismes $f$, $f'$, $g$, $g'$ est une immersion;
\item\label{prop:ProdFibEspAlg24} $S$ est localement s\'epar\'e.\cqfd
\end{romlist}
\end{numlist} 
\end{fpropstar}

\begin{fssrems}\label{rem:ProdFibEspAlg1} 
\begin{numlist}
\item\label{rem:ProdFibEspAlg11} Bien entendu, un cas particulier de \ref{prop:ProdFibEspAlg24} est celui o\`u $S=\Spec{(F)}$. 
En d'autres termes, le foncteur $X\mapsto X_{\Top}$ \emph{commute aux produits} de $F$-espaces.
\item\label{rem:ProdFibEspAlg12} Le cas \ref{prop:ProdFibEspAlg21} n'utilise pas l'hypoth\`ese que $F$ est topologiquement hens\'elien. 
\item\label{rem:ProdFibEspAlg13}Un cas particulier utile de \ref{prop:ProdFibEspAlg21} est le fait que pour $f:X\to S$ quelconque
 et $s\in S(F)$, la fibre $f_{\Top}^{-1}(s)$ s'identifie \`a $(f^{-1}(s))_{\Top}$, l'inclusion d'un point rationnel \'etant toujours une immersion.
\item\label{rem:ProdFibEspAlg14} On verra en \ref{ssec:ExNonLocSep} un exemple o\`u les conditions de \ref{prop:ProdFibEspAlg1} ne sont pas satisfaites.
\end{numlist}
\end{fssrems}
\begin{fsprop}\label{prop:topLisseEtale} \emph{(morphismes lisses et \'etales)} On suppose $F$ topologiquement hens\'elien.
Soit $f:X\to Y$ un morphisme \emph{lisse} de $F$-espaces. Alors:
\begin{numlist}
\item\label{prop:topLisseEtale1}  $f_{\Top}:X_{\Top}\to Y_{\Top}$ est ouverte.
\item\label{prop:topLisseEtale2} On suppose que  $f$ est \'etale et que, de plus,  la bijection naturelle 
\break $(X\times_{Y}X)_{\Top}\to X_{\Top}\times_{Y_{\Top}}X_{\Top}$ est un hom\'eomorphisme 
\emph{(condition v\'erifi\'ee notamment, d'apr\`es \ref{prop:ProdFibEspAlg}, si $Y$ est localement s\'epar\'e)}. Alors
 $f_{\Top}$ est un ho\-m\'eo\-mor\-phisme local.
\item\label{prop:topLisseEtale3} Si $Y$ est localement s\'epar\'e, alors $f_{\Top}$ \og a des sections locales\fg\ au 
sens suivant: pour tout $x\in X_{\Top}$, il existe un voisinage $V$ de $f(x)$ dans $Y_{\Top}$ et une section continue 
$s:V\to X_{\Top}$ de $f$ sur $V$ telle que $s(f(x))=x$.
\cqfd
\end{numlist}
\end{fsprop}

\Subsection{Corps valu\'es hens\'eliens: utilisation de mod\`eles entiers}
On suppose dans cette section que $F$ est un corps muni d'une valuation (non triviale) $v$ d'anneau $A$, et de la topologie associ\'ee. 
On notera $\FI$ l'ensemble des id\'eaux stricts et non nuls de $A$, $\Fm\in\FI$ son id\'eal maximal, $k=A/\Fm$ son corps r\'esiduel.

On se propose de d\'ecrire l'espace $X_{\Top}$, o\`u $X$ est une $F$-vari\'et\'e, en termes des \og mod\`eles entiers\fg\ de $X$, c'est-\`a-dire 
des $A$-espaces alg\'ebriques $\cX$ de type fini munis d'un isomorphisme $\cX_{F}\flis X$. Comme dans la section \ref{sec:espalg}, 
les d\'emonstrations sont omises.

\subsubsection{Topologie sur les points entiers. }\label{not:top-entiere}
 Soit $\cX$ un $A$-espace alg\'ebrique quasi-s\'epar\'e  
 de type fini. (Le lecteur pourra constater que les constructions qui suivent ont un sens pour des foncteurs plus g\'en\'eraux).
 
 On notera
$$\begin{array}{rccll}
\rho_{\cX,J}:&\cX(A)&\ffl&\cX(A/J) &(J\in\FI)\\
\rho_{\cX}:&\cX(A)&\ffl&\varprojlim\limits_{J\in\FI}\cX(A/J)
\end{array}$$
les applications \'evidentes; on utilisera aussi la notation \og$x\bmod J$\fg\ pour $\rho_{\cX,J}(x)$. 

On munira chaque $\cX(A/J)$ de la topologie discr\`ete, $\varprojlim\limits_{J\in\FI}\cX(A/J)$ de la topologie limite projective, et  l'ensemble $\cX(A)$ de la topologie induite via $\rho_{\cX}$. L'espace topologique obtenu sera not\'e $\cX(A)_{A,\Top}$, ou $\cX_{A,\Top}$.\smallskip

\begin{fssprop}\label{TopEnt:injective} L'application $\rho_{\cX}$ est injective. En particulier, l'espace $\cX_{A,\Top}$ est \emph{s\'epar\'e}.
\end{fssprop}
\dem soient $x$ et $x'$ dans $\cX(A)$ ayant m\^eme image par $\rho_{\cX}$. Le noyau de la double fl\`eche $\Spec\,(A)\underset{x'}{\overset{x}{\rightrightarrows}}\cX$ est alors re\-pr\'e\-sen\-t\'e par un mono\-mor\-phisme de $A$-sch\'emas $Z\to\Spec\,(A)$; d'autre part, le fait que 
$\rho_{\cX}(x)=\rho_{\cX}(x')$ entra\^{\i}ne l'existence, pour chaque $J\in\FI$, d'un $A$-morphisme $s_{J}:\Spec\,(A/J)\to Z$. Comme $Z$ est un sch\'ema (et les $s_{J}$ automatiquement compatibles), on en d\'eduit un $A$-morphisme $\Spec\,(\wh{A})\to Z$. Ainsi, $Z\to\Spec\,(A)$ est \`a la fois un monomorphisme de sch\'emas et un \'epimorphisme de faisceaux fpqc (rappelons que $\wh{A}$ est fid\`element plat sur $A$); c'est donc un isomorphisme, de sorte que $x=x'$.\cqfd\medskip

\noindent Si $x\in\cX(A)$ et $J\in\FI$, on notera $B_{\cX}(x,J)$ la \og boule de rayon $J$\fg
$$B_{\cX}(x,J)=\rho_{X,J}^{-1}\left(\rho_{X,J}(x)\right)=\left\{x'\in\cX(A)\mid x'\bmod J=x\bmod J\right\}$$
(de sorte que la famille $\left(B_{\cX}(x,J)\right)_{J\in\FI}$ est une base de voisinages ouverts de $x$ dans $\cX_{A,\Top}$).\smallskip

Si $\varphi:\cX\to\cY$ est un morphisme d'espaces alg\'ebriques quasi-s\'epar\'es de type fini sur $A$, l'application naturelle $\varphi_{A,\Top}:\cX_{A,\Top}\to\cY_{A,\Top}$ est continue, et plus pr\'ecis\'ement envoie $B_{\cX}(x,J)$ dans  $B_{\cY}(\varphi(x),J)$ pour tout $x\in \cX(A)$ et tout $J\in\FI$.

\begin{fsprop}\label{SoritesTopEnt} On garde les notations de \rref{not:top-entiere}; en particulier les lettres $\cX$, $\cY$ d\'esignent des $A$-espaces alg\'ebriques quasi-s\'epar\'es de type fini, et $\varphi:\cX\to\cY$ un $A$-morphisme.
\begin{numlist}
\item\label{SoritesTopEnt1} Le foncteur $\cX\mapsto\cX_{A,\Top}$ (de la cat\'egorie des $A$-espaces alg\'ebriques quasi-s\'epar\'es de type fini dans celle des espaces topologiques) commute aux produits fibr\'es finis.
\item\label{SoritesTopEnt3} Si $\varphi$ est un monomorphisme, alors $\varphi_{A,\Top}$ est un plongement topologique ferm\'e.
\item\label{SoritesTopEnt4} Supposons $A$ hens\'elien et  $\varphi$  \'etale. Alors $\varphi_{A,\Top}$ est un hom\'eomorphisme local; plus pr\'ecis\'ement, pour tout $x\in\cX(A)$ et tout $J\in\FI$, $\varphi_{A,\Top}$  induit un ho\-m\'eo\-mor\-phisme de $B_{\cX}(x,J)$ sur $B_{\cY}(f(x),J)$.\cqfd
\end{numlist}
\end{fsprop}
\begin{fssrem} (non utilis\'ee dans la suite)

Avec les notations de \rref{not:top-entiere}, il est clair que l'application $\rho_{\cX}$ se factorise par $\cX(\wh{A})$, puisque $A/J\flis\wh{A}/J\wh{A}$ pour tout $J\in\FI$. En fait, on voit facilement que l'ensemble  $\varprojlim\limits_{J\in\FI}\cX(A/J)$ \emph{s'identifie} \`a $\cX(\wh{A})$ si $\wh{A}$ est hens\'elien, ou de fa\-\c{c}on \'equivalente si chacun des quotients $A/J$ ($J\in\FI$) est hens\'elien; ce sera  le cas notamment si $A$ est hens\'elien, ou si $\Gamma$ est de rang $1$ (c'est-\`a-dire si $\dim A=1$).
\end{fssrem}
\subsubsection{Comparaison des topologies sur les $A$-points et les $F$-points. }\label{ssec:compartop}
Gardons les hypoth\`eses et notations de \ref{not:top-entiere}; en outre, notons $j:\Spec{(F)}\inj\Spec{(A)}$ le morphisme canonique, et consid\'erons l'application naturelle
$$\begin{array}{rcl}
j^*:\cX_{A,\Top}&\ffl&\cX_{F,\Top}\\
x&\mapsto&x\circ j. 
\end{array}$$
Rappelons que $\cX_{A,\Top}$ est l'ensemble $\cX(A)$ muni de la topologie de \rref{not:top-entiere}; par ailleurs, $\cX_{F,\Top}$ d\'esigne l'ensemble $\cX(F)$, identifi\'e \`a $\cX_{F}(F)$ et muni \`a ce titre  de  la topologie du \S\,\rref{sec:espalg} (puisque $\cX_{F}=\cX\otimes_{A}F$ est un $F$-espace alg\'ebrique de type fini).
\begin{fsprop}\label{EntRat} Avec les notations de \rref{ssec:compartop}, on suppose  en outre que $A$ est \emph{hens\'elien}, 
ou bien que $\cX$ est un \emph{sch\'ema}. Alors:
\begin{numlist}
\item\label{EntRat0} $j^*$ est continue et ouverte.
\item\label{EntRat1} Si $\cX_{F}$ est localement s\'epar\'e, $j^*$ est un hom\'eomorphisme local.
\item\label{EntRat2} Si $\cX$ est localement s\'epar\'e, alors, pour tout $x\in\cX(A)$, $j^*$ induit un hom\'eomorphisme de $B_{\cX}(x,\Fm)$ sur un ouvert de $\cX(F)$.
\item\label{EntRat3} Si $\cX$ est s\'epar\'e (resp. propre sur $A$), alors $j^*$ est un plongement topologique ouvert (resp. un hom\'eomorphisme).\cqfd
\end{numlist}
\end{fsprop}
Donnons quelques indications sur la preuve de \ref{EntRat}.

\begin{trivlist}
\item[(a)] On commence par montrer que si $\cX$ est localement s\'epar\'e, alors (sans hy\-po\-th\`ese hens\'elienne sur $A$) la restriction de $j^*$ \`a $B_{\cX}(x,\Fm)$ est injective (pour $y$ et $y'$ ayant m\^eme image, consid\'erer le noyau de $\Spec\,(A)\underset{y}{\overset{y'}\rightrightarrows}\cX$). 
\item[(b)] Il en r\'esulte que pour $\cX$ localement s\'epar\'e, \ref{EntRat1} et \ref{EntRat2} sont \'equivalentes. 
\item[(c)] On montre  \ref{EntRat2} lorsque $\cX$ est un sch\'ema par r\'eduction au cas affine.
\item[(d)] Pour $A$ hens\'elien et $\cX_{F}$ localement s\'epar\'e, il y a un $A$-sch\'ema affine $\cX'$ et un $A$-morphisme \'etale $\pi:\cX'\to\cX$ tels que $x$ se rel\`eve en $x'\in\cX'(A)$. On montre alors \ref{EntRat1} en utilisant \ref{SoritesTopEnt}\,\ref{SoritesTopEnt4} (appliqu\'e \`a $\pi$),
\ref{prop:topLisseEtale}\,\ref{prop:topLisseEtale2} (appliqu\'e \`a $\pi_{F}$), et l'assertion \ref{EntRat2} pour les sch\'emas (appliqu\'ee \`a $\cX'$). 
\item[(e)] L'assertion \ref{EntRat3} est cons\'equence \'evidente de \ref{EntRat1} et du fait que, pour $\cX$ s\'epar\'e (reps. propre), $j^*$ est injective (resp. bijective).
\item[(f)] Enfin,  \ref{EntRat0} est imm\'ediate (et r\'esulte d'ailleurs de  \ref{EntRat1}) si $\cX$ est un sch\'ema affine; pour $A$ hens\'elien, on se ram\`ene \`a ce cas en utilisant \ref{SoritesTopEnt}\,\ref{SoritesTopEnt4} comme en (d).
\end{trivlist}
\Subsection{Torseurs sous un groupe lisse}\label{tors-lisse}

\begin{fsprop}\label{prop:tors-groupe-lisse}
Soit $F$ un corps  \emph{topologiquement hens\'elien}. Soient $G$ un $F$-groupe alg\'ebrique 
\emph{lisse}, $Y$ un $F$-espace, et $f:X\to Y$ un $G$-torseur 
au-dessus de $Y$ (de sorte que $X$ est aussi un $F$-espace). 

Alors $f_{\Top}: X_{\Top}\to Y_{\Top}$ est ouverte, et son image $I\subset Y_{\Top}$ est ouverte et ferm\'ee. 

Si de plus $Y$ est localement s\'epar\'e, l'application induite $X_{\Top}\to I$ est  une $G_{\Top}$-fibration principale.
\end{fsprop}
\dem comme $f$ est lisse, $f_{\Top}$ est  ouverte 
(\ref{prop:topLisseEtale}\,\ref{prop:topLisseEtale1}), et son image est ouverte. 
Cette image est l'ensemble des $y\in Y(F)$ tels que le $G$-torseur $X_{y}:=f^{-1}(y)$ sur
 $\Spec{(F)}$ soit trivial. Elle est donc ferm\'ee en vertu du lemme \ref{lem:tors-groupe-lisse} ci-dessous.

La derni\`ere assertion vient du fait que $X_{\Top}\to Y_{\Top}$ est automatiquement 
un pseudo-torseur sous $G_{\Top}$ par fonctorialit\'e (c'est-\`a-dire que l'appli\-ca\-tion $G_{\Top}\times X_{\Top}\to X_{\Top}\times_{Y_{\Top}} X_{\Top}$ donn\'ee par $(g,x)\mapsto(gx,x)$ est un ho\-m\'eo\-mor\-phisme; on utilise ici la compatibilit\'e aux produits fibr\'es), et que $f_{\Top}$ admet des sections 
locales au voisinage de tout point de $I$ (\ref{prop:topLisseEtale}\,\ref{prop:topLisseEtale3}). \cqfd

\begin{fsslem}\label{lem:tors-groupe-lisse} Sous les hypoth\`eses de 
la proposition \rref{prop:tors-groupe-lisse}, l'application
$$\begin{array}{rcl}[f]:Y_{\Top}&\ffl&\mathrm{H}^1_{\mathrm{et}}(F,G)\\
y&\longmapsto & [X_{y}]
\end{array}
$$
est localement constante.
\end{fsslem} 
\dem fixons une classe $\xi\in\mathrm{H}^1_{\mathrm{et}}(F,G)$, et soit $T$ un $G$-torseur 
sur $\Spec{(F)}$ de classe $\xi$. 
Consid\'erons le $G$-torseur \og constant\fg\ $T_{Y}= {T\times_{F} Y}$ au-dessus de $Y$. 
Le faisceau \'etale $I:=\ul{\mathrm{Isom}}_{G_{Y}}(T_{Y},X)$ des morphismes (de $G$-torseurs sur $Y$) 
de $T_{Y}$ vers $X$ est un torseur sous le groupe des $G$-automor\-phismes de $T$; ce groupe est
 une forme int\'erieure de $G$ et donc un $F$-groupe lisse, de sorte que l'image de $I(F)$ 
dans $Y_{\Top}$ est ouverte. Or celle-ci n'est autre que $[f]^{-1}(\xi)$, d'o\`u la conclusion.\cqfd

\begin{fsprop}\label{prop:tors-strict}
Soit $F$ un corps  \emph{topologiquement hens\'elien}. Soient $G$ un $F$-groupe 
alg\'ebrique, $Y$ un $F$-espace, et $f:X\to Y$ un $G$-torseur au-dessus de $Y$. 

On note $I \subset Y_\Top$ l'image de $f_\Top$. 
Consid\'erons les propri\'et\'es suivantes:
\begin{romlist}
\item\label{prop:tors-strict1} l'application induite $X_{\Top}\to I$ est une $G_{\Top}$-fibration principale;
\item\label{prop:tors-strict2} l'application induite $X_{\Top}\to I$ est ouverte;
\item\label{prop:tors-strict3} l'application induite $X_{\Top}/G(F)\to I$ est un hom\'eomorphisme;
\item\label{prop:tors-strict4} $f_{\Top}:X_{\Top}\to Y_{\Top}$ est stricte \emph{(\ref{conventions})}.
\end{romlist}
Alors on a les implications \ref{prop:tors-strict1}$\Rightarrow$\ref{prop:tors-strict2}$\Leftrightarrow$\ref{prop:tors-strict3}$\Leftrightarrow$\ref{prop:tors-strict4},
 et les quatre conditions 
sont \'equivalentes si $Y$ est localement s\'epar\'e.
\end{fsprop} 
\begin{fssrem}\label{rem:prop:tors-strict} Lorsque $F$ est un corps valu\'e admissible, cette proposition 
est une cons\'equence triviale du th\'eor\`eme principal \ref{ThPpal}: les propri\'et\'es \ref{prop:tors-strict2}  
\`a \ref{prop:tors-strict4} sont vraies, et \ref{prop:tors-strict1} l'est si $Y$ est localement s\'epar\'e. 
Comme le th\'eor\`eme  \ref{ThPpal} sera d\'emontr\'e sans utiliser \ref{prop:tors-strict}, cette proposition est 
donc redondante dans ce cas. Ceci \'etant, elle donne 
un r\'esultat  partiel dans la direction de \ref{ThPpal} sans faire intervenir les compactifications. 
\end{fssrem}
\noindent\textsl{D\'emonstration de la proposition:} les implications
 \ref{prop:tors-strict1}$\Rightarrow$\ref{prop:tors-strict2}$\Rightarrow$\ref{prop:tors-strict3}$\Leftrightarrow$\ref{prop:tors-strict4} sont 
 imm\'ediates; 
en outre  \ref{prop:tors-strict3}$\Rightarrow$\ref{prop:tors-strict2} vient du fait que la relation d'\'equivalence d\'eduite 
d'une action continue de groupe est toujours ouverte. 

Supposons maintenant \ref{prop:tors-strict2} v\'erifi\'ee et  $Y$ localement s\'epar\'e, et montrons \ref{prop:tors-strict1}. Posons 
$X'=X/G^\bec$, o\`u $G^\bec$ est d\'efini en \ref{ssec:pgsgl}; noter en particulier que $G^\bec_{\Top}=G_{\Top}$. Alors  
$f$ se factorise suivant 
$$
f : X \buildrel f' \over \lgr  X' \buildrel g \over \lgr  Y.
$$ 
Notons $I'\subset X'_{\Top}$ l'image de $f'_{\Top}$. On a un diagramme d'espaces topologiques
$$X_{\Top}\buildrel f'_{1} \over \ssurj I' \buildrel g_{1} \over \ssurj  I\inj Y_{\Top}$$
o\`u $f'_{1}$ est surjective par d\'efinition de $I'$,  et $g_{1}$ \emph{bijective} par  \ref{quotient-torseur}. 
Par hypoth\`ese  la compos\'ee $g_{1}\circ f'_{1}$ est ouverte. On en d\'eduit  que  $g_{1}$ est ouverte, donc est un hom\'eomorphisme
de $I'$ sur $I$. 
Il suffit donc de voir que $f'_{1}$ est une $G^\bec_{\Top}$-fibration principale. Comme $f'$ est un $G^\bec$-torseur et 
que $G^\bec$ est lisse, il suffit pour cela d'appliquer  \ref{prop:tors-groupe-lisse} en remarquant que $X'$ est 
\emph{localement s\'epar\'e}: 
cette derni\`ere propri\'et\'e r\'esulte de l'hypoth\`ese sur $Y$ et du fait que $g:X'\to Y$ est un morphisme s\'epar\'e 
(il est localement isomorphe,
pour la topologie fppf sur sur $Y$, \`a la projection $Y\times (G/G^\bec)\to Y$).\cqfd 

\Subsection{Corps valu\'es hens\'eliens: approximation faible et applications}
\subsubsection{Rappels sur les corps valu\'es hens\'eliens. }\label{rapp:hens}\ 

Soit $(K,v)$ un corps valu\'e, d'anneau de valuation $A$. Rappelons que pour que $(K,v)$ soit hens\'elien 
(c'est-\`a-dire pour que $A$ soit un anneau local hens\'elien) il faut et il suffit que pour toute extension finie $L$ de $K$,
la valuation $v$ se prolonge de mani\`ere unique en une valuation de $L$ \cite[th. 32.8 p. 348]{Wa}. Dans ce cas, on supposera toujours $L$ muni de cette valuation,
et de la topologie associ\'ee. 

On prendra garde qu'un corps valu\'e complet n'est pas n\'ecessairement hen\-s\'e\-lien; c'est toutefois vrai si la valuation est de rang $1$. 

Supposons $(K,v)$  hens\'elien; notons $\wh{K}$ son compl\'et\'e, et  $K_s$ une cl\^oture s\'eparable de $K$.
Alors $\wh{K}$ est hens\'elien (c'est imm\'ediat) et de plus  $K$ est s\'eparablement ferm\'e dans 
$\widehat K$ \cite[th. 32.19 p. 357]{Wa}, de sorte que   $\wh{K} \otimes_K  K_s$ est un corps. 
Le corps $\wh{K} \otimes_K  K_s$ est m\^eme une  \emph{cl\^oture s\'eparable}
de $\wh{K}$. En effet, en vertu du th\'eor\`eme de l'\'el\'ement primitif, une extension
finie s\'eparable  de $\widehat K$ est
isomorphe \`a $\widehat K[t]/P(t)$  o\`u
$P$ d\'esigne un $\widehat K$-polyn\^ome unitaire s\'eparable; si  $Q$ est   un $K$-polyn\^ome unitaire (s\'eparable)
assez proche de  $P$, alors  $P(T)$ et  $Q(T)$ ont m\^eme 
m\^eme corps de d\'ecomposition sur $\widehat K$ \cite[32.20]{Wa},  donc $P$ se d\'ecompose sur  $\wh{K} \otimes_K  K_s$. 
On a montr\'e que $\wh{K} \otimes_K  K_s$ est une cl\^oture s\'eparable de $\widehat K$, 
 on a donc un isomorphisme   $\Gal( \wh{K} \otimes_K K_s/ \wh{K}) \simlgr \Gal( K_s/K)$. 

\begin{fsprop}\label{hens-app-faible} \textup{(approximation faible)} Soient  $F$ un corps   valu\'e hen\-s\'e\-lien et $X$ un $F$-espace. 
On suppose que $X$ est  \emph{lisse} ou que $F$ est \emph{admissible (\ref{def:admissible})}. Alors $X(F)$ est dense dans $X(\widehat F)$. 
\end{fsprop}

\dem dans les deux cas, les \'enonc\'es \ref{cor:FiltrationCLO} et \ref{prop:caractTopEspAlg} permettent de se ramener au cas o\`u $X$ est une \emph{vari\'et\'e}.

Lorsque $F$ est admissible, le r\'esultat est alors d\'emontr\'e dans \cite[corollaire 1.2.1]{MB2} (comme cons\'equence du th\'eor\`eme d'approximation fort, voir \ref{th:Greenberg} plus bas). 

Supposons d\'esormais $X$ lisse. On peut en outre supposer $X$ affine et pu\-re\-ment de dimension $d \geq 1$. 
Soient $\Omega$ un ouvert non vide de $X(\widehat F)_\Top$  et $x \in \Omega$. Rempla\c{c}ant $X$ par un voisinage affine de $x$, on peut supposer qu'il existe un morphisme \'etale $f: X \to {\bf A}^d_{F}$. Comme $\widehat F$ est hens\'elien, la pro\-po\-si\-tion \ref{Hens=hens} montre
que  l'application   $f_{\Top}:X(\widehat F)_{\Top} \to \wh{F}^d$ est 
ouverte. En particulier $f_{\Top}(\Omega)$ est ouvert dans $\wh{F}^d$, donc rencontre  $F^d$ car $F$ est dense dans $\wh{F}$. Autrement dit, $\Omega$ rencontre $f_{\Top}^{-1}(F^d)$. Mais 
comme $F$ est s\'eparablement ferm\'e  dans $\widehat F$ et $f$ \'etale, on a
$f_{\Top}^{-1}(F^d)= X(F)$, d'o\`u $X(F) \cap \Omega \not =\emptyset$.
\cqfd

\begin{fsprop}\label{krasner} Soit $(K,v)$ un corps valu\'e hens\'elien, de compl\'et\'e $\wh{K}$. 
 Soit $G$ un $K$-sch\'ema en groupes localement de type fini. On notera
 $$\rho^i:\rH^i(K,G) \to \mathrm{H}^i(\wh{K},G)$$
 l'application naturelle, d\'efinie pour $i\in\{0,1\}$ en g\'en\'eral et pour $i\in\NN$ si $G$ est commutatif. 

 \begin{numlist}
\item\label{krasner1}  Si $K$ est admissible, $\rho^1$ est injective. 
\item\label{krasner2}   Si $G$ est  lisse, $\rho^1$  est bijective. 
\item\label{krasner4}   Si $G$ est commutatif, $\rho^i$ est bijective pour tout $i\geq2$.
\end{numlist}
\end{fsprop}

\dem \ref{krasner1} L'argument de torsion habituel nous ram\`ene \`a \'etablir la 
trivialit\'e du noyau de  l'application $H^1(K,G) \to H^1(\wh{K},G)$. 
Soit donc $E$ un $G$-torseur tel que $E(\wh{K}) \not = \emptyset$. Alors $E$ est 
automatiquement un $K$-espace alg\'ebrique s\'epar\'e et localement de type fini; 
il admet donc un sous-espace ouvert $U$, de type fini sur $K$ et v\'erifiant $U(\wh{K})\neq\emptyset$. 
Le th\'eor\`eme d'approximation faible \ref{hens-app-faible} montre que $E(K)\neq\emptyset$. On conclut que  $E/K$ est le $G$-torseur trivial.
\smallskip

\noindent \ref{krasner2}  Comme  $G$ est lisse sur $K$, les cohomologies \'etale et fppf de $G$ co\"{\i}ncident. D'autre part tout 
$G$-torseur est lisse sur $K$, de sorte que l'injectivit\'e se d\'emontre comme en \ref{krasner1} en utilisant le cas lisse de \ref{hens-app-faible}. 
Pour la surjectivit\'e,  on se limite dans  un premier temps au cas o\`u $G$ est affine.
On plonge alors $G$ dans
 un $K$-groupe lin\'eaire  $\GL(V)$ et on pose $X= \GL(V)/G$.
Suivant le th\'eor\`eme 90 de Hilbert, $\mathrm{H}^1_{\mathrm{et}}(K,\GL(V))=1$ et   
$\mathrm{H}^1_{\mathrm{et}}(\wh{K},\GL(V))=1$. La suite exacte courte 
de cohomologie galoisienne \cite[I.5.4]{S}  
produit donc un carr\'e commutatif exact d'ensembles point\'es 
$$ 
\begin{CD}
X(K) @>>> \mathrm{H}^1_{\mathrm{et}}(K,G) @>>> 1\\
@VVV @VVV \\
X(\wh{K}) @>>> \mathrm{H}^1_{\mathrm{et}}(\wh{K},G) @>>> 1 .
\end{CD}
$$ 
On a vu que l'application  $X(\wh{K})_\Top \to  \mathrm{H}^1_{\mathrm{et}}(\wh{K},G)$ est localement constante (lemme \ref{lem:tors-groupe-lisse}).
Comme $X(K)$ est dense dans $X(\wh{K})_\Top$  d'apr\`es \ref{hens-app-faible}, il suit que
le compos\'e $X(K) \to X(\wh{K}) \to  \mathrm{H}^1_{\mathrm{et}}(\wh{K},G) $ est surjectif. Le diagramme ci-dessus permet de conclure que
l'application $\mathrm{H}^1(K,G) \to \mathrm{H}^1_{\mathrm{et}}(\wh{K},G)$ est 
surjective, donc bijective.
\smallskip

Il reste \`a montrer la surjectivit\'e de
la restriction $\mathrm{H}^1_{\mathrm{et}}(K,G) \to \mathrm{H}^1_{\mathrm{et}}(\wh{K},G)$
pour $G/K$ lisse quelconque. Soit $\wh{\gamma} \in \mathrm{H}^1_{\mathrm{et}}(\wh{K},G)$, et soit $L$ une extension finie s\'eparable de $\wh{K}$ qui trivialise $\wh{\gamma}$. D'apr\`es \ref{rapp:hens}, il existe une extension finie s\'eparable $K'$ de $K$ telle que $L= \wh{K} \otimes_K K'$, qui est aussi le compl\'et\'e $\wh{K'}$ de $K'$.

Le quotient $X=\bigl(\prod\limits_{K'/K}G \bigr)/G$  est repr\'esentable
par un $K$-sch\'ema \cite[VI$_A$.3.2]{SGA3}. On observe que $X$ est une 
$K$-forme de $G^{[K':K]-1}$ donc est lisse.  On a le diagramme commutatif d'ensembles point\'es 
$$
\begin{CD}
X(K) @>>> \mathrm{H}^1_{\mathrm{et}}(K,G) @>>> \mathrm{H}^1_{\mathrm{et}}\bigl(K , \prod\limits_{K'/K} G \bigr) &\, \simla \, & \mathrm{H}_{\mathrm{et}}^1(K',G)  \\
@VVV @VVV @VVV @VVV \\
X(\wh{K}) @>>> \mathrm{H}^1_{\mathrm{et}}(\wh{K},G)
 @>>> \mathrm{H}^1_{\mathrm{et}}\bigl(\wh{K} , \prod\limits_{K'/K} G \bigr)
&\, \simla \, & \mathrm{H}^1_{\mathrm{et}}(\wh{K'},G) , 
\end{CD} 
$$  
o\`u dans chaque ligne la suite des deux premi\`eres fl\`eches est exacte,
 et o\`u les bijections de droite sont  celles  \og de Shapiro\fg\ \cite[XXIV.8.4]{SGA3}.
L'application $X(\wh{K})_\Top \to  H^1_{\mathrm{et}}(\wh{K},G)$ est localement constante et 
$X(K)$ est dense dans $X(\wh{K})_\Top$, 
donc $X(K)$ se surjecte sur $\mathrm{H}^1(\wh{K'}/ \wh{K},G)$, qui contient $\wh{\gamma}$ vu le choix de $K'$. 
Le diagramme ci-dessus montre que $\wh{\gamma}$ provient de $\mathrm{H}^1_{\mathrm{et}}(K,G)$. 
\smallskip

\noindent \ref{krasner4} \noindent{\it Premier cas: $G$ est lisse.} 
On va utiliser l'argument classique de d\'ecalage par r\'ecurrence sur $i \geq 2$ en notant 
que  le cas $i=1$  a \'et\'e trait\'e.  Pour tout corps $F/K$, 
on a $\mathrm{H}^i_{\mathrm{et}}(F,G) \simlgr \mathrm{H}^i_{\mathrm{fppf}}(F,G)$
pour tout $i \geq 0$ (voir \cite[11.7]{Brauer3} ou \cite[III.3.9 et 3.11.(b)]{Mi}).   
En particulier $\mathrm{H}^i(K,G)$ co\"{\i}ncide avec  la cohomologie galoisienne $\mathrm{H}^i(\Gal(K_s/K),G(K_s))$ et commute
aux limites inductives filtrantes de $K$-corps.  
Soit donc $K'/K$ une extension finie s\'eparable et consid\'erons le 
$K$-groupe quotient (lisse) $H=\bigl(\prod\limits_{K'/K}G \bigr)/G$.
La suite exacte de $K$-groupes $0 \to G \to \prod\limits_{K'/K}G \to H \to 0$
donne lieu au diagramme commutatif exact
{\small
$$
\begin{CD}
 \mathrm{H}^{i-1}_{\mathrm{et}}(K, \prod\limits_{K'/K} G) @>>> \mathrm{H}_{\mathrm{et}}^{i-1}(K,H)
@>>> \mathrm{H}^i_{\mathrm{et}}(K,G) @>>> \mathrm{H}^i_{\mathrm{et}}(K,\prod\limits_{K'/K}G) \\
@VVV @VV{\wr}V @VVV @VVV \\
 \mathrm{H}^{i-1}_{\mathrm{et}}(\wh{K}, \prod\limits_{K'/K} G) @>>> \mathrm{H}_{\mathrm{et}}^{i-1}(\wh{K},H)
@>>> \mathrm{H}^i_{\mathrm{et}}(\wh{K},G) @>>> \mathrm{H}^i_{\mathrm{et}}(\wh{K},\prod\limits_{K'/K}G). \\
\end{CD}
$$
}
La seconde fl\`eche verticale (en partant de la   gauche) est bijective suivant l'hypoth\`ese de r\'ecurrence.
La premi\`ere l'est \'egalement en tenant compte du \og lemme de Shapiro\fg\ \cite[I.2.5]{S}.   
Par suite l'application $$
\ker\bigl( \mathrm{H}^i_{\mathrm{et}}(K,G) \to  \mathrm{H}^i_{\mathrm{et}}(K',G) \bigr)
\, \to 
\ker\bigl( \mathrm{H}^i_{\mathrm{et}}(\wh{K},G) \to  \mathrm{H}^i_{\mathrm{et}}(\wh{K} \otimes_K K',G) \bigr)
$$ 
est bijective. En passant \`a limite sur les sous-extensions finies $K'/K$ de $K_s$, on 
obtient que  $\mathrm{H}^i_{\mathrm{et}}(K,G) \simlgr 
\ker\bigl( \mathrm{H}^i_{\mathrm{et}}(\wh{K},G) \to  
\mathrm{H}^i_{\mathrm{et}}(\wh{K} \otimes_K K_s,G) \bigr)$.
Suivant \ref{rapp:hens}, le corps  $\widehat K \otimes_K K_s$ est s\'eparablement clos
d'o\`u la bijection souhait\'ee $\mathrm{H}^i_{\mathrm{et}}( K,G) \simlgr
  \mathrm{H}^i_{\mathrm{et}}(\wh{K},G)$. 
\smallskip

\noindent{\it Second cas: $G$ est fini.} Soit $A$ l'alg\`ebre affine du dual de Cartier de $G$. Alors $G$ se plonge canoniquement dans le groupe des unit\'es $U$ de $A$, qui est un ouvert de Zariski du $K$-espace affine sous-jacent \`a $A$ est est donc lisse. On a donc une suite exacte $1\to G\to U\to U/G\to 1$. Les groupes $U$ et $U/G$ \'etant lisses commutatifs, ils sont justiciables du cas pr\'ec\'edent et l'on conclut par d\'evissage.
\smallskip

\noindent{\it Cas g\'en\'eral.} En caract\'eristique nulle, 
$G$ est lisse et le premier cas suffit. Si $K$ est de caract\'eristique $p>0$, 
il existe un entier $n \geq 1$ tel que le quotient de $G'= G/ {_{\mathrm{Fr}^n}G}$ par 
le noyau du morphisme de Frobenius it\'er\'e soit lisse \cite[VII$_A$.8.3]{SGA3}. 
Ce noyau $N= {_{\mathrm{Fr}^n}G}$ est fini sur $K$.
En tenant compte des cas pr\'ec\'edents, la suite exacte
$0 \to N \to G \to G' \to 1$ donne donc le r\'esultat par d\'evissage.
\cqfd

\begin{fsrems} 
\noindent (a) 
L'assertion \ref{krasner1} est fausse si $K$ est hens\'elien mais non admissible. 
Plus g\'en\'eralement, soit $E\subset F$ une extension de corps \emph{non s\'eparable}, 
de caract\'eristique $p$. Alors il existe 
une suite  $(a_1,\dots, a_n)$ dans $E$, libre sur $E^p$
mais  li\'ee sur $F^p$, avec $n$ minimal. Consid\'erons le $E$-groupe alg\'ebrique $G$ 
noyau de l'homomorphisme $f:(x_1, \dots, x_{n-1}) \mapsto \sum\limits_{i=1}^{n-1} x_i^p \, a_i$ de $ \GG_a^{n-1}$ dans $\GG_a$.
Alors l'\'equation $f(x_1,\dots,x_{n-1})=a_n$ d\'efinit un $G$-torseur non trivial qui est trivialis\'e par $F$.
Ainsi l'application naturelle $\rH^1(E,G) \to \mathrm{H}^1(F,G)$ n'est pas injective.
\smallskip

\noindent (b) L'assertion \ref{krasner2} est fausse en g\'en\'eral pour un $K$-groupe non lisse
comme l'indique le contre-exemple suivant. On note $A$ l'hens\'elis\'e de $\FF_p[t]$ en 
$0$ et $K$ son corps de fractions. Alors le compl\'et\'e de $K$ est $\FF_p(\!(t)\!)$.
La fl\`eche $K/K^p \to \wh{K}/(\wh{K})^p$ n'est 
pas surjective car le membre de gauche est d\'enombrable alors
que le membre de droite ne l'est pas. 
\smallskip

\noindent (c) Lorsque $G$ est lisse ou $K$ admissible (cas \ref{krasner1} et \ref{krasner2}), le th\'eor\`eme d'approximation faible \ref{hens-app-faible}, 
appliqu\'e \`a $G$, peut \^etre vu comme le cas $i=0$: l'application $\rho^0:G(K)\to G(\wh{K})$ (\'evidemment injective) est d'image dense
pour la topologie naturelle sur $G(\wh{K})$.
\end{fsrems}

\section{Corps valu\'es admissibles; le cas d'un groupe $G$ tel 
que $G^\circ_{\red}$ soit lisse}\label{sec:admissible}

\Subsection{Corps valu\'es admissibles: g\'en\'eralit\'es}\label{subsec:genadmissible}

Donnons d'abord  quelques propri\'et\'es \'el\'ementaires des corps valu\'es admissibles, d\'efinis en \ref{def:admissible}; elles ne seront
pas utilis\'ees ici, et nous en laissons au lecteur les d\'etails des d\'emonstrations:

\begin{fsprop}\label{prop:pptesElemAdm}
Soit $(K,v)$ un corps valu\'e admissible, de compl\'et\'e $\wh{K}$ et d'exposant caract\'eristique $p$. 
\begin{numlist}
\item $K$ est alg\'ebriquement ferm\'e dans $\wh{K}$.
\item L'endomorphisme $x\mapsto x^p$ de $K$ est ferm\'e.
\item\label{prop:pptesElemAdm3} Si $L$ est une extension finie de $K$ \emph{(munie de sa valuation prolongeant $v$, cf. \ref{rapp:hens})}, alors:
\begin{subromlist}
\item\label{prop:pptesElemAdm3a} $L$ est admissible;
\item\label{prop:pptesElemAdm3b} l'homomorphisme canonique $\wh{K}\otimes_{K}L\to\wh{L}$ est un isomorphisme;
\item\label{prop:pptesElemAdm3c} comme $K$-espace vectoriel topologique, $L$ est topologiquement libre 
(c'est-\`a-dire isomorphe \`a $K^{\oplus[L:K]}$); en particulier tout sous-espace vectoriel de $L$ (et notamment $K$) est ferm\'e dans $L$.
\end{subromlist}
\end{numlist}
\end{fsprop}

Donnons seulement quelques indications sur la preuve de \ref{prop:pptesElemAdm3}. On sait d'apr\`es \cite[VI.8.2]{BAC} que 
l'homomorphisme de \ref{prop:pptesElemAdm3b} est surjectif et a pour noyau le nilradical de $\wh{K}\otimes_{K}L$; puisque $\wh{K}/K$ est s\'eparable,
\ref{prop:pptesElemAdm3b} en r\'esulte (et entra\^{\i}ne imm\'ediatement \ref{prop:pptesElemAdm3a}). Pour montrer \ref{prop:pptesElemAdm3c}, 
rappelons d'abord que \ref{prop:pptesElemAdm3c} est vrai si $K$ est \emph{complet}: sur un corps valu\'e complet non discret, tout espace  
vectoriel topologique s\'epar\'e de dimension finie est libre \cite[VI.5.2, prop. 4]{BAC}. 
Soit maintenant $\varphi:L\to K$ une forme $K$-lin\'eaire, et montrons que $\varphi$ est continue 
(ce qui suffit pour montrer \ref{prop:pptesElemAdm3c}). L'assertion \ref{prop:pptesElemAdm3b} 
implique que $\varphi$ se prolonge en une forme $\wh{K}$-lin\'eaire $\varphi_{\wh{K}}:\wh{L}\to \wh{K}$; 
celle-ci est automatiquement continue d'apr\`es le cas complet, donc $\varphi$ l'est aussi par restriction.\cqfd
\Subsection{Le th\'eor\`eme d'approximation fort; applications}\label{subsec:approx}
\subsubsection{Notations. }\label{propre:notations} On d\'esigne par $(K,v)$ un corps valu\'e \emph{admissible}, par $A$ l'anneau de $v$, par $k$ son corps r\'esiduel et par $\wh{K}$ le compl\'et\'e de $K$.

\begin{fsteo}\label{th:Greenberg} \textup{(\og th\'eor\`eme d'approximation fort\fg)} Sous les hy\-po\-th\`eses de  \textup{\ref{propre:notations}}, soit  $\cZ$ un $A$-espace alg\'ebrique de pr\'esentation finie. Pour tout id\'eal non nul $J$ de $A$, il existe un id\'eal non nul $J'\subset J$ tel que $\cZ(A)$ et $\cZ(A/J')$ aient m\^eme image dans $\cZ(A/J)$.
\end{fsteo}
\dem lorsque $\cZ$ est un sch\'ema, le th\'eor\`eme est  \'etabli dans \cite{Gre} pour les valuations discr\`etes et dans \cite{MB2} pour le cas g\'en\'eral.

Sinon, il existe d'apr\`es \ref{cor:FiltrationCLO} un morphisme \'etale $\pi:\cZ'\to \cZ$, o\`u $\cZ'$ 
est un sch\'ema affine, tel que (notamment) l'application induite $\pi_{k}:\cZ'(k)\to\cZ(k)$ soit surjective. 

Or, comme $\pi$ est \'etale et que $A$ (et donc tout quotient de $A$) est hens\'elien, le syst\`eme projectif 
 $(\cZ'(A/J))_{J}$ (o\`u $J$ parcourt les id\'eaux stricts de $A$) se d\'eduit du syst\`eme projectif  $(\cZ(A/J))_{J}$ 
par le changement de base surjectif $\pi_{k}$. L'\'enonc\'e est  une cons\'equence formelle de  cette remarque 
et du cas des sch\'emas (appliqu\'e \`a $\cZ'$).\cqfd

\begin{fssrem} Dans l'\'enonc\'e de \ref{th:Greenberg}, on peut prendre $J'$ de la forme $cJ^n$, o\`u $n\in \NN_{>0}$
et $c\in A\smallsetminus\{0\}$ sont ind\'ependants de $J$; nous n'aurons pas \`a utiliser ce fait.
\end{fssrem}
Le th\'eor\`eme \ref{th:Greenberg} a la cons\'equence suivante, essentielle pour la preuve de \ref{ThPpal}:
\begin{fsteo}\label{th:propre} Sous les hypoth\`eses de \textup{\ref{propre:notations}}, soit $f:X\to Y$ un morphisme
 \emph{propre} de $K$-espaces, et soit  $y$ un point de $Y(K)$. On suppose que la fibre 
$C:=f_{\Top}^{-1}(y)\subset X_{\Top}$ est \emph{compacte} \textup{(observer qu'elle est a priori 
s\'epar\'ee puisqu'elle s'identifie \`a $(X_{y})_{\Top}$, cf. \rref{rem:ProdFibEspAlg1}\,\rref{rem:ProdFibEspAlg13})}.

Alors  tout voisinage $U$ de $C$ dans $X_{\Top}$ contient un voisinage de la forme $f_{\Top}^{-1}(V)$, o\`u $V$ 
est un voisinage de $y$ dans $Y_{\Top}$.

De fa\c{c}on \'equivalente, si $\Phi\subset X_{\Top}$ est un ferm\'e disjoint de $C$, alors $y$ n'est pas adh\'erent 
\`a $f_{\Top}(\Phi)$.
\end{fsteo}
\dem commen\c{c}ons par le cas o\`u $Y$ est une $K$-vari\'et\'e, que l'on peut supposer  \emph{affine}, la question 
\'etant locale sur $Y_{\Top}$. Choisissons alors un $A$-sch\'ema affine $\cY$, de pr\'esentation finie, de fibre
 g\'en\'erique $Y$, tel que $y$ se prolonge en un point $\til{y}\in \cY(A)$. Alors $\cY(A)$ est un
voisinage de $y$ (identifi\'e \`a $\til{y}$) dans $Y(K)$; quitte \`a restreindre $U$, on peut supposer
 que $U\subset f_{\Top}^{-1}(\cY(A))$. 

\begin{fsslem}\label{Nagata} Il existe un $A$-espace $\cX$ de pr\'esentation finie, de fibre g\'en\'erique $X$, 
et un $A$-morphisme \emph{propre} $\varphi:\cX\to\cY$ tel que $\varphi_{K}:\cX_{K}\to\cY_{K}$ s'identifie \`a $f$. 
\end{fsslem}
\dem c'est essentiellement le  th\'eor\`eme de compactification de Nagata (voir \cite[Theorem 1.2.1]{CLO}, 
ou \cite[Theorem 4.1]{ConNagata} pour le cas des sch\'emas), mais il faut prendre garde que le morphisme 
naturel $X\to \cY$ n'est pas en g\'en\'eral de type fini. Il faut commencer par prolonger $X$ en un $\cY$-espace 
de pr\'esentation finie $\cX_{1}$, de fibre g\'en\'erique $X$:  il suffit ensuite d'appliquer le th\'eor\`eme 
cit\'e \`a $\cX_{1}\to \cY$. 

Lorsque $X$ est un sch\'ema, l'existence de $\cX_{1}$ r\'esulte des th\'eor\`emes g\'e\-n\'e\-raux de \cite[\S\:8]{EGA4}. 
Pour un espace alg\'ebrique, on \'ecrit $X$ comme quotient d'un $Y$-sch\'ema $X'$ par une $Y$-relation 
d'\'equivalence \'etale $R \rightrightarrows X$, et l'on applique les r\'esultats de {\it loc. cit.} au 
diagramme en question.\cqfd
\medskip

On fixe d\'esormais un \og $\cY$-mod\`ele propre\fg\ $\varphi:\cX\to\cY$ comme dans le lemme \ref{Nagata}.
Le crit\`ere valuatif de propret\'e implique alors que $f_{\Top}^{-1}(\cY(A))=\cX(A)$ (vus comme sous-ensembles de $X(K)$); 
nous savons en outre par \ref{EntRat}\,\ref{EntRat3} que les sous-espaces de $X_{\Top}$ et $Y_{\Top}$ ainsi d\'efinis
 co\"{\i}ncident avec les espaces $\cX_{A,\Top}$ et  $\cY_{A,\Top}$ d\'efinis en \ref{not:top-entiere}. Nous pouvons donc,
 dans l'\'enonc\'e, remplacer $X_{\Top}$ et $Y_{\Top}$ par $\cX_{A,\Top}$ et  $\cY_{A,\Top}$; en particulier, 
une base d'ouverts de $\cX_{A,\Top}$ est fournie par les boules  $B_{\cX}(\xi,J)$ ($\xi\in\cX(A)$, $J$ id\'eal non nul de $A$) 
d\'efinies en  \ref{not:top-entiere}. 

Nous noterons $\cX_{\til{y}}$ le produit fibr\'e $\cX\times_{\cY,\til{y}}\Spec{(A)}$: c'est un sous-espace ferm\'e de $\cX$ et
 un $A$-sch\'ema propre, et  $\cX_{\til{y}}(A)$ s'identifie \`a $C$. 
La compacit\'e de $C$ implique donc que, quitte \`a restreindre $U$, on peut supposer que celui-ci est, pour $J$ convenable, 
de la forme
$$\begin{array}{rcl} 
U&=&\bigcup\limits_{\xi\in C}B_{\cX}(\xi,J)\\
&=&\left\{\,\xi'\in\cX(A)\;\mid\; \exists\xi\in \cX_{\til{y}}(A), \;\xi\bmod J=\xi'\bmod J\,\right\}.
\end{array}$$
Le th\'eor\`eme d'approximation \ref{th:Greenberg}, appliqu\'e avec  $\cZ=\cX_{\til{y}}$,  fournit un id\'eal 
$J'\subset J$ de $A$ tel que $\cX_{\til{y}}(A/J')$ et $\cX_{\til{y}}(A)$ aient la m\^eme image dans $\cX_{\til{y}}(A/J)$. 

Posons alors $V:=B_{\cY}(\til{y},J')$; soit $\xi\in \varphi^{-1}(V)$, et v\'erifions que $\xi\in U$. Par hypoth\`ese,
 $\varphi(\xi\bmod J')=\til{y}\bmod J'$, de sorte que $\xi\bmod J'\in\cX_{\til{y}}(A/J')$. Vu le choix de $J'$, son image 
dans $\cX_{\til{y}}(A/J)$, qui est \'evidemment $\xi\bmod J$, se rel\`eve en un point de $\cX_{\til{y}}(A)$; autrement dit,
 on a $\xi\in U$, comme annonc\'e.

\medskip

Ne supposant plus que $Y$ est un sch\'ema, choisissons une $K$-vari\'et\'e affine $Y'$ et un morphisme \'etale $\pi:Y'\to Y$
 tel que $y$ se rel\`eve en $y'\in Y'(F)$. Consid\'erons le produit fibr\'e $X':=X\times_{Y}Y'$ et le diagramme correspondant
 d'espaces topologiques
$$\xymatrix{C'\,\ar@{^{(}->}[r]\ar[d]& X'_{\Top} \ar[r]^{\pi'_{\Top}}\ar[d]^{f'_{\Top}} & X_{\Top}\ar[d]^{f_{\Top}}\\
\{y'\}\, \ar@{^{(}->}[r] & Y' \ar[r]^{\pi_{\Top}} & Y
}
$$
o\`u l'on a pos\'e $C'=f_{\Top}^{\prime-1}(y')\subset X'$; on sait (\ref{rem:ProdFibEspAlg1}) que $\pi'_{\Top}$ induit un
 hom\'eomorphisme de $C'$ sur $C$, de sorte que $C'$ est compact. 

Posons $U'=\pi^{\prime-1}_{\Top}(U)$: c'est un voisinage de $C'$ dans $X'_{\Top}$. Il contient donc, d'apr\`es le cas d\'ej\`a \'etabli, un voisinage de la forme $f_{\Top}^{\prime-1}(W')$, o\`u $W'\subset Y'_{\Top}$ est un ouvert contenant $y'$. On a donc
$$C=\pi'_{\Top}(C')\subset \pi'_{\Top}(f_{\Top}^{\prime-1}(W'))\subset \pi'_{\Top}(U')\subset U.$$
Bien que le diagramme ci-dessus ne soit pas n\'ecessairement cart\'esien, le diagramme sous-jacent \emph{d'ensembles} l'est,
 de sorte que 
$$\pi'_{\Top}(f_{\Top}^{\prime-1}(W'))=f_{\Top}^{-1}(\pi_{\Top}(W')).$$
Mais $\pi_{\Top}$ est ouverte puisque $\pi$ est \'etale (\ref{prop:topLisseEtale}) donc $V:=\pi_{\Top}(W')$ est un voisinage 
ouvert de $y$ dans $Y_{\Top}$, et les relations qui pr\'ec\`edent montrent que $f_{\Top}^{-1}(V)\subset U$, ce qui ach\`eve
 la d\'emonstration.\cqfd

\begin{fscor}\label{ImageFermee} \textup{\cite[1.3]{MB2}} Sous les hypoth\`eses de \textup{\ref{propre:notations}},
 soit $f:X\to Y$ un morphisme \emph{propre} de $K$-espaces. Alors l'image de $f_{\Top}$ est ferm\'ee dans $Y_{\Top}$.
\end{fscor}
\dem appliquer le th\'eor\`eme \ref{th:propre} au ferm\'e $\Phi=X_{\Top}$ et \`a un point  $y\notin \im(f_{\Top})$.\cqfd

\begin{fscor}\label{cor:fibrescompactes} Sous les hypoth\`eses de \textup{\ref{propre:notations}}, soit $f:X\to Y$ 
un morphisme \emph{propre} de $K$-espaces.
\begin{numlist}
\item\label{cor:fibrescompactes1} Posons $Z:=\left\{z\in Y_{\Top}\mid f_{\Top}^{-1}(z)\textup{ est compact}\right\}$.
Alors la restriction $f_{\Top,Z}:f_{\Top}^{-1}(Z)\to Z$ de $f_{\Top}$ au-dessus de $Z$ est topologiquement propre.\smallskip
\item\label{cor:fibrescompactes2} Posons $Z_{1}:=\left\{z\in Y_{\Top}\mid \mathrm{Card\,}\left(f_{\Top}^{-1}(z)\right)=1\right\}$.
Alors la restriction $f_{\Top,Z_{1}}$ de $f_{\Top}$ au-dessus de $Z_{1}$ est un hom\'eomorphisme.
\end{numlist}
\end{fscor}
\dem pour \ref{cor:fibrescompactes1}, il suffit de montrer que $f_{\Top,Z}$ est ferm\'ee (puisqu'elle est \`a fibres compactes, 
cf. \cite[I, \S\,10, \no 2, th. 1]{BTG}). Or c'est une cons\'equence imm\'ediate de la derni\`ere assertion de \ref{th:propre}. 
On en d\'eduit \ref{cor:fibrescompactes2} puisque $f_{\Top,Z_{1}}$ est bijective par construction, et  elle est propre d'apr\`es
 \ref{cor:fibrescompactes1} car $Z_{1}\subset Z$. \cqfd

\begin{fscor}\label{fini} Sous les hypoth\`eses de \textup{\ref{propre:notations}}, soit $f:X\to Y$ un morphisme \emph{fini} de $K$-espaces.  
Alors $f_{\Top}$ est ferm\'ee (et donc propre).\cqfd
\end{fscor}

\begin{fsrem}
Le corollaire \ref{fini} peut s'obtenir de mani\`ere plus directe et \'el\'ementaire: pour le cas d'une valuation de rang 1, voir
 par exemple \cite[proposition 2.2.1]{MB1}.
\end{fsrem}

\Subsection{Groupes $G$ tels que $G^\circ_{\red}$ soit lisse}\label{tors-lisse2}
On d\'esigne par $F$ un corps  topologiquement hens\'elien et 
 v\'erifiant la propri\'et\'e suivante: pour tout morphisme \emph{fini} $f:X\to Y$ de $F$-vari\'et\'es, l'application 
$f_{\Top}:X_{\Top}\to Y_{\Top}$ est ferm\'ee.
 
 D'apr\`es \ref{fini}, ces conditions sont v\'erifi\'ees notamment si $F$ est un corps valu\'e admissible.
 
\begin{fsteo}\label{SsGpesLissesQuotFini} 
Soit $G$ un $F$-groupe alg\'ebrique tel que le $F$-sch\'ema $G^{\circ}_{\red}$ soit lisse 
\emph{(c'est donc automatiquement un sous-$F$-groupe de $G^\circ$)}. 
Soient $Y$ un $F$-espace et $f:X\to Y$ un $G$-torseur. 

Alors $f_{\Top}:X_{\Top}\to Y_{\Top}$ 
est ouverte sur son image; en outre celle-ci 
est ferm\'ee dans $Y_{\Top}$, et elle est  ouverte si $F$ est parfait.

Si de plus $Y$ est localement s\'epar\'e, alors $f_{\Top}$ fait de $X_{\Top}$ 
un $G_{\Top}$-fibr\'e principal sur son image.

\end{fsteo}

\begin{fssrem}\label{rem:SsGpesLissesQuotFini} (cf. remarque \ref{rem:prop:tors-strict}) Dans le cas d'un corps valu\'e admissible, 
cet \'enonc\'e est, tout comme \ref{prop:tors-strict}, cons\'equence du th\'eor\`eme principal \ref{ThPpal}, qui sera \'etabli sans
 en faire usage.
\end{fssrem}

\medskip

\noindent\textsl{D\'emonstration de \rref{SsGpesLissesQuotFini}:} la question \'etant locale sur $Y_{\Top}$, 
 nous pouvons supposer que $Y$ est un sch\'ema. 

On a  $G^{\circ}_{\red}\subset G^\bec$, de sorte que $G/G^\bec$ est fini sur $F$.

Posons $Z:=X/G^\bec$, de sorte que $f$ se d\'ecompose en $X\xrightarrow{h}Z\xrightarrow{g}Y$.
Comme $G^\bec$ est lisse et $F$ topologiquement hens\'elien, la proposition \ref{prop:tors-groupe-lisse} 
montre que $h_{\Top}$ est ouverte 
et que son image $Z'_{\Top}$ est (ouverte et) ferm\'ee dans $Z_{\Top}$.

D'autre part $g$ est un morphisme fini (et en particulier $Z$ est un sch\'ema), donc $g_{\Top}$ est ferm\'ee; en outre, 
le corollaire \ref{quotient-torseur} nous dit qu'elle est injective. C'est donc un hom\'eomorphisme sur un ferm\'e de $Y_{\Top}$, ainsi
 que sa restriction \`a $Z'_{\Top}$. Enfin, si $F$ est parfait, alors $G^\bec=G_{\red}$ et $g$ est radiciel, donc ($F$ \'etant parfait)
 $g_{\Top}$ est bijective et est donc un hom\'eomorphisme: l'image de $f_{\Top}$ s'identifie via $g_{\Top}$ \`a celle, ouverte, de $h_{\Top}$. 

Enfin, si $Y$ est localement s\'epar\'e, alors $Z$ l'est aussi puisque $g$ 
est un mor\-phisme s\'epar\'e, donc $X_{\Top}$ est un $G_{\Top}$-fibr\'e principal au-dessus de $Z'_{\Top}$ (\ref{prop:tors-groupe-lisse}
 \`a nouveau).\cqfd

\begin{fscor}\label{cor:strictCasMult}
Soient  $G$ un $F$-groupe alg\'ebrique et $f:X\to Y$ un 
$G$-torseur, o\`u $X$ et $Y$ sont des $F$-espaces. On suppose que $G^\circ$ est  \emph{de type multiplicatif}.
 Alors $f_{\Top}:X_{\Top}\to Y_{\Top}$ est ouverte sur son image, qui est ferm\'ee dans $Y_{\Top}$; si $Y$ est localement s\'epar\'e, 
c'est m\^eme un $G_{\Top}$-fibr\'e principal sur cette image.
\end{fscor}
\dem cela r\'esulte de \ref{SsGpesLissesQuotFini} et du fait que $G^\circ_{\red}$ est lisse.\cqfd

\section{Un th\'eor\`eme de compactification }\label{section-compactification}
\modif{%
On d\'esigne par $k$ un corps quelconque, et par $k_{s}$ une cl\^oture 
s\'eparable de $k$;  les r\'esultats de cette section n'ont d'int\'er\^et que si $k$ n'est pas parfait.

\begin{fdefi}\label{defi:compactification} Soit $X$ une $k$-vari\'et\'e munie d'une action d'un $k$-sch\'ema en groupes $G$. Nous appellerons \emph{compactification $G$-\'equivariante} de $X$ une immersion ouverte $G$-\'equivariante
 $j:X\inj X^c$ o\`u $X^c$ est une $k$-vari\'et\'e \emph{propre} munie d'une action de $G$. 
 
Dans cette situation, les points et sous-sch\'emas de $X^c\setminus X$ seront dits \og \`a l'infini\fg.
\end{fdefi}

On observera que, dans cette d\'efinition,  l'immersion $j$ n'est pas suppos\'ee sch\'ematiquement dense. Nous avons fait ce choix pour la commodit\'e de la r\'edaction; pour les applications o\`u la densit\'e est requise, il suffit de faire appel au lemme suivant:

\begin{fslem}\label{lem:densite} Soit $Z$ une $k$-vari\'et\'e munie d'une action d'un $k$-sch\'ema en groupes $G$, et soit $Y\subset Z$ un sous-sch\'ema stable par $G$. Alors l'adh\'erence sch\'ematique $\ol{Y}$ de $Y$ dans $Z$ est stable par $G$.
\end{fslem}
\dem consid\'erons le diagramme commutatif
$$\xymatrix{G\times_{k}Y \ar[d]^{\beta} \ar@{}[r]|{\subset} & G\times_{k}\ol{Y}\ar[d]^{\alpha}\\
Y\ar@{}[r]|{\subset} & X
}$$
o\`u les fl\`eches $\alpha$ et $\beta$ sont d\'eduites de l'action de $G$. Il s'agit de voir que $\alpha$ se factorise par $\ol{Y}$. Or $\alpha^{-1}(\ol{Y})$ est un sous-sch\'ema ferm\'e de $G\times_{k}\ol{Y}$ qui contient $G\times_{k}Y$, et $G\times_{k}Y$ est sch\'ematiquement dense dans $G\times_{k}\ol{Y}$ puisque l'immersion $Y\inj \ol{Y}$ est sch\'ematiquement dense et que $G$ est plat sur $\Spec(k)$. Donc  $\alpha^{-1}(\ol{Y})= G\times_{k}\ol{Y}$, d'o\`u la conclusion.\cqfd
\medskip
}

Le premier auteur a annonc\'e des  r\'esultats de compactification g\'en\'eraux pour les groupes
alg\'ebriques \cite{Ga} qui g\'en\'eralisent un r\'esultat de Borel-Tits pour les  groupes alg\'ebriques affines 
sur un corps parfait \cite[th. 8.2]{BT1} ainsi que  le cas d'un $k$-groupe commutatif \cite[\S\,10.2, th. 7]{BLR}.
Nous nous  int\'eressons ici \`a la variante suivante.
\modif{%
\begin{fteo}\label{ofer} Soient $G$ un $k$-groupe alg\'ebrique et $H=G^\bec$ le plus grand sous-groupe lisse 
de $G$ \textup{(cf.\ \ref{ssec:pgsgl})}. Soit $J$ un $k$-groupe alg\'ebrique lisse agissant sur $G$ par automorphismes
de groupe.
Alors $X:=G/H$ admet  une compactification $G \rtimes J$-\'equivariante   $X^c$ munie d'un fibr\'e 
en droites ample $G \rtimes J$-lin\'earis\'e et
satisfaisant $X(E)= X^c(E)$ pour toute extension s\'eparable $E$ de 
$k$ \emph{(ce qui \'equivaut \`a dire que $X^c(k_{s})=\{x_{0,k_{s}}\}$,
 o\`u $x_{0}$ d\'esigne la classe neutre de $X(k)$)}. 

Si de plus $G$ satisfait la condition $(*)$ de \rref{condition:star}, 
on peut imposer en outre la condition que pour toute extension s\'eparable $K$ de $k$,
$X_{K}^c$ n'admette aucune $K$-orbite \`a l'infini pour l'action de $G_{K}$.
\end{fteo}

Suivant le \S \ref{ssec:pgsgl}, l'action du $k$-groupe lisse $J$ sur $G$ normalise  $H$ si bien que l'on a une action
de $G \rtimes_k J$ sur $G/H$.
Si $G$ est affine et $J=1$, alors l'\'enonc\'e ci-dessus est un cas particulier de \cite[Th. A]{Ga}. Il se 
trouve que ce cas particulier est bien plus simple que le cas g\'en\'eral,
et nous en donnons une d\'emonstration ci-dessous.

\Subsection{D\'emonstration du th\'eor\`eme \ref{ofer}: d\'evissage}

\begin{fsdefi}\label{defi:compactif} Soient $k$ un corps, $k_{s}$ une cl\^oture s\'eparable de $k$, $G$ un $k$-groupe alg\'ebrique, 
$H$ un sous-groupe alg\'ebrique de $G$. Soit $J$ un $k$-groupe alg\'ebrique qui agit sur $G$ par automorphismes de groupes
et qui normalise $H$. On note  $X$ le $G$-espace homog\`ene $G/H$. 
\smallskip

\noindent Une \emph{bonne $G\rtimes_k J$-compactification} de $X$ est la donn\'ee:
\begin{itemize}
\item d'une $k$-vari\'et\'e projective $X^c$, munie d'une action \`a gauche de $G \rtimes_k J$,
\item d'un fibr\'e en droites ample $G \rtimes_k J$-lin\'earis\'e $L$ sur $X^c$,
\item d'une immersion ouverte $G \rtimes_k J$-\'equivariante $j:X\inj X^c$,
\end{itemize}
tels que le groupe $G(k_{s})$ op\`ere \emph{transitivement} sur $X^c(k_{s})$ 
\emph{(ce qui implique notamment que $X^c(k_{s})=j(X(k_{s}))$).}

Si de plus  $X^c$ n'admet aucune  $G$-orbite \`a l'infini d\'efinie sur une extension s\'eparable de $k$, 
on  dit que $X^c$ est une \emph{tr\`es bonne compactification}. 
\end{fsdefi}
}

\begin{fsrems}\label{rem:strategie} 

\noindent\textup{(1)} Dans l'\'enonc\'e du th\'eor\`eme \ref{ofer}, la condition sur $X^c(k_{s})$ \'equivaut trivialement 
 \`a dire que l'action de $G(k_{s})$ sur $X^c(k_{s})$ est transitive. 
On peut donc reformuler  \ref{ofer} en disant que \emph{$X$ admet une bonne compactification, et une tr\`es bonne
 compactification si $G$ v\'erifie $(*)$.}\smallskip

\noindent\textup{(2)}  Pour d\'emontrer \ref{ofer}, nous allons donc proc\'eder comme suit. Si le r\'esultat est en d\'efaut
 pour $(G,H)$ 
 \modif{%
 et $J$, il existe un plus petit sous-groupe $G'$ (normalis\'e par $J$) de $G$ contenant $H$ tel qu'il 
 soit encore en d\'efaut pour $(G',H)$.
 Rempla\c{c}ant $G$ par $G'$, nous pouvons donc supposer que pour tout sous-groupe $H'$ de $G$ normalis\'e par $J$
 tel que $H\subset H'\varsubsetneq G$,
 le $H'$-espace homog\`ene $H'/H$ admet une bonne $H' \rtimes J$-compactification. La construction fondamentale 
 de la d\'emonstration consiste  alors \`a exhiber un tel sous-groupe $H'$ tel que $G/H'$ ait une bonne 
 $G \rtimes J$-compactification; 
 }
 on conclut ensuite gr\^ace \`a 
l'\'enonc\'e de d\'evissage \ref{prop:devissage} ci-dessous, inspir\'e des constructions de \cite[\S\,10.2]{BLR}.
\end{fsrems}

\modif{%
\begin{fsprop}\label{prop:devissage} Soient $k$ un corps, $G$ un $k$-groupe alg\'ebrique, $H\subset H'$ des sous-groupes
 alg\'ebriques de $G$. Soit $J$ un \modif{$k$-groupe alg\'ebrique  agissant} sur $G$ par automorphismes de groupe
  \modif{en laissant stables} $H$ et $H'$.

On suppose que $G/H'$  admet une bonne (resp. tr\`es bonne) $G \rtimes J$-compac\-ti\-fi\-ca\-tion et 
que $H'/H$  admet une bonne (resp. tr\`es bonne) $H' \rtimes_k J$-compac\-ti\-fi\-ca\-tion.
Alors $G/H$ admet une bonne (resp. tr\`es bonne) $G \rtimes_k J$-compac\-ti\-fi\-ca\-tion. 
\end{fsprop}

\dem posons $X=G/H'$, $Y=G/H$, et consid\'erons les morphismes canoniques (et $G \rtimes_k J$-\'equivariants) $G\to Y\xrightarrow{q} X$. 
Le compos\'e  $G\to X$ fait de $G$ un $H'$-torseur \`a droite au-dessus de $X$. 

Posons $Z=H'/H$ (qui s'identifie \`a $q^{-1}(x_{0})$ o\`u $x_{0}\in X(k)$ est la classe neutre). Il existe par hypoth\`ese 
une bonne $H'\rtimes_k J$-compactification $Z\inj Z^c$ de $Z$, munie notamment d'un fibr\'e en droites ample
$H'\rtimes_k J$-lin\'earis\'e $N$. 

On consid\`ere alors le produit contract\'e $Y^c:=G\PC{H'}Z^c\to X$. D'apr\`es 
le \S\,\ref{subsec-prodcontract}, 
$Y^c$ est repr\'esentable par un $k$-sch\'ema projectif sur $X$; de plus,  $M:=G\PC{H'}N$ est un fibr\'e en droites 
$G \rtimes_k J$-lin\'earis\'e sur $Y^c$, ample relativement \`a $X$.

Par hypoth\`ese, $X$ admet une bonne $G \rtimes_k J$-compactification $X^c$; nous avons donc un diagramme commutatif de $k$-vari\'et\'es 
munies d'actions \`a gauche de $G \rtimes_k J$:
\begin{equation}\label{diag:devissage}
\xymatrix{Y=G/H
\ar@{^{(}->}^(.45){i}[r] \ar[dr]_(.45){q} 
& Y^c=G\PC{\smash{H'}}Z^c \ar[d]^{q^c}\\
& X=G/{H'} \enskip \ar@{^{(}->}^(.6){j}[r] & X^c
}
\end{equation}
o\`u tous les morphismes sont $G \rtimes_k J$-\'equivariants, $i$ et $j$ \'etant \modif{des immersions ouvertes}; 
de plus $Y^c$ (resp. $X^c$) 
est projectif sur $X$ (resp. sur $k$) et est muni d'un faisceau inversible $G\rtimes J$-lin\'earis\'e $M$ (resp. $L$),
ample relativement 
\`a $X$ (resp. $k$). On a alors le r\'esultat de prolongement suivant:
}

\begin{fsslem}\label{lem:devissage1} 
Consid\'erons un diagramme $P\xrightarrow{\pi}U\overset{j}{\hookrightarrow}  V$ 
de $k$-vari\'et\'es munies d'actions \modif{d'un $k$-groupe alg\'ebrique $\Gamma$}, les morphismes \'etant \modif{$\Gamma$}-\'equi\-va\-riants. On suppose que:
\begin{itemize}
\item $j$ est une immersion ouverte;
\item $\pi$ est projectif et $P$ admet un faisceau inversible \modif{$\Gamma$}-lin\'earis\'e $\calL$, ample relativement \`a $\pi$.
\end{itemize}
Alors il existe un diagramme cart\'esien \modif{$\Gamma$}-\'equivariant
$$
\xymatrix{P\,\ar@{^{(}->}[r]\ar[d]_{\pi}  &  P'\ar[d]^{\pi'}\\
U\,\ar@{^{(}->}[r]^{j} &  V
}$$
o\`u $\pi'$ est projectif et o\`u $P'$ admet un faisceau inversible \modif{$\Gamma$}-lin\'earis\'e, ample relativement \`a $\pi'$ et
 prolongeant une puissance de $\calL$.
\end{fsslem}
\dem quitte \`a remplacer $\calL$ par une puissance convenable, on peut le supposer tr\`es ample relativement \`a $\pi$, 
de sorte que l'on a une $U$-immersion ferm\'ee canonique $P\inj \PP(\cE)$ o\`u $\cE=\pi_{*}\calL$; de plus, comme $\calL$ 
est \modif{$\Gamma$}-lin\'earis\'e, l'action de \modif{$\Gamma$} sur $P$ est induite par une \modif{$\Gamma$}-lin\'earisation de $\cE$ sur $U$. D'apr\`es un  
lemme de prolongement de Thomason \cite[2.2]{T}, il existe un $\OOO_{V}$-module coh\'erent \modif{$\Gamma$}-lin\'earis\'e 
$\cE'$ prolongeant $\cE$. Il suffit de prendre pour $P'$ l'adh\'erence sch\'ematique de $P$ dans $\PP(\cE')$\modif{, qui est stable sous $\Gamma$ d'apr\`es \ref{lem:densite}}.\cqfd\medskip

\modif{%
Si l'on applique le lemme  \ref{lem:devissage1} 
\modif{au diagramme $Y^c\to X\to X^c$, avec les actions de  $\Gamma=G \rtimes_k J$}, 
on voit que le diagramme (\ref{diag:devissage}) se compl\`ete en 
$$\xymatrix{Y\ar@{^{(}->}^{i}[r]\ar[dr]_{q} & Y^c\ar@{}[dr]|{\square}\ar[d]^{q^c} \ar@{^{(}->}^{i^c}[r] & Y^{cc}\ar[d]^{q^{cc}} \\
& X\ar@{^{(}->}^{j}[r] & X^c
}$$
dans lequel le carr\'e est cart\'esien, $q^{cc}$ est projectif et $Y^{cc}$ est muni d'un faisceau $q^{cc}$-ample 
$G\rtimes_k J$-lin\'earis\'e 
$M^c$, prolongeant une puissance de $M$. Dans ces conditions, pour $m\in\NN$ assez grand, le faisceau inversible 
(\'evidemment ${G \rtimes_k J}$-\-lin\'earis\'e) $M^c\otimes (q^{cc})^*L^{\otimes m}$ est ample sur $Y^{cc}$ \cite[(4.6.13)(ii)]{EGA2}. 
}

Il reste \`a montrer que $G(k_{s})$ op\`ere transitivement sur $Y^{cc}(k_{s})$. Soit donc $y\in Y^{cc}(k_{s})$: montrons
 que $y\in G(k_{s}).y_{0,k_{s}}$ o\`u $y_{0}$ est la classe neutre de $Y$. L'image de $y$ dans $X^c$ appartient \`a $X^c(k_{s})$
 donc \`a l'orbite sous $G(k_{s})$ de l'origine $x_{0,k_{s}}$ de $X(k_{s})$, puisque $X^c$ est une bonne compactification. 
On peut donc, en faisant op\'erer $G(k_{s})$, supposer que $q^{cc}(y)=x_{0,k_{s}}$. Or $(q^c)^{-1}(x_{0})$ s'identifie \`a $Z^c$,
 et l'hypoth\`ese de bonne compactification pour $Z$ entra\^{\i}ne que $H'(k_{s})$ op\`ere transitivement sur $Z^c(k_{s})$, de
 sorte que $y$ est bien dans l'orbite de $y_{0,k_{s}}$. 

On suppose maintenant que $X^c$ (resp. $Z^c$) est une tr\`es bonne compactification 
de $G/H'$ (resp. $H'/H$). Soit $E/k$ une extension s\'eparable 
et  $W \subset Y^{cc} \times_k E$ un  sous-espace localement ferm\'e qui soit une $G$-orbite.
On veut montrer que $W \subset Y \times_k E$. Sans perte de g\'en\'eralit\'e, il est loisible de 
supposer que $E=k$.  La  projection de $W$ par  $q^{cc}$ est une  $G$-orbite de $X^c$, donc
est un sous-espace  de $X$ selon notre hypoth\`ese. Par  suite
$W$ est un sous-espace de $Y^c$. Par homog\'en\'eit\'e, l'application  $W \to X$ est surjective.
On note $x_0$ le point privil\'egi\'e de $X=G/H$. Alors
la fibre $W_{x_0}$ est un sous-espace localement ferm\'e de $Z^c$, qui est une 
$H'$-orbite.  Notre hypoth\`ese implique que $W_{x_0} \subset Z \subset Y$.
Par homog\'en\'eit\'e, on conclut que  $W \subset Y$. 
\cqfd

\Subsection{D\'emonstration du th\'eor\`eme \ref{ofer}: construction et fin}

Posons $G_{\mathrm{af}}=\Spec{\Gamma(G,\cO_{G})}$: alors \cite[VI$_{\mathrm{B}}$, th\'eor\`eme 12.2]{SGA3},
 le morphisme naturel $G\to G_{\mathrm{af}}$ induit un isomorphisme $G/N\cong G_{\mathrm{af}}$, o\`u $N\subset G$
 est distingu\'e et lisse, et en particulier contenu dans $H$. On peut donc remplacer $G$ par $G_{\mathrm{af}}$ et
  supposer que $G$ est affine 
  \modif{%
  (en gardant bien s\^ur l'action de $J$).

On note $k'$ le corps de d\'e\-fi\-ni\-tion du $\overline k$-groupe lisse
${(G \times_k \overline k)_{\red}}$; c'est une extension finie radicielle de $k$
dont on note $d$ le degr\'e. On dispose d'une immersion ferm\'ee $G \rtimes_k J$-\'equivariante
$j: G \to  \W_{k'/k} G_{k'}$ \cite[A.5.7]{CGP}; 
on identifiera $G$ et $j(G)$. 
On pose $\widetilde H = {G\cap \W_{k'/k} G_{k', \red} }\,  \subset \,  \W_{k'/k} G_{k'}$:
ce  $k$-groupe contient $H$ et satisfait  $H(k_s)=\widetilde H(k_s)$.
On observe que \modif{l'action de $J$ respecte} $\widetilde H$. En effet, le morphisme d'action ${J_{k'} \times_{k'} G_{k'}} \to G_{k'}$ donne lieu 
\`a un $k$-morphisme ${\bigl(J_{k'} \times_{k'} G_{k', \red} \bigr)_{\red}} \simlgr  \bigl(J_{k'} \times_{k'} G_{k'} \bigr)_{\red} \to G_{k', \red}$.
Comme $J_{k'}$ est \modif{g\'eo\-m\'e\-tri\-que\-ment} r\'eduit,  $J_{k'} \times_{k'} G_{k',\red}$ est r\'eduit
\cite[5.49.ii]{GW}, d'o\`u un  morphisme $J_{k'} \times_{k'} G_{k', \red} \to G_{k', \red}$. Ainsi
$J_{k'}$ normalise  $G_{k', \red}$ et $J$ normalise $\W_{k'/k} G_{k', \red}$ et $\widetilde H$.
}

\begin{fslem} \label{arret} $\widetilde H = G$ si et seulement si $G$ est lisse.  
\end{fslem}

\noindent{\it D\'emonstration.} Si $G$ est lisse, alors $k=k'$ et 
$\widetilde H=G =  \W_{k'/k} G_{k'} = \W_{k'/k} G_{k', \red}$. 
R\'eciproquement, on suppose que  $\widetilde H = G$, 
c'est-\`a-dire que $j(G)\subset\W_{k'/k} G_{k', \red}$. 
On a donc un diagramme commutatif
$$\xymatrix{& \W_{k'/k} G_{k', \red} \ar@{_{(}->}[d]  \\
G\: \ar@<.4ex>@{-->}[ur]\ar@<-.4ex>@{^{(}->}[r] & \W_{k'/k} G_{k'}.
}$$
Effectuant le changement de base $k'/k$, on obtient un diagramme commutatif
$$\xymatrix{& \left(\W_{k'/k} G_{k', \red}\right)_{k'} \ar@{_{(}->}[d] \ar[r]^(.6){q_1} & G_{k',\red}\ar@{_{(}->}[d] \\
G_{k'}\: \ar@<.2ex>@{-->}[ur]\ar@{^{(}->}[r] & \left(\W_{k'/k} G_{k'}\right)_{k'} \ar[r]^(.6){q} & G_{k'}
}$$
o\`u $q$ et $q_{1}$ sont les morphismes d'adjonction \cite[A.5.7]{CGP}. Or, la compos\'ee des fl\`eches inf\'erieures
 est l'identit\'e de $G_{k'}$: on conclut que $\id_{G_{k'}}$ se factorise par l'inclusion de $G_{k',\red}$ dans $G_{k'}$, 
ce qui signifie que $G_{k',\red}=G_{k'}$, donc que $G$ est lisse.
\cqfd
\medskip

Le $k'$-sch\'ema $Q':=G_{k'}/ G_{k', \red}$ est le spectre d'une $k'$-alg\`ebre locale artinienne
$A'$ de corps r\'esiduel $k'$; on note  $s:\Spec{(k')}\inj Q'$ son point ferm\'e. 
Nous allons appliquer \`a cette situation les constructions du \S\,\ref{orbite:bord}, dont nous reprenons 
les notations: ainsi, $V$ d\'esigne le $k$-sch\'ema sous-jacent \`a $Q'$, $W= \W_{k'/k}Q'$ la restriction de Weil 
de $Q'$ sur $k$, et $V^{[d]}=\Hilb^d_{V/k}$ le sch\'ema de Hilbert des sous-sch\'emas finis de longueur $d$ de $V$.
\modif{%
 Noter qu'ici $V$ est fini, de sorte que  $V^{[d]}$ est projectif sur $k$ et admet un fibr\'e en droites ample $G\rtimes_k J$-lin\'earis\'e
 (lemme \ref{lem:linearisation}), et que $W$ est affine comme restriction de Weil d'un $k'$-sch\'ema affine.

On a de plus une action naturelle de $G \rtimes_k J$ sur $V$ compatible au
morphisme $\sigma:V\to\Spec(k')$ d\'efinissant le $k'$-sch\'ema  $Q'$. 

Selon le \S\,\ref{orbite:bord}, on dispose d'un $k$-morphisme qui est  une immersion ouverte $G \rtimes_k J$-\'equi\-va\-riante
$$ 
u: W  \, \to \, V^{[d]}.
$$
En particulier, le point $s\in W(k)=Q'(k')$ a pour image dans $V^{[d]}(k)$ le 
sous-sch\'ema ferm\'e r\'eduit $\{s\}$ de $V$, qui est bien de degr\'e $d$ sur $k$.  
Par abus, on notera encore $s$ le point correspondant de $W(k_{s})$.

Observons que $V^{[d]}(k_{s})=W(k_{s})=\{s\}$: ceci r\'esulte du fait que $V\otimes_{k}k_{s}$ est local de corps 
r\'esiduel $k'\otimes_{k}k_{s}$, extension de degr\'e $d$ de $k_{s}$. 

Le groupe $G$ agit sur le couple $(V^{[d]},W)$, et le stabilisateur
de l'unique $k$-point $s$ de  $V^{[d]}(k)$ est $\widetilde H$.
Ceci produit  un morphisme d'orbite
$i: G/ \widetilde H \to W$ qui est une  immersion (\ref{cor:orbites}). 
Utilisons  la m\^eme notation $X$ pour le quotient $G/ \widetilde H$ et pour son image par 
le compos\'e  $i: G/\widetilde H \to W \to V^{[d]}$, et notons $X^c$ l'adh\'erence sch\'ematique de $X$ dans $V^{[d]}$.
 Alors $X^c$ est une compactification $G\rtimes_j J$-\'equivariante de $X$, et 
 c'est m\^eme une bonne $G\rtimes_k J$-compactification puisque $X^c(k_{s})=X(k_{s})=\{s\}$. 
 Enfin le lemme \ref{arret} nous dit que $\til{H}\varsubsetneq G$, 
sauf dans le cas trivial o\`u $G$ est lisse. Comme expliqu\'e dans la remarque \ref{rem:strategie}\,(2), ceci 
suffit \`a conclure compte tenu de la proposition \ref{prop:devissage}.
}
\smallskip

Nous allons raffiner l'argument dans le cas o\`u le groupe $G$ satisfait la condition $(*)$
afin de montrer que la compactification construite est tr\`es bonne.
Tout d'abord, le $k$-groupe $\til{H}$ satisfait lui aussi l'hypoth\`ese $(*)$ puisqu'il
con\-tient $G^\bec$ (lemme \ref{lem:star}\,\ref{lem:star2}) et la r\'ecurrence fonctionne 
bien en appliquant le raffinement de la proposition \ref{prop:devissage}. 
L'unique chose \`a v\'erifier est que $X^c$ n'a aucune $G$-orbite  \`a l'infini d\'efinie sur une extension s\'eparable de $k$.
\modif{Soient donc  $E$ une extension s\'eparable de $k$ et $I\subset (X^c)_{E}$ une $G_{E}$-orbite. Alors 
$k'\otimes_{k}E$ est un corps, donc le th\'eor\`eme \ref{th:pas-d'orbite}  indique que $I \subset W_{E}$. 
}
Pour conclure, il suffit donc de montrer que  $X=X^c\cap W$, 
c'est-\`a-dire que $X$ est \emph{ferm\'e} dans $W$. Comme $W$ est affine, cette assertion r\'esulte 
(compte tenu de l'hypoth\`ese (*)) du lemme \ref{Rosenlicht}, ce qui ach\`eve la d\'e\-mons\-tration.\cqfd

\section{D\'emonstration du th\'eor\`eme \ref{ThPpal}}\label{caslocal}

On suppose d\'esormais que $K$ est un corps valu\'e admissible, et l'on se donne $G$ et $f:X\to Y$ comme dans \ref{ThPpal}, dont nous 
allons achever la preuve. Les cas particuliers o\`u  $G$ est lisse (\ref{prop:tors-groupe-lisse})  et o\`u $K$ est parfait 
(\ref{SsGpesLissesQuotFini}) ont d\'ej\`a \'et\'e trait\'es. Il reste \`a montrer les assertions \ref{ThPpal1}\ref{ThPpal1a}  et  \ref{ThPpal1}\ref{ThPpal1c} ($I\subset Y_{\Top}$ est localement ferm\'e, et est ferm\'e sous la condition $(*)$), 
\ref{ThPpal2} (${f_{\Top}}$ est ouverte sur  $I$), et \ref{ThPpal3} ($X_{\Top}\to I$ est un $G_{\Top}$-fibr\'e principal si 
$Y$ est localement s\'epar\'e).
On fixe une cl\^oture s\'eparable de $K$, not\'ee  $K_{s}$.

Notons $H=G^\bec$ le plus grand sous-groupe lisse de $G$ (\ref{ssec:pgsgl}). 
Suivant le th\'eor\`eme \ref{ofer}, $G/H$ admet une compactification $G$-\'equi\-va\-riante $(G/H)^c$, ayant l'origine comme unique point 
s\'eparable. Sous la condition (*), nous supposerons en outre que $(G/H)^c$ est une tr\`es bonne compactification.

Consid\'erons le produit contract\'e $Z^c:=X \PC{G} (G/H)^c$: c'est une compactification relative \`a $Y$ (c'est-\`a-dire propre sur $Y$) 
de $Z:=X \PC{G} (G/H)$; ce dernier n'est autre que le quotient $X/H$. Posons $Z^\infty:=Z^c \setminus Z$ (avec sa structure r\'eduite, 
par exemple). On a un diagramme commutatif

$$\xymatrix{%
X\ar[drr]_(.3){f}\ar[r]^{\pi} & Z\: \ar[dr]^(.55){h}\ar@{^{(}->}[r]^{j} & Z^c \ar[d]^(.4){h^c} \ar@{<-^{)}}[r]^{i} & 
\:Z^\infty \ar[dl]^(.25){h^\infty}\\
&&Y
}$$
dans lequel $\pi$ est un torseur sous le groupe \emph{lisse} $H$, $j$ est une immersion ouverte, 
$i$ est l'immersion ferm\'ee compl\'ementaire, $h^c$ et $h^\infty$ sont propres.

\begin{flem}\label{clef}
 $h_{\Top}: Z(K)\to Y(K)$ est injective, et les images de $h_{\Top}$ et de $h^\infty_{\Top}$ sont disjointes. 
 Si de plus $(G/H)^c$ est une tr\`es bonne compactification \emph{(et notamment si $G$ v\'erifie (*), d'apr\`es nos conventions)}, 
alors on a $Z^c(K)=Z(K)$, de sorte que $\im(h_{\Top})=\im(h^c_{\Top})$.
 \end{flem}

\dem l'injectivit\'e r\'esulte du corollaire \ref{quotient-torseur}. Notons $L$ l'image de $h_{\Top}$. Pour  $y\in L$, consid\'erons les
 fibres $Z_{y}\subset Z^c_{y}$ de $h$ et $h^c$ en $y$ de sorte que $Z_{y}(K)$ a  un seul \'el\'ement $z$, qui est m\^eme le seul
 point s\'eparable de $Z_{y}$. Puisque $\pi$ est lisse et surjectif, $\pi^{-1}(z)(K_s)\neq\emptyset$, et donc le $G$-torseur $X_{y}$ de\-vient
 trivial sur $K_{s}$, de sorte que $(Z^c_{y})^{\mathstrut}_{K_{s}}$ est isomorphe \`a $(G/H)^c_{K_{s}}$ et a donc, lui aussi, un seul
 point s\'eparable; autrement dit,  $Z^\infty_{y}(K_{s})=\emptyset$. Ceci montre que $L\cap\im(h^\infty_{\Top})=\emptyset$. 
 
Supposons que $(G/H)^c$ soit une tr\`es bonne compactification. Soit 
$z \in Z^c(K)$ et montrons que $z\in Z(K)$. Posons $y=h^c(z)\in Y(K)$.  
Alors $Z^c_{y}$ s'i\-den\-ti\-fie
 \`a $X_{y}\PC{G}(G/H)^c$ qui est muni d'une action \`a gauche du $K$-groupe al\-g\'e\-brique  $G'=\ul\Aut_G(X_y)$ (\ref{torsion}).
 Consid\`erons la $K$-orbite $T$ de $z$ dans $Z^c_y$ sous $G'$  (remarque \ref{rem:orbites}.\ref{rem:orbites2}). 
 Celle-ci correspond canoniquement, d'a\-pr\`es \ref{torsion}, \`a une $K$-orbite $T_{0}$ de $(G/H)^c$ sous l'action de $G$. Vu notre hy\-po\-th\`ese, on a $T_{0}\subset G/H$, et donc $T\subset Z_{y}$ et  $z\in Z(K)$.
\cqfd

\medskip

Ceci implique d\'ej\`a que
$$L= \im(h^c_{\Top}) \setminus \im(h^\infty_{\Top})$$
et, puisque $h^c$ et $h^\infty$ sont propres, les deux images au membre de droite sont ferm\'ees 
dans $Y_{\Top}$ (corollaire \ref{ImageFermee}), de sorte que $L$ est localement ferm\'e. Sous l'hypoth\`ese (*), on a m\^eme
 $L=\im(h^c_{\Top})$ qui est ferm\'e.

En outre,  l'assertion \ref{cor:fibrescompactes2} du  corollaire \ref{cor:fibrescompactes} entra\^{\i}ne que $h_{\Top}$ est un
 ho\-m\'eo\-mor\-phisme sur son image (c'est en effet la restriction de $h^c_{\Top}$ au-dessus de $L$, qui est contenu dans 
l'ensemble not\'e $Z_{1}$ dans loc. cit.).

En r\'esum\'e, $f_{\Top}$ se d\'ecompose comme suit:
$$ X_{\Top}\xrightarrow{\pi_{\Top}} U\inj Z_{\Top} \fflis L \hookrightarrow Y_{\Top}$$
o\`u l'on a pos\'e $U=\im(\pi_{\Top})\subset Z_{\Top}$ et o\`u:
\begin{itemize}
\item la premi\`ere fl\`eche $X_{\Top}\to U$ est surjective et ouverte, et est 
une $G_{\Top}$-fibration principale si $Y$, et donc $Z$, est localement s\'epar\'e;
\item $U\inj Z_{\Top}$ est un plongement ouvert et ferm\'e;
\item $Z_{\Top} \fflis L$ est un hom\'eomorphisme;
\item $L \hookrightarrow Y_{\Top}$ est  un plongement localement ferm\'e, et ferm\'e sous l'hypoth\`ese (*);
\end{itemize}
(bien entendu, les deux premi\`eres propri\'et\'es r\'esultent  de la proposition \ref{prop:tors-groupe-lisse} 
puisque $\pi$ est un torseur sous le groupe lisse $H$, et que $H_{\Top}=G_{\Top}$). Ceci ach\`eve la d\'emonstration du th\'eor\`eme.
\cqfd 

\section{Exemples et compl\'ements}\label{sec:exemples}

\Subsection{Un exemple d'orbite topologique non ferm\'ee}\label{ex:OrbiteNonFermee}

Soit $K$ un corps topologique (s\'epar\'e et non discret), de caract\'eristique $p>0$ et \emph{non parfait}; 
nous allons donner un exemple d'un $K$-groupe alg\'ebrique $G$ et d'un sous-groupe $H$ tels que l'image de $G(K)$ 
dans $(G/H)(K)$ ne soit pas ferm\'ee. 

On notera $\sigma:K\to K$ l'endomorphisme de Frobenius, et $K_{0}\subset K$ son image (que l'on munira de sa topologie de sous-espace de $K$).
 On consid\`ere l'action du $K$-groupe $G:= \GG_a \rtimes \GG_m$ (produit semi-direct pour l'action standard de $\GG_m$ sur $\GG_a$)
sur la droite affine $\AAA^1_K$ donn\'ee par
$$
(x,y). z= x^p + y^p z \qquad (x \in \GG_a, \;y \in \GG_m\,, \;z \in \AAA^1).
$$
Ainsi $\AAA^1_K$ est un $K$-espace homog\`ene \`a gauche sous $G$ et
le stabilisateur de $0$ est le $K$-sous-groupe ferm\'e $\alpha_p \rtimes_K \GG_m$.   

Fixons  $z\in K\smallsetminus K_{0}$ et consid\'erons l'application d'orbite $\omega_{z}:g\mapsto g.z$ 
 de $G(K)$ dans $K$. Elle est injective et son image est $G(K).z=K_{0}+K_{0}^\times\,z\subset K$, qui ne contient pas $0$; 
mais puisque $K^\times$ est dense dans $K$ ($K$ est non discret), l'adh\'erence de $G(K).z$ contient $0$, donc $G(K).z$ n'est pas ferm\'ee.

Le stabilisateur $G_z$ du point $z$ satisfait $G_z(K_s)=1$ et donc $G_z^\bec=1$.
Par ailleurs, vu que $G$ agit transitivement sur  $\AAA^1_K$,
$G_z$ contient un $\overline k$-tore de rang un. Ainsi le groupe  $G_z$ 
ne satisfait pas la condition $(*)$ 
ce qui est coh\'erent avec le fait que le $G_z$-torseur $G \to G/G_z$ ne satisfait pas le th\'eor\`eme 
\ref{ThPpal} \ref{ThPpal1} \ref{ThPpal1c}.

\begin{fsrems} Le lecteur pourra v\'erifier les compl\'ements suivants:
\begin{numlist}
\item
\begin{subromlist}
\item $G(K).z$ est  localement ferm\'ee dans $K$ si et seulement si $K_{0}$ et $E_{z}:=K_{0}+K_{0}z$
 sont ferm\'es dans $K$; 
\item pour que $\omega_{z}$ soit stricte (donc un hom\'eomorphisme sur son image), il faut et il suffit que $\sigma$ 
soit un hom\'eomorphisme sur $K_{0}$  (ce qui est toujours le cas si $K$ est un corps valu\'e) et que $E_{z}$ soit 
topologiquement libre (comme $K_{0}$-espace vectoriel).
\end{subromlist}
\item L'orbite de $0$ sous $G(K)$ est $K_{0}$; en particulier elle est localement ferm\'ee si et seulement si $K_{0}$ est 
localement ferm\'e (donc ferm\'e) dans $K$. 

L'application d'orbite correspondante $g\mapsto g.0$ est stricte si et seulement si $\sigma$ est un hom\'eomorphisme sur $K_{0}$; 
dans ce cas, $G(K)\to K_{0}$ a une section continue donn\'ee par $t\mapsto (t^{1/p},1)$ donc est un torseur trivial sous $K^\times_{\Top}$.
\end{numlist}
\end{fsrems}

\Subsection{Contre-exemples sur un corps valu\'e hens\'elien non admissible}

Soit $A$ un anneau de valuation discr\`ete hens\'elien, de caract\'eristique $p>0$, et 
soit $v\in\wh{A}$ tel que $v\notin A$ et $v^p\in A$; des exemples d'une telle situation 
ont \'et\'e construits par Nagata \cite[A1, (E2.1)]{N} et par F.K. Schmidt (voir  \cite[11.40]{Ku}, ou \cite[\S\,3.6, Example 11]{BLR}).
 On note $K$ le corps des fractions de $A$, muni de la topologie de la valuation. 

\subsubsection{Une orbite non localement ferm\'ee. }Consid\'erons  le morphisme de 
Frobenius $\Phi: \Aa^1_{K}\to \Aa^1_{K}$: l'image $K_{0}$ de $\Phi_{\Top}$ n'est pas ferm\'ee
 puisque $v^p$ est adh\'erent \`a $K_{0}$, mais n'est pas dans $K_{0}$. Puisque c'est un sous-groupe de $K$,
 elle n'est pas non plus localement ferm\'ee. 

Comme $\Phi$ est fini, ceci montre aussi que l'on ne peut pas supprimer l'hypoth\`ese admissible dans \ref{fini}.

\subsubsection{Une application d'orbite non stricte. }Reprenons maintenant l'exemple de \ref{ex:OrbiteNonFermee} avec
 ce m\^eme corps $K$ et en prenant $z=v^p$. L'application d'orbite (injective) $\omega_{z}$ est 
$$\begin{array}{rcl}
\omega_{z}: K\times K^\times & \ffl & K\\
(x,y) & \longmapsto & x^p+y^p\,z
\end{array}
$$
qui n'est \emph{pas un hom\'eomorphisme sur son image} $K_{0}+K_{0}^\times\,z$: en effet, soit $(t_{n})_{n\in\NN}$ une suite d'\'el\'ements de $K^\times$ 
tels que $\lim\limits_{n\to+\infty}t_{n}^p=z$. Alors la suite $(0,\frac{1}{t_{n}})$ ne converge pas dans $K\times K^\times$, mais son image par $\omega_{z}$
 est la suite $(\frac{z}{t_{n}^p})$ qui converge vers $1=\omega_{z}(0,1)$.

\Subsection{Espaces non localement s\'epar\'es. }\label{ssec:ExNonLocSep} Soit $F$ un corps  topologiquement hens\'elien non discret de 
caract\'eristique dif\-f\'e\-rente de $2$. Consid\'erons le morphisme
$$r:\quad R:=\Aa^1_{F}\amalg\GG_{m,F} \to \Aa^2_{F}$$
dont la restriction \`a $\Aa^1_{F}$ (resp. \`a $\GG_{m,F}$)  est le morphisme diagonal (resp. le morphisme $x\mapsto(x,-x)$). Alors $R$ est une relation 
d'\'equivalence \'etale sur $\Aa^1_{F}$; consid\'erons le faisceau \'etale quotient 
$$\pi: \quad \Aa^1_{F}=:L\to X=L/R$$ 
qui est un $F$-espace. Le morphisme $r$ co\"{\i}ncide avec le monomorphisme naturel $L\times_{X}L\to L\times_{F}L$, de sorte que $r_{\Top}$ n'est pas un plongement 
topologique (son image est la r\'eunion des deux diagonales de $F^2$). On en d\'eduit par \ref{prop:ProdFibEspAlg}\,\ref{prop:ProdFibEspAlg1} que 
la bijection $(L\times_{X}L)_{\Top}\to L_{\Top}\times_{X_{\Top}}L_{\Top}$ n'est \emph{pas un hom\'eomorphisme}.\smallskip

Notons $l_{0}\in L(F)$ l'origine de $L$. La projection $\pi_{\Top}$ \emph{n'est pas un ho\-m\'eo\-mor\-phisme local} en $l_{0}$ puisqu'elle n'est pas localement 
injective (tout voisinage de $l_{0}$ contient un point $l\neq l_{0}$ et son oppos\'e).

\subsubsection{Description de $X_{\Top}$. } Le morphisme $f:x\mapsto x^2$ de $\Aa^1_{F}$ dans $\Aa^1_{F}$ se factorise par $\pi$, donnant naissance \`a un
 diagramme commutatif
$$\xymatrix{ L\ar[r]^{\pi}\ar@/_/@<-5pt>[rr]_{f}& X \ar[r]^{q} & \text{\rlap{$M:=\Aa^1_{F}$}}
}$$
dans lequel $\pi$ est \'etale et surjectif et $q$ induit un isomorphisme au-dessus de $\GG_{m,F}$. Notons $l_{0}\in L(F)$ l'origine de $L$ et $x_{0}$ et $m_{0}$ ses
 images dans $X$ et dans $M$, et $L^*$ (resp. $X^*$, $M^*$) les ouverts de $L$  (resp. $X$, $M$) compl\'ementaires de ces points. 

Il est clair que $X_{\Top}$ est r\'eunion de l'ouvert $X^*_{\Top}$ (qui s'identifie \`a  $M^*_{\Top}\cong F^\times$ par $q_{\Top}$) et de l'image $I$ de $\pi_{\Top}$. 
Comme $\pi_{\Top}$ est ouverte (\ref{prop:topLisseEtale}), $I$ est aussi un ouvert de $X_{\Top}$, et $q$ induit une bijection de $I$ sur l'ensemble $Q$ des carr\'es de $F$. Cette bijection est m\^eme un hom\'eomorphisme (remarquer que $f_{\Top}$ est ouverte sur son image) de sorte que $X_{\Top}$ peut se d\'ecrire comme la somme topologique $X_{\Top}=Q\amalg_{Q^\times} F^\times$, o\`u $Q^\times\subset F^\times$ est l'ensemble des carr\'es non nuls. Vu les hypoth\`eses faites sur $F$, $Q^\times$ est ouvert et ferm\'e dans $F^\times$,
 qui est donc somme disjointe de $Q^\times$ et de $F\setminus Q$: en conclusion,
$$X_{\Top}=Q\amalg (F\setminus Q).$$

\begin{fssrem}\label{rem:Kollar}
 Le morphisme $q:X\to M$ est un exemple de \og bug-eyed cover\fg\ au sens de \cite[\S4]{Ko}. 
 
 L'espace $X$ appara\^{\i}t aussi dans \cite{Sch} comme exemple d'un espace al\-g\'e\-brique admettant un $\GG_{m}$-torseur non localement trivial pour
 la topologie de Zariski; ce dernier s'obtient, par le changement de groupe $\mu_{2}\inj\GG_{m,F}$, \`a partir du $\mu_{2}$-torseur que nous allons  
d\'ecrire plus bas (\ref{ssec:mu2-torseur}).
\end{fssrem}
\subsubsection{Le cas o\`u $F=\RR$. } Dans ce cas, la description ci-dessus montre que $X_{\Top}$  s'identifie (via $q:X_{\Top}\to\RR$) \`a la somme disjointe
 des intervalles $\RR_{<0}$ et $\RR_{\geq0}$.
L'application $\pi_{\Top}$ s'identifie \`a l'application $t\mapsto t^2$ de $\RR$ dans $\RR_{\geq0}$; ici, elle  admet en tout point de $L_{\Top}$ 
des sections locales continues.

\subsubsection{Le cas o\`u $F=\CC$. } Alors $q_{\Top}: X_{\Top}\to \CC$ est un hom\'eomorphisme, et $\pi_{\Top}$ n'admet pas de sections 
locales continues au voisinage de $x_{0}$.
\medskip

On remarquera que l'injection naturelle $X(\RR)_{\Top}\inj X(\CC)_{\Top}$ n'est pas un plongement topologique.

\subsubsection{Un $\mu_{2}$-torseur sur $X$. }\label{ssec:mu2-torseur}Nous allons construire un $\mu_{2}$-torseur $\alpha:\til{X}\to X$ tel que 
l'application induite $\til{X}_{\Top}\to \im(\alpha_{\Top})$ ne soit pas un $\{\pm1\}$-fibr\'e principal. Consid\'erons deux exemplaires $L_{1}$ et
$L_{2}$ de $L$ (avec origines $l_{1}$ et $l_{2}$) et recollons-les en identifiant les ouverts  $L_{1}^*$ et $L_{2}^*$  de la mani\`ere \'evidente. 
On obtient une $F$-vari\'et\'e non s\'epar\'ee $\til{X}$, qui est la \og droite \`a point d\'edoubl\'e\fg\ bien connue. Il existe une unique involution
 $\sigma$ de $\til{X}$ qui induit sur $\til{X}^*\cong \GG_{m}$ l'application $x\mapsto -x$, et qui \emph{\'echange} les images de $l_{1}$ et $l_{2}$.
 On en d\'eduit une action libre de $\mu_{2}$ sur $\til{X}$, et il est facile de v\'erifier que l'application \'evidente $\alpha:\til{X}\to X$ h\'erit\'ee
 des projections $L_{i}\cong L\xrightarrow{\pi}X$ ($i=1,\,2$) identifie $X$ au quotient de $\til{X}$ par cette action et fait de $\til{X}$ un $\mu_{2}$-torseur sur $X$. L'image $I$ de $\alpha_{\Top}$ est celle de $\pi_{\Top}$, et $\til{X}_{\Top}$ n'est pas un $\{\pm1\}$-fibr\'e principal 
sur $I$ car $I$ est s\'epar\'e et $\til{X}_{\Top}$ ne l'est pas.

Si $F=\RR$, $I$ s'identifie \`a  $\RR_{\geq0}$ et $\til{X}_{\Top}$ \`a la r\'eunion de deux exemplaires $L_{1,\Top}$ et $L_{2,\Top}$  de $\RR$,
 identifi\'es le long de $\RR^\times$; 
l'application $\pi_{\Top}$ envoie  $L_{i,\Top}$ ($i=1,2$) sur $I$ par l'application $t\mapsto t^2$. 
On voit donc que $\pi_{\Top}$ a \emph{quatre} sections continues $s_{i,\varepsilon}$ ($i=1,2; \varepsilon=\pm1$), o\`u $s_{i,\varepsilon}$ 
 envoie le point d'abscisse $u$ de $I$ sur 
le point d'abscisse $\varepsilon\sqrt{u}$ de $L_{i,\Top}$. 

Si $F=\CC$, $\pi_{\Top}$ n'a pas de section continue.

\begin{fssrem}
Le $\GG_{m}$-torseur de Schr\"oer \cite{Sch}, d\'ej\`a mentionn\'e en \ref{rem:Kollar}, s'obtient comme quotient de $L\times_{F}\GG_{m,F}$ par la 
relation d'\'equivalence qui identifie $(x,\lambda)\in L^*\times_{F}\GG_{m,F}$ \`a $(-x,-\lambda)$. Il donne un autre exemple o\`u la conclusion de \ref{ThPpal}\,\ref{ThPpal3} est en d\'efaut; cette fois le groupe $G$ est lisse et connexe.
\end{fssrem}

\subsubsection{Encore un torseur sur $X$. }Pour le lecteur qui jugerait artificiels les exemples pr\'ec\'edents, nous allons d\'ebusquer le m\^eme 
espace $X$ \og dans la nature\fg, comme quotient d'une honn\^ete vari\'et\'e (quasi-affine) par une action de groupe. Soit $U\subset\Aa^2_{F}=\Spec F[x,y]$ l'ouvert compl\'ementaire de l'origine. Soit $T$ le tore maximal de $\SL_{2,F}$ form\'e des matrices $\begin{pmatrix}
\lambda&0\\0&\lambda^{-1}\end{pmatrix}$ ($\lambda\neq0$), et soit $G=T\cup\rho T$ o\`u $\rho=\begin{pmatrix}
0&-1\\1&0\end{pmatrix}$. En d'autres termes, $G$ est le sous-groupe de $\SL_{2,F}$ stabilisant la paire $\{h,-h\}$ o\`u $h$ est la forme 
quadratique $(x,y)\mapsto xy$, et $T=G^\circ=\mathrm{SO}(h)$. 

Il est imm\'ediat que $G$ op\`ere librement sur $U$, mais qu'il a deux types d'orbites:
\begin{itemize}
\item les r\'eunions de deux hyperboles \og oppos\'ees\fg, d'\'equation de la forme $(xy-a)(xy+a)=0$ (avec $a\neq0$);
\item la r\'eunion des deux axes, priv\'ee de l'origine.
\end{itemize}
Consid\'erons d'abord le quotient de $U$ par $T$. Soient $U_{1}$ et $U_{2}$ les ouverts de $U$ compl\'ementaires des deux axes: alors la fonction $xy$ induit un isomorphisme de $T\backslash U_{1}$ (resp. $T\backslash U_{2}$) sur un exemplaire $L_{1}$ (resp. $L_{2}$) de la droite affine, et il est facile d'en d\'eduire un isomorphisme $T\backslash U\flis \til{X}$, o\`u $\til{X}$ est d\'efini en \ref{ssec:mu2-torseur}. Noter que $U\to \til{X}$ est localement trivial pour la topologie de Zariski et qu'en  particulier l'application induite $U_{\Top}\to \til{X}_{\Top}$ est un $T_{\Top}$-fibr\'e principal.

La matrice $\rho$ induit sur $\til{X}$ l'involution $\sigma$ de \ref{ssec:mu2-torseur}:
elle transforme l'hyperbole $xy=a$ en $xy=-a$, et elle \'echange les deux axes. Par suite le quotient $G\backslash U$ 
n'est autre que $X$, et $U$ est un $G$-torseur sur $X$. 
Ici encore, $U_{\Top}$ n'est pas un $G_{\Top}$-fibr\'e principal sur son image, puisqu'il se factorise par $\til{X}_{\Top}\to I$ qui n'en est pas un, et que $U_{\Top}\to \til{X}_{\Top}$ est surjectif.

\Subsection{Cas d'un sch\'ema en groupes non constant}

Soient $K$ et $Y$ comme dans le th\'eor\`eme \ref{ThPpal}. Il est naturel de se demander ce que l'on peut dire de $f_{\Top}:X_{\Top}\to Y_{\Top}$ lorsque $f:X\to Y$ est un torseur sous un $Y$-espace en groupes ${\mathfrak G}$, non n\'ecessairement \og constant\fg\ sur $Y$ (par abus, nous appellerons constant un groupe $Y\times_{K}G$ o\`u $G$ est un $K$-groupe). 

Un cas particulier important est celui d'un groupe ${\mathfrak G}$ \emph{localement constant} (au sens \'etale ou fppf). M\^eme dans ce cas, la strat\'egie de d\'emonstration utilis\'ee pour  \ref{ThPpal} se heurte imm\'ediatement \`a l'absence d'un analogue du sous-groupe $G^\bec$.

Dans la suite, nous supposerons en g\'en\'eral que la base $Y$ est localement s\'e\-pa\-r\'ee; ceci assure que $\mathfrak{G}_{\Top}\to Y_{\Top}$ est un \og groupe topologique relatif\fg\ op\'erant sur $X_{\Top}$  (\ref{prop:ProdFibEspAlg}\,\ref{prop:ProdFibEspAlg24}).

\begin{fsprop}\label{prop_lisse_non_constant}
Soient $F$ un corps  topologiquement hens\'elien,  $Y$ un $F$-espace localement s\'epar\'e, ${\mathfrak G}$ un $Y$-espace alg\'ebrique en groupes lisse de type fini \modif{et quasi-s\'epar\'e},  $f:X \to Y$ un  ${\mathfrak G}$-torseur. 
Alors $f_\Top$ admet des sections locales en tout point de $X_{\Top}$; en particulier elle est ouverte 
 et l'application induite $X_{\Top}\to\im{f_{\Top}}$ est une $\mathfrak{G}_{\Top}$-fibration principale. 
\end{fsprop}

\dem le morphisme $f$ est lisse et $Y$ est localement s\'epar\'e donc $f_\Top$ admet des sections locales  (\ref{prop:topLisseEtale}\,\ref{prop:topLisseEtale3}). \modif{(L'hypoth\`ese de quasi-s\'eparation est l\`a pour assurer que $\FG$ est un $F$-espace au sens du pr\'esent article)}.\cqfd

\begin{fsrem}
Dans la situation de \ref{prop_lisse_non_constant}, l'image de $f_{\Top}$ n'est pas n\'ecessairement ferm\'ee, comme le montre l'exemple suivant. Supposant $F$ de caract\'eristique diff\'erente de $2$, prenons pour $Y$ la droite affine $\Spec{F[t]}$. Consid\'erons la $\cO_{Y}$-alg\`ebre $\cA:=\cO_{Y}[Z]/(Z^2-t^2)$, qui est libre de rang $2$. Soit $\FG$ le groupe des unit\'es de norme $1$ de $\cA$: on peut le voir comme le sous-groupe de $\GL_{2,\cO_{Y}}$ form\'e des matrices de la forme 
$M(x,y)={\begin{pmatrix}
x& t^2y\\y&x
\end{pmatrix}}
$ et de d\'eterminant $1$. C'est un $Y$-groupe affine, lisse et commutatif dont la  restriction  \`a  $U:=\Spec{F[t,t^{-1}]}$ est isomorphe \`a $\GG_{m,U}$ (par $M(x,y)\mapsto x+ty$) et dont la fibre \`a l'origine est isomorphe \`a $\mu_{2,F}\times\GG_{a,F}$ (par $M(x,y)\mapsto (x,xy)$).

Soit $d\in F^\times\smallsetminus F^{\times2}$. Alors les matrices $M(x,y)$ de d\'eterminant $d$ forment un $\FG$-torseur $X\to Y$. Il est trivial sur $U$ comme tout $\GG_{m}$-torseur; une section explicite est donn\'ee par $M\left(\frac{1+d}{2},\frac{1-d}{2t}\right)$. Sa fibre \`a l'origine est isomorphe \`a $\Spec{F\bigl(\sqrt{d}\bigr)}\times \Aa^1_{F}$ donc n'a pas de point rationnel. Ainsi l'image de $f_{\Top}:X\to F$ est $F^\times$.
\end{fsrem}

\begin{fslem}\label{lemme:locconstant} Soient $k$ un corps, $Y$ un $F$-espace, $\FG$ un \modif{$Y$-espace al\-g\'e\-brique en groupes quasi-s\'epar\'e de type fini}. On suppose que  ${\mathfrak G}$ est localement isomorphe pour la topologie \'etale \`a un $Y$-groupe constant $G \times_k Y$; il correspond donc \`a une classe $\alpha\in\rH^1(Y,\ul{\Aut}\,G)$. 

On suppose en outre qu'il existe un $k$-sch\'ema en groupes \emph{lisse} $\Gamma$ et un $k$-morphisme $\varphi:\Gamma\to \ul{\Aut}\,G$ tels que $\alpha$ soit l'image d'une classe $\gamma\in\rH^1(Y,\Gamma)$ par $\rH^1(Y,\varphi):\rH^1(Y,\Gamma)\to\rH^1(Y,\ul{\Aut}\,G)$.

Alors, pour tout $y\in Y(k)$, il existe un voisinage \'etale \emph{point\'e} $(Y',y')$ de $y$ tel que le $Y'$-groupe $\FG\times_{Y}Y'$ soit constant \emph{(et donc isomorphe \`a $Y'\times_{k}\FG_{y}$ o\`u $\FG_{y}$ est la fibre de $\FG$ en $y$)}.
\end{fslem}
\dem la classe $\gamma$ correspond \`a un $\Gamma$-torseur $T\to Y$. Soit $T_{1}\to Y$ le $\Gamma$-torseur $Y\times_{k}T_{y}$ o\`u $T_{y}$ est la fibre de $T$ en $y$. Alors $U:=\ul{\mathrm{Isom}}_{\Gamma}(T,T_{1})$ est un torseur sous $\ul{\Aut}_{\Gamma}(T)$ qui est une forme de $\Gamma$ et donc un \modif{$Y$-espace} en groupes lisse. Donc $U$ est un $Y$-espace alg\'ebrique lisse sur $Y$. Comme il a  un point rationnel au-dessus de $y$, il admet donc une section sur un voisinage \'etale point\'e $(Y',y')$ de $y$. Donc $T$ devient constant sur $Y'$, c'est-\`a-dire que $\gamma_{Y'}$ provient d'une classe dans $\rH^1(k,\Gamma)$, d'o\`u il suit que $\alpha_{Y'}$ provient d'une classe de $\rH^1(k,\ul{\Aut}\,G)$, d'o\`u la conclusion.\cqfd
\medskip

\begin{fscor}\label{cor:locconstant} Sous les hypoth\`eses du lemme \rref{lemme:locconstant}, on suppose en outre que $Y$ est localement s\'epar\'e, que $k$ est un corps  topologiquement hen\-s\'e\-lien, et que:
\begin{itemize}
\item ou bien $\FG$ est lisse sur $Y$ \emph{(ou que $G$ est un $k$-groupe lisse, ce qui revient au m\^eme si $Y\neq\emptyset$)};
\item ou bien $k$ est un corps valu\'e admissible.
\end{itemize}
Soit $f:X\to Y$ un $\FG$-torseur. Alors $\im\,(f_{\Top})\subset Y_{\Top}$ est localement ferm\'ee, et l'application induite $X_{\Top}\to\im\,(f_{\Top})$ est une $\FG_{\Top}$-fibration principale. 

De plus $\im\,(f_{\Top})$ est ouverte et  ferm\'ee si $\FG$ est lisse, et est ferm\'ee si $G$ v\'erifie \textup{(*)}.
\end{fscor}
\dem la question \'etant locale sur $Y_{\Top}$, le lemme \ref{lemme:locconstant} ram\`ene la situation au cas d'un $Y$-groupe constant (l'hypoth\`ese sur $Y$ assure qu'avec les notations du lemme, $Y'_{\Top}\to Y_{\Top}$ est un hom\'eomorphisme local). On applique alors le th\'eor\`eme \ref{ThPpal} dans le cas admissible, et la proposition \ref{prop:tors-groupe-lisse} dans le cas lisse.\cqfd

\begin{fsrems}
La condition du lemme \ref{lemme:locconstant} sur l'existence de $\varphi:\Gamma\to  \ul{\Aut}\,G$ est tr\`es restrictive. Elle est trivialement v\'erifi\'ee si $\ul{\Aut}\,G$  est re\-pr\'e\-sen\-table et lisse, notamment lorsque $\FG$ est \emph{r\'eductif} \cite[XXIV, th\'eo\-r\`eme 1.3]{SGA3}; noter que dans ce cas  $\ul{\Aut}\,G$  n'est pas toujours de type fini.\smallskip

Un autre cas utile est celui des \og $Y$-formes fortement int\'erieures\fg\ (voir par exemple \cite[2.2.4.9]{C-F}; certains auteurs les appellent \og formes  int\'e\-rieures pures\fg) d'un $k$-groupe al\-g\'e\-brique $G$: partant d'un $G$-torseur (\`a droite) $X_{1}\to Y$, on consid\`ere le $Y$-groupe $\FG:=\ul{\Aut}_{G}(X_{1})$. 
Sa classe dans $\rH^1(Y,\ul{\Aut}\,G)$ est l'image de la classe de $X_{1}$ par $\rH^1(Y,G)\xrightarrow{\,\rH^1(\mathrm{int})\,}\rH^1(Y,\ul{\Aut}\,G)$ o\`u $\mathrm{int}:G\to\ul{\Aut}\,G$ est le morphisme de conjugaison. 
Ce dernier passe au quo\-tient par le centre $Z(G)$ de $G$, 
de sorte que le lemme  \ref{lemme:locconstant} et le corollaire  \ref{cor:locconstant} s'ap\-pliquent \`a $\FG$ chaque fois que $G/Z(G)$ est lisse. Comme tout $\FG$-torseur \`a droite est de la forme $\ul{\mathrm{Isom}}_{G}(X_{1},X_{2})$ o\`u $X_{2}$ est un autre $G$-torseur sur $Y$, on obtient la proposition \ref{prop:Isom} qui suit.
\end{fsrems}

\begin{fsprop}\label{prop:Isom} Soient $(K,v)$ un corps valu\'e  admissible et $Y$ un $K$-espace localement s\'epar\'e.
Soit $G$ un $K$-groupe alg\'ebrique. On  suppose  que  $G/Z(G)$ est lisse. 

 Soient $X_1 \to Y$ et $X_2 \to Y$ des $G$-torseurs. On note ${\mathfrak G}_i:=\ul\Aut_{G}(X_{i})\to Y$
 le tordu int\'erieur de $G$ par le torseur $X_i$ $(i=1,2)$. On consid\`ere le $({\mathfrak G}_1, {\mathfrak G}_2)$-bitorseur  $f:T \to Y$ d\'efini par
$$ 
T:= \underline {\rm Isom}_G(X_1,X_2) 
$$ 
On d\'esigne par $I$ l'image de $f_\Top: T_\Top \to Y_\Top$. 
\begin{numlist}

\item $I$ est localement ferm\'e dans $Y_\Top$, et est ferm\'e si $G$ v\'erifie \textup{(*)};

\item pour tout $y \in I$,   il existe un voisinage ouvert $U_y$ de $y$ dans $I$ de sorte que l'application $f_\Top^{-1}(U_y) \to U_y$ soit 
munie d'une structure de  $(G_y)_\Top$-fibr\'e principal, o\`u  $G_y= {\mathfrak G}_2 \times_Y y$.\cqfd
\end{numlist}
\end{fsprop}

\addvspace{\bigskipamount}

\end{document}